\documentclass[stsy,nonblindrev]{informs-stsy}

\OneAndAHalfSpacedXI

\usepackage{latexsym}
\usepackage{epsfig}
\usepackage{subfig}
\usepackage{amsmath}
\usepackage{amssymb}
\usepackage{graphicx}
\usepackage{color}
\usepackage{multirow}
\usepackage{endnotes}
\usepackage{enumerate}
\usepackage{bm}
\usepackage{booktabs}
\usepackage[compact]{titlesec}
\usepackage{lscape}
\usepackage{natbib}
\usepackage[normalem]{ulem}
\usepackage{enumitem}
\usepackage{caption}

\usepackage{natbib}
 \bibpunct[, ]{(}{)}{,}{a}{}{,}%
 \def\bibfont{\small}%
\usepackage[colorlinks=true,breaklinks=true,bookmarks=true,urlcolor=blue,
     citecolor=blue,linkcolor=blue,bookmarksopen=false,draft=false,hypertexnames=false]{hyperref}

\def\EMAIL#1{\href{mailto:#1}{#1}}

\newcommand{\ot}[1]{\tilde{#1}}

\newcommand{\ru}[1]{\big\lceil #1 \big\rceil}
\newcommand{\rd}[1]{\lfloor #1 \rfloor}
\newcommand{\ep}{\epsilon}
\newcommand{\al}{\alpha}
\newcommand{\ga}{\gamma}
\newcommand{\la}{\lambda}

\newcommand{\E}{\mathbf{E}}
\newcommand{\pr}{\mathbf{P}}

\newcommand{\dr}{\mathrm{d}}
\newcommand{\e}{\mathrm{e}}

\newcommand{\D}{\mathbb{D}}

\newcommand{\I}{\mathbb{I}}
\newcommand{\N}{\mathbb{N}}
\newcommand{\R}{\mathbb{R}}

\newcommand{\mA}{\mathcal{A}}
\newcommand{\mB}{\mathcal{B}}

\newcommand{\mE}{\mathcal{E}}
\newcommand{\mF}{\mathcal{F}}

\newcommand{\mI}{\mathcal{I}}
\newcommand{\mJ}{\mathcal{J}}

\newcommand{\mN}{\mathcal{N}}
\newcommand{\mP}{\mathcal{P}}
\newcommand{\mR}{\mathcal{R}}
\newcommand{\mS}{\mathcal{S}}
\newcommand{\mX}{\mathcal{X}}

\newcommand{\mZ}{\mathcal{Z}}

\TheoremsNumberedThrough     
\ECRepeatTheorems

\EquationsNumberedThrough    

\MANUSCRIPTNO{}

\begin{document}



\RUNTITLE{Resource Collaboration and Multitasking}

\TITLE{On the Optimal Control of Parallel Processing Networks with Resource Collaboration and Multitasking}

\ARTICLEAUTHORS{%
\AUTHOR{Erhun \"Ozkan}
\AFF{College of Administrative Sciences and Economics, Ko\c c University, Istanbul, Turkey, \EMAIL{erhozkan@ku.edu.tr}}
} 

\ABSTRACT{%
We study scheduling control of parallel processing networks in which some resources need to simultaneously collaborate to perform some activities and some resources multitask. Resource collaboration and multitasking give rise to synchronization constraints in resource scheduling when the resources are not divisible, that is, when the resources cannot be split. The synchronization constraints affect the system performance significantly. For example, because of those constraints, the system capacity can be strictly less than the capacity of the bottleneck resource. Furthermore, the resource scheduling decisions are not trivial under those constraints. For example, not all static prioritization policies retain the maximum system capacity and the ones that retain the maximum system capacity do not necessarily minimize the delay (or in general the holding cost). We study optimal scheduling control of a class of parallel networks and propose a dynamic prioritization policy that retains the maximum system capacity and is asymptotically optimal in diffusion scale and conventional heavy-traffic regime with respect to the expected discounted total holding cost objective.
}%



\maketitle

%


\section{Introduction}\label{s_intro}

We study control of processing networks in which some resources need to collaborate to perform some activities and some resources multitask. Resource collaboration implies the cases in which multiple resources need to simultaneously collaborate in order to perform an activity. For example, a doctor, a nurse, and an anesthetist may need to simultaneously collaborate in order to perform a surgery and a worker may need to use a specific machine to perform a manufacturing operation. Multitasking means that some resources perform multiple activities (resource sharing) and is prevalent in many application domains. For example, a doctor makes both patient diagnosis and discharge decisions; a research \& development team can work on different types of projects; a commodity cluster can be used by multiple diverse cluster computing frameworks in a data network. 

\begin{figure}[ht]
  \centering
  \subfloat[A parallel network with three job types and three multitasking resources (R1, R2, R3).]{\label{f_bc2}\includegraphics[width=0.42\textwidth]{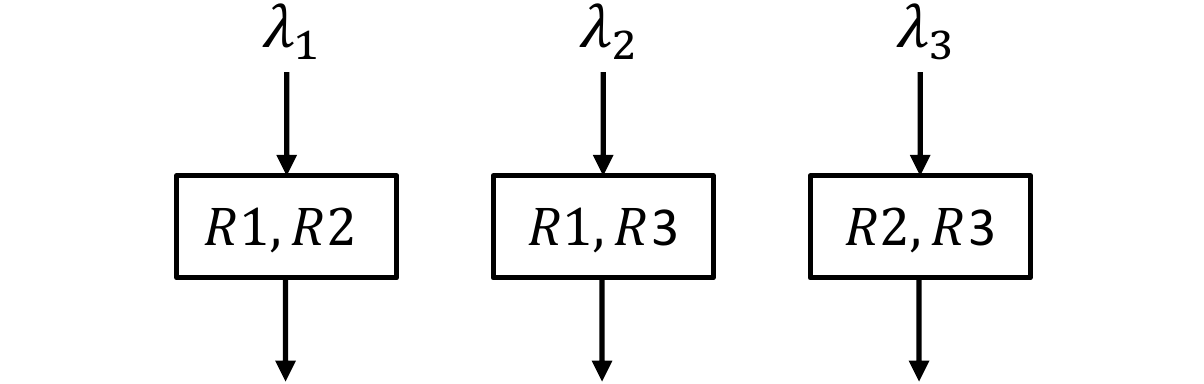}}\hspace{1cm}
  \subfloat[A parallel network with three job types and two multitasking resources (R1, R2).]{\label{f_bc}\includegraphics[width=0.42\textwidth]{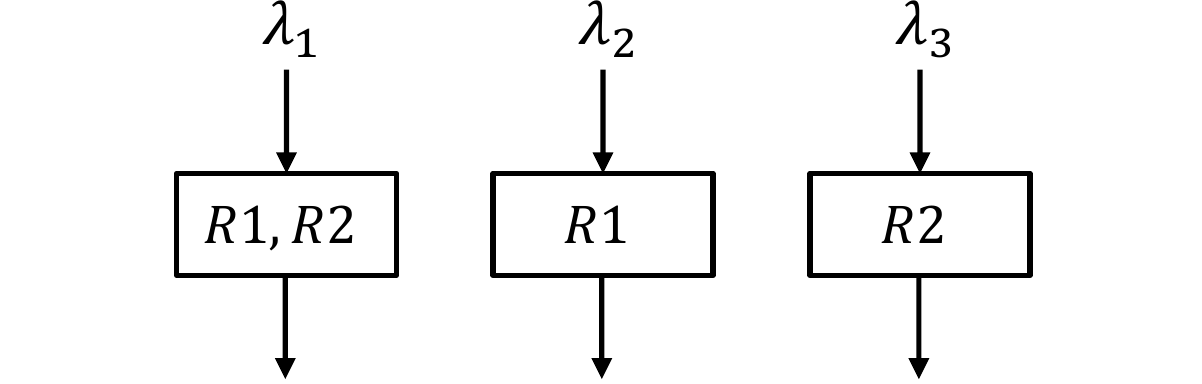}}\\
  \subfloat[The graph associated with Figure \ref{f_bc2}.]{\label{f_bc2g}\includegraphics[width=0.42\textwidth]{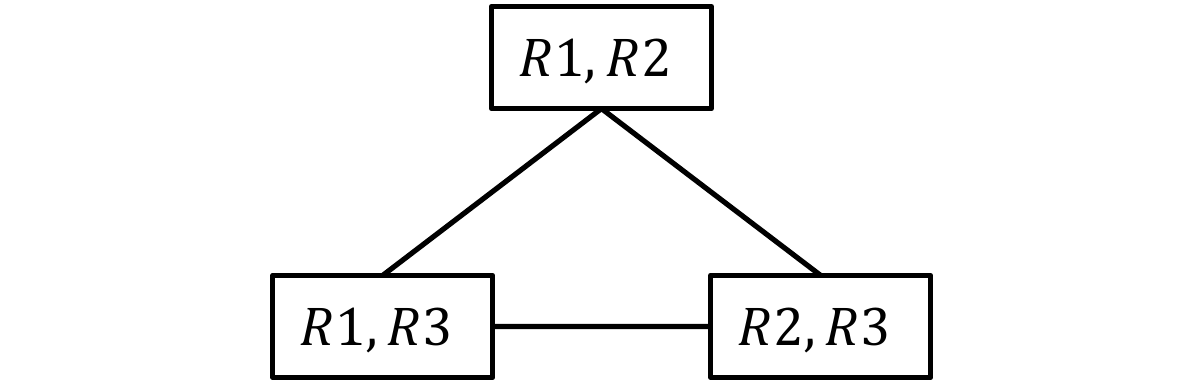}}\hspace{1cm}
  \subfloat[The graph associated with Figure \ref{f_bc}.]{\label{f_bcg}\includegraphics[width=0.42\textwidth]{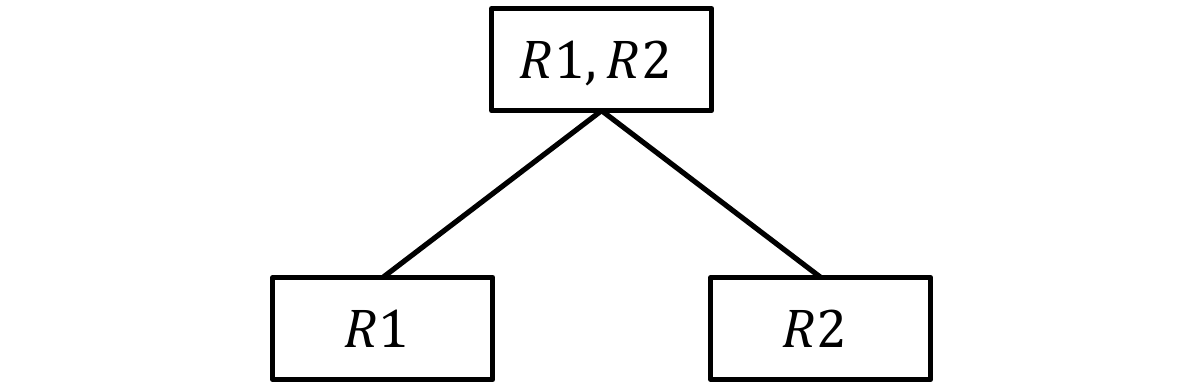}}
  \caption{Two examples of parallel networks with resource collaboration and multitasking.}
  \label{f_network}
\end{figure}

Resource collaboration and multitasking give rise to synchronization constraints in resource scheduling when the resources are not divisible, that is, when the resources cannot be split. For example, a doctor cannot make diagnosis and discharge decisions simultaneously. The synchronization constraints affect the system performance significantly. According to the traditional capacity analysis, system capacity is equal to the capacity of the bottleneck resource. However, in the presence of collaboration and multitasking, the system capacity can be strictly less than the bottleneck capacity, that is, there can be a capacity loss. To illustrate this, let us consider the simple parallel network in Figure \ref{f_bc2}, in which there are three job types arriving in the system and three multitasking resources. Type 1 jobs are processed by simultaneous collaboration of resources 1 and 2 (R1 and R2). Type 2 (3) jobs are processed by simultaneous collaboration of R1 (R2) and R3. Observe that whenever a job is processed, two resources are utilized and one resource is forced to idle. Therefore, the resources are not fully utilized and so there is a capacity loss.

Next, let us consider the parallel network in Figure \ref{f_bc}, in which there are only two resources (R1 and R2). Consider the static prioritization policy under which type 1 jobs have preemptive priority over the type 2 and 3 jobs. Under the aforementioned policy, both resources are utilized in a work-conserving fashion, that is, a resource does not idle as long as there is a job in the system that it can process alone or in collaboration with other resources. Consequently, both resources can be fully utilized and so there is no capacity loss in the network in Figure \ref{f_bc}.

The main difference between the networks in Figures \ref{f_bc2} and \ref{f_bc} is their underlying collaboration architecture which specifies the required resources in the service processes associated with each job type. \cite{gur15} derive conditions related to the collaboration architecture that are necessary and sufficient for the system capacity to be equal to the bottleneck capacity. For example, they show that there is no capacity loss in networks with \textit{hierarchical collaboration architecture}, in which the network can be represented by a graph where nodes are the job types, there is an edge between two nodes (job types) if there exists a resource required in the service of both job types, nodes are arranged in collaboration levels, and resources at each node are a subset of those used by the node above it. As can be seen in Figures \ref{f_bc2g} and \ref{f_bcg}, the network in Figure \ref{f_bc} has hierarchical collaboration architecture but the one in Figure \ref{f_bc2} has not. In networks with hierarchical collaboration architecture, there exist scheduling policies under which all resources are utilized in a work-conserving fashion and thus there is no capacity loss under those policies. For example, consider the \textit{collaborative hierarchical preemptive priority policy} proposed by \cite{gur18} (see Section 5 therein and we will call this policy \textit{the GVM policy} henceforth) under which the jobs requiring the highest number of resources in their service processes receive preemptive priority over the remaining jobs.

Although there is no capacity loss under the GVM policy, it does not necessarily minimize the delay (or in general holding cost). To illustrate this, let us consider the objective of minimizing expected discounted total holding cost in the network in Figure \ref{f_bc}. Let $h_j$ and $\mu_j$ denote the holding cost rate per type $j$ job per unit time and the service rate for type $j$ jobs, respectively, for all $j\in\{1,2,3\}$. On the one hand, if $h_1\mu_1\geq h_2\mu_2 + h_3\mu_3$, that is, if holding a type 1 job is more expensive than jointly holding a type 2 job and a type 3 job in the system, then the GVM policy makes the right scheduling decisions. On the other hand, if $h_1\mu_1< h_2\mu_2 + h_3\mu_3$, then, intuitively, under a good scheduling policy, the resources should give priority to the type 2 and 3 jobs as much as possible, which is not the case under the GVM policy (see Section \ref{s_ex_2} for a numerical experiment under which the GVM policy does not perform well). Hence, let us consider the \textit{priority to individual activities} (PIA) policy introduced by \cite{gur18}, under which the resources give preemptive priority to the job types that they process individually. Although the PIA policy is tempting especially when $h_1\mu_1 << h_2\mu_2 + h_3\mu_3$, it causes a capacity loss due to its non-work-conserving structure (see Section \ref{s_ex_2} for details). Consequently, there is a fine line between cost minimization and instability in networks with resource collaboration and multitasking.

Both the GVM and the PIA policies are static prioritization policies, neither of them necessarily minimize the holding cost, and the latter policy does not even retain the maximum system capacity. These results bring out the following question: \textit{Does there exist a scheduling policy which retains the maximum system capacity and minimizes the expected discounted total holding cost?} In this paper, we answer the aforementioned question affirmatively. Our main contribution is a dynamic prioritization policy for parallel networks with hierarchical collaboration architecture that retains the maximum system capacity and is asymptotically optimal in diffusion scale in the conventional heavy-traffic regime. We also show that our proposed policy is implementable and asymptotically optimal in a class of parallel networks with capacity loss (e.g., the network in Figure \ref{f_bc2}).

Deriving an exact optimal control policy is very challenging even for the simple network in Figure \ref{f_bc}. A potential approach is to use Markov Decision Process (MDP) techniques under the assumption that the inter-arrival and service times are exponentially distributed. However, because the associated state is the number of each job type in the system and the allocation of the resources to the job types (that is, a 4-dimensional state space for the simple network in Figure \ref{f_bc}), curse of dimensionality arises. Therefore, a more efficient solution approach is to derive asymptotically optimal control policies in the conventional heavy-traffic regime. We assume that there exists at least one resource in heavy traffic, that is, there exists at least one resource whose processing capacity is barely enough to process all incoming jobs requiring the service of that resource. Otherwise, if all resources are in light traffic, that is, if the processing capacities of all resources are more than enough, then any work-conserving policy will perform reasonably well. Our assumption that there exists at least one resource in heavy traffic is not restrictive in application domains where some of the resources are expensive. For example, in hospitals, resources such as doctors, nurses, CT-scanners generally have tight processing capacity.

We use Harrison's classical scheme (see \cite{har88,har00}): We formulate a control problem for the original stochastic processing network and define a notion of heavy traffic. Next, we use diffusion approximations to formulate a Brownian control problem (BCP). Then, we solve the BCP and interpret a control policy from the optimal BCP solution. Finally, we prove the asymptotic optimality of the proposed policy (see Section 1 of \cite{pes16} for more information about Harrison's scheme). A major technical challenge for our paper is that the resulting BCP is multidimensional. Specifically, the dimension of the resulting BCP is equal to the number of resources in heavy traffic and thus the resulting workload process is multidimensional. Because we consider a parallel network, in which there is no network effect, there is no control on the amount of workload arriving at each resource. Furthermore, because we focus on control policies under which the resources in heavy traffic are utilized in work-conserving fashion, we are able to prove that under any sequence of admissible policies, the multidimensional workload process weakly converges to a stochastic process with invariant distribution (that is, the distribution of the limiting stochastic process is independent of the underlying sequence of admissible policies). By utilizing the aforementioned convergence result, we are able to obtain a path-wise numerical solution for the multidimensional BCP. Specifically, the BCP is time-decomposable, and we can convert it to a linear program (LP) and solve it numerically at discrete time epochs. By utilizing the Lipschitz continuity of the optimal LP solutions, we prove that the optimal LP objective function values over time provide an asymptotic lower bound on the performance of any admissible policy (see Theorem \ref{t_lb}).

The proposed policy is a dynamic and state-dependent prioritization policy under which the LP is solved at discrete time epochs. The parameters of the LP are holding cost rates and service rates of the job types that are processed by at least one resource in heavy traffic and the weighted total number of jobs waiting to be processed by each resource in heavy traffic (that is, the ``workload'' for each resource in heavy traffic). The objective of the LP is to minimize the instantaneous holding cost rate by allocating the workloads for the resources in heavy traffic to the associated job types. In other words, the LP determines the optimal workload splits for the resources in heavy traffic. After the LP is solved, the system controller compares the numbers of the job types in the system with the optimal LP solution. If those numbers are different, the proposed policy allocates the resources to the job types based on an index rule for a ``fixed'' amount of time. The length of that ``fixed'' time interval is chosen such that at the end of the interval, the actual numbers of the job types in the system will be close to the most recently computed optimal LP solution with a high probability. Then, the system controller resolves the LP and follows the same procedure. Because the time between successive LP solution epochs is sufficiently short and the optimal LP solutions are Lipschitz continuous, the optimal LP solutions do not change significantly between successive LP solution epochs. Therefore, under the proposed policy, the number of each job type that is processed by at least one resource in heavy traffic always tracks an optimal LP solution sufficiently closely. Furthermore, because the resources in light traffic have excess capacity, the number of each job type that is processed by resources only in light traffic never grows a lot. Consequently, we prove the asymptotic optimality of the proposed policy (see Theorem \ref{t_ao}).

We present a literature review in Section \ref{s_lit} and some notation in Section \ref{s_not}. We present the model description, asymptotic framework, and the objective in Section \ref{s_model}. We present the proposed policy in Section \ref{s_policy}.  Next, we present the asymptotic optimality result in Section \ref{s_th}. We present some examples in Section \ref{s_ex}. Then, we relax a modeling assumption in Section \ref{s_ext}. Finally, we make concluding remarks in Section \ref{s_conc}. All the proofs are presented in the electronic companion.

\subsection{Literature Review}\label{s_lit}

There is a vast literature about resource sharing (multitasking) in networks. A significant amount of those studies are about communication networks. For example, \cite{cza01} consider grid networks, \cite{yu11} consider cellular communication networks, \cite{hin11} consider data networks, and \cite{sha14} consider switched networks.

\cite{har14} consider the control of a very general processing network with resource collaboration and multitasking. They present an open problem of devising a dynamic resource allocation policy that achieves what they call hierarchical greedy ideal (HGI) performance in the heavy traffic limit. \cite{bud20} propose a rate-allocation policy under which the server allocation rates are determined by certain thresholds for the queue-length processes and prove that their policy achieves the corresponding HGI performance in the heavy traffic limit under some conditions. The rate-allocation policy of \cite{bud20} requires resources to be divisible, that is, a resource must be able to process multiple job types simultaneously. For example, that policy requires the resources to be divisible even in the simple network in Figure \ref{f_bc} (see Section 5 of \cite{bud20} for details and observe that the network in Figure 2 of \cite{bud20} is equivalent to the network in Figure \ref{f_bc}). In contrast, our proposed policy does not require divisible resources. 


There are also studies about networks with resource collaboration and multitasking in which resources are indivisible. An early example is \cite{cou87}, where reward maximization in a specific parallel network is considered. They prove optimality of a threshold policy by using MDP techniques. \cite{dai05} consider throughput maximization in processing networks (see specifically Section 7 therein for indivisible resources). \cite{dob12} consider a tandem queueing network with two distinct resources and three stages and derive control policies that maximize throughput by using MDP techniques. They are motivated by physician scheduling in emergency departments. \cite{gur15} study capacity loss in networks with resource collaboration and multitasking. They derive necessary and sufficient conditions for the system capacity to be equal to the bottleneck capacity. In a follow-up study, \cite{gur18} study impact of static prioritization policies on the capacities of the networks with parallel servers and hierarchical collaboration architecture. Our study complements the latter study by considering a cost minimization objective. A very recent study considering cost minimization in parallel networks with resource collaboration and multitasking is \cite{zyc19}. In their setting, the resources are identical and a service process can require the collaboration of multiple resources. They prove asymptotic optimality of a static policy and a state-dependent policy in many-server heavy-traffic regime in a parallel network in fluid scale. In contrast, the resources are distinct and static policies can perform poorly in our setting. Finally, similar to \cite{gur15}, \cite{bo19} study the system capacity in networks with resource collaboration and multitasking. They prove that not only the exact computation of the system capacity of those networks is a strongly NP-hard problem but also approximating the system capacity accurately is an NP-hard problem. Nevertheless, they prove that the system capacity can be efficiently computed in networks with a special collaboration architecture.

There are studies that provide empirical evidence about importance of coordination of collaborating resources in hospitals. For example, \cite{gur20} analyze and quantify how collaboration of resources affects total processing time in hospitals and show that a hospitalist spends on average $20\%$ of his/her time for the coordination of collaborative activities. 

Resource collaboration and multitasking also appear in project management literature, specifically, in resource-constrained project scheduling (see for example \cite{bru99}). The aforementioned literature considers a static problem in the sense that finite number of jobs (projects) are considered. In contrast, jobs arrive in the system dynamically over time in our case.

\subsection{Notation}\label{s_not}

The set of nonnegative and strictly positive integers are denoted by $\N$ and $\N_+$, respectively. For all $n\in\N_+$, $\R^n$ denotes the $n$-dimensional Euclidean space and $\R_+^n$ denotes the nonnegative orthant in $\R^n$. For any $x,y\in \R$, $x\vee y:= \max\{x,y\}$, $x\wedge y:= \min\{x,y\}$, and $(x)^+:=x\vee 0$. For any $x:=(x_1,x_2,\ldots,x_n)\in\R^n$ and $y:=(y_1,y_2,\ldots,y_n)\in\R^n$, we let $|x-y|_{\infty}:=\max_{i\in\{1,2,\ldots,n\}}|x_i-y_i|$. For any $x\in \R$, $\rd{x}$ ($\lceil x \rceil$) denotes the greatest (smallest) integer which is smaller (greater) than or equal to $x$ and $[x]$ denotes the closest integer to $x$. For any given set $\mX$, $|\mX|$ denotes the cardinality of $\mX$.

For all $n\in\N_+$, $\D^n$ denotes the set of functions $f:\R_+\rightarrow\R^n$ that are right continuous with left limits. We let $\bm{0},\bm{e}\in \D$ be such that $\bm{0}(t)=0$ and $\bm{e}(t)=t$ for all $t\in\R_+$. For $x,y\in \D$, $x\vee y$, $x\wedge y$, and $(x)^+$ are elements of $\D$ such that $(x\vee y)(t):=x(t)\vee y(t)$, $(x\wedge y)(t):= x(t)\wedge y(t)$, and $(x)^+(t):=(x(t))^+$ for all $t\in\R_+$. For any $x\in \D$, we define the mappings $\Psi, \Phi:\D\rightarrow\D$ such that for all $t\in\R_+$,
\begin{equation*}
\Psi(x)(t) := \sup_{0\leq s\leq t} (-x(s))^+,\quad\quad \Phi(x)(t) := x(t)+\Psi(x)(t),
\end{equation*}
where $\Phi$ is the one-dimensional reflection map (see Chapter 13.5 of \cite{whi02}). For $x\in \D$ and $t\in \R_+$, we let $\Vert x \Vert_t := \sup_{0\leq s \leq t} |x(s)|$. We consider $\D^n$ endowed with the usual Skorokhod $J_1$ topology (see Chapter 3 of \cite{bil99}). Let $\mB(\D^n)$ denote the Borel $\sigma$-algebra on $\D^n$ associated with the Skorokhod $J_1$ topology. For stochastic processes $\{W^r$, $r\in\N_+\}$ and $W$ whose sample paths are in $\D^n$ for some $n\in\N_+$, ``$W^r \Rightarrow W$'' means that the probability measures induced by $\{W^r$, $r\in\N_+\}$ on $(\D^n,\mB(\D^n))$ converge weakly to the one induced by $W$ on $(\D^n,\mB(\D^n))$ as $r\rightarrow \infty$. All the convergence results hold as $r\rightarrow\infty$.

Let $\mN=\{1,2,\ldots,n\}$ and $x_i\in\D$ for all $i\in \mN$. Then $(x_i,i\in \mN)$ denotes the process $(x_1,x_2,\ldots,x_n)$ in $\D^n$. We abbreviate the phrase ``uniformly on compact intervals'' by ``u.o.c.'' and ``almost surely'' by ``a.s.''. We let $\xrightarrow{a.s.}$ denote almost sure convergence and $\overset{d}{=}$ denote ``equal in distribution''. We repeatedly use the fact that convergence in the $J_1$ metric is equivalent to u.o.c. convergence when the limit process is continuous (see page 124 in \cite{bil99}). Let $\{x^r,r\in\N\}$ be a sequence in $\D$ and $x\in\D$. Then $x^r\rightarrow x$ u.o.c., if $\Vert x^r-x\Vert_{t}\rightarrow 0$ for all $t\in\R_+$. We let ``$\circ$'' denote the composition map and $\I$ denote the indicator function. We assume that all the random variables and stochastic processes are defined in the same complete probability space $(\Omega, \mF, \pr)$, $\E$ denotes expectation under $\pr$, and $\pr(A, B):=\pr(A\cap B)$.

\section{Model Description}\label{s_model}

We consider a parallel network with $J\in\N_+$ different job types and $I\in\N_+$ distinct resources. We let $\mJ:=\{1,2,\ldots,J\}$ denote the set of job types and $\mI:=\{1,2,\ldots,I\}$ denote the set of resources. Each job type has an associated buffer and each job is processed only once.  A job may require collaboration of multiple resources in its service process and a resource can process at most a single job at any given time. Let $A$ be an $I\times J$ dimensional resource-job incidence matrix such that $A_{ij}:=1$ if resource $i$ is required in the service process of type $j$ jobs and $A_{ij}:=0$ otherwise. We let $\mI_j:=\{i\in\mI:A_{ij}=1\}$ for all $j\in\mJ$ and $\mJ_i:=\{j\in\mJ:A_{ij}=1\}$ for all $i\in\mI$. Then $\mI_j$ denotes the set of resources required in the process of type $j$ jobs and $\mJ_i$ denotes the set of job types that are processed by resource $i$. We assume that $|\mI_j|\geq 1$ and $|\mJ_i|\geq 1$ for all $i\in\mI$ and $j\in\mJ$. If $|\mI_j|>1$ for some $j\in\mJ$, then type $j$ jobs require collaboration of multiple resources. If $|\mJ_i|>1$ for some $i\in\mI$, then resource $i$ multitasks. Without loss of generality, we assume that for all $i,k\in\mI$ such that $i\neq k$, both $\mJ_i\not\subset\mJ_k$ and $\mJ_k\not\subset\mJ_i$. In other words, there does not exist a resource pair such that one of the resources processes all of the job types that the other one processes. Otherwise, we can identify that resource pair as a single resource.
 
\subsection{Collaboration Architecture}\label{s_ca}

Collaboration architecture of a network is characterized by the resource-job incidence matrix $A$ and is independent of the arrival and service rates. \cite{gur15} prove that the system capacity is equal to the bottleneck capacity in open networks with hierarchical collaboration architecture (see Theorem 4 therein).

\begin{definition}\label{d_hca}
A network has hierarchical collaboration architecture if it satisfies the following condition: If $\mI_j\cap\mI_l\neq\emptyset$ for some $j,l\in\mJ$, then $\mI_j\subset\mI_l$ or $\mI_l\subset\mI_j$.
\end{definition}

\begin{figure}[htbp]
  \centering
  \subfloat{\label{h1}\includegraphics[width=0.7\textwidth]{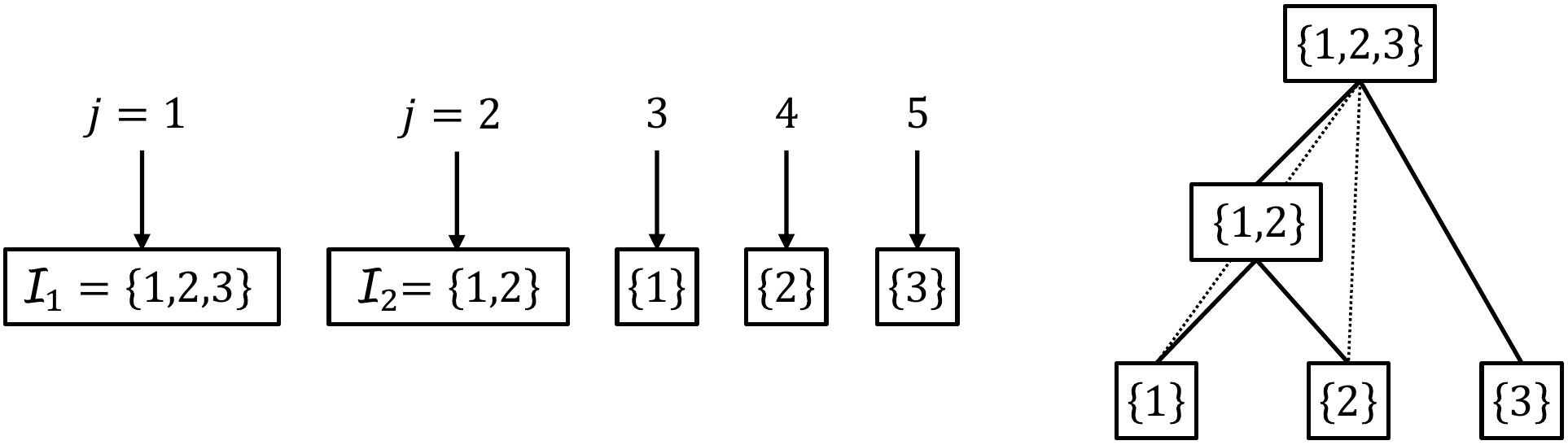}}\\
  \subfloat{\label{h2}\includegraphics[width=0.8\textwidth]{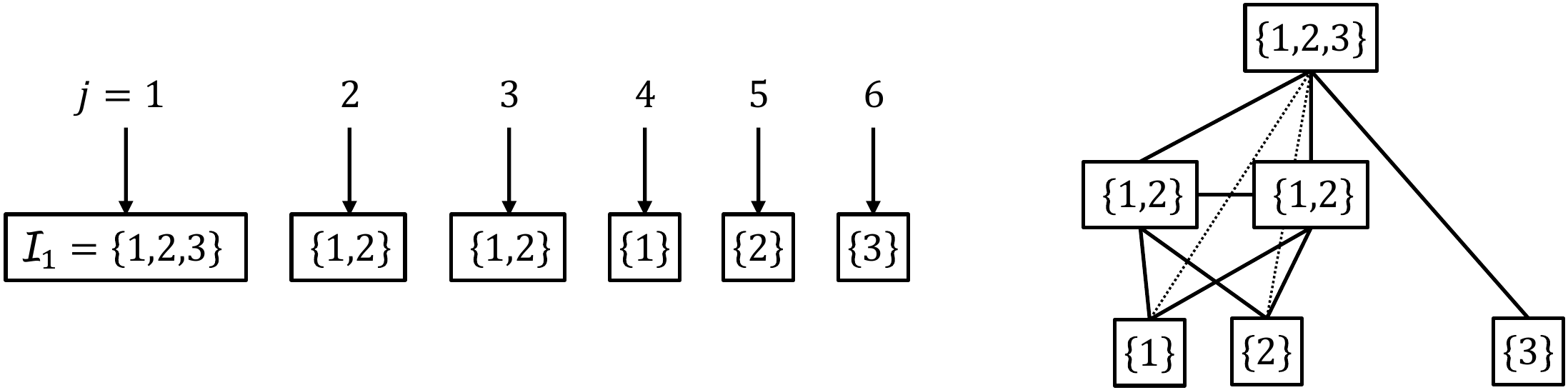}}\\
  \caption{Examples of parallel networks with hierarchical collaboration architecture.}
  \label{ex_hca}
\end{figure}

As \cite{gur18} state, ``hierarchical collaboration architectures can be represented by a graph where nodes are arranged in ``collaboration levels'' and resources at each node are a subset of those used by the node above it''. For example, Figure \ref{ex_hca} presents two parallel networks with hierarchical collaboration architecture together with the associated collaboration graphs.

\begin{assumption}\label{a_h} (Architecture) The parallel network that we consider has hierarchical collaboration architecture defined in Definition \ref{d_hca}.
\end{assumption}

Without loss of generality, we assume that the graph associated with the network architecture is connected. Otherwise, we can decompose the network and study each subnetwork separately.

\subsection{Stochastic Primitives and Asymptotic Framework}\label{s_primitives}

We consider a sequence of parallel networks indexed by $r\in\N_+$. We use the same sequence of service times but a different sequence of inter-arrival times in each network.

For all $j\in\mJ$, let $\{v_{jn},n\in\N_+\}$ be a strictly positive and independent and identically distributed (i.i.d.) sequence of random variables with mean $1/\mu_{j}$ and coefficient of variation $\sigma_{j}$. We let $v_{jn}$ denote the service time of the $n$th type $j$ job for all $j\in\mJ$ and $n\in\N_+$. For all $j\in\mJ$, $n\in\N_+$, and $t\in\R_+$, let $V_{j}(0):=0$ and
\begin{equation*}
V_{j}(n):= \sum_{k=1}^n v_{jk},\hspace{1cm} S_{j}(t):= \sup\left\{n\in\N: V_{j}(n)\leq t\right\}.
\end{equation*}
Then, $S_j$ is a renewal process such that $S_j(t)$ is the number of type $j$ jobs processed up to time $t$ given that type $j$ jobs are processed at all times in $[0,t]$.

In the $r$th system, we associate the inter-arrival times of type $j\in\mJ$ jobs with strictly positive and i.i.d. sequence of random variables $\{\bar{u}_{jn},n\in\N_+\}$ and the constant $\la_j^r>0$. For all $j\in\mJ$, $n\in\N_+$, and $r\in\N_+$, we have $\E[\bar{u}_{jn}]=1$, the variance of $\bar{u}_{jn}$ is $\beta_j^2$, and $u_{jn}^r:=\bar{u}_{jn}/\la_j^r$ denotes the inter-arrival time between the $(n-1)$st and $n$th type $j$ job in the $r$th system. Then, for all $j\in\mJ$ and $r\in\N_+$, $\{u_{jn}^r,n\in\N_+\}$ is an i.i.d. sequence of random variables with mean $1/\la_j^r$ and coefficient of variation $\beta_j$. For all $j\in\mJ$, $n\in\N_+$, $r\in\N_+$, and $t\in\R_+$, we let $U_j^r(0):=0$ and
\begin{equation*}
U_j^r(n):= \sum_{k=1}^n u_{jk}^r,\hspace{1cm} E_j^r(t):= \sup\left\{n\in\N: U_j^r(n)\leq t\right\}.
\end{equation*}
Then, $E_j^r$ is a renewal process such that $E_j^r(t)$ is the number of external type $j$ job arrivals up to time $t$ in the $r$th system.

For all $j\in\mJ$, we assume that the sequences $\{v_{jn},n\in\N_+\}$ and $\{\bar{u}_{jn},n\in\N_+\}$ are mutually independent of each other and all other stochastic primitives. Furthermore, we make the exponential moment assumption for the inter-arrival and service times.
\begin{assumption}\label{a_moment} (Moment)
There exists an $\bar{\al}>0$ such that for all $\al\in(-\bar{\al},\bar{\al})$,
\begin{equation*}
\E\left[ \e^{\al\bar{u}_{j1}}\right]<\infty,\qquad\E\left[ \e^{\al v_{j1}}\right]<\infty,\qquad\forall j\in\mJ.
\end{equation*}
\end{assumption}
Exponential moment assumption is common in the queueing literature, see for example \cite{har98,bel01,mag03,mey03,ozk19}.

Let $\rho_i^r:= \sum_{j\in\mJ_i} \la_j^r/ \mu_{j}$ for all $i\in\mI$, that is, $\rho_i^r$ denotes the load on resource $i$. The following assumption sets up the asymptotic regime.
\begin{assumption}\label{a_regime} (Asymptotic Regime)
\begin{enumerate}
\item $\la_j^r \rightarrow \la_j >0$  for all $j\in\mJ$.
\item $r\left( \rho_i^r-1\right) \rightarrow \theta_i\in\R\cup\{-\infty\}$ for all $i\in\mI$.
\end{enumerate}
\end{assumption}

Part 2 of Assumption \ref{a_regime} states that the resources can be in either light or heavy traffic. On the one hand, if $\theta_i=-\infty$ for some $i\in\mI$, then resource $i$ is in light traffic. On the other hand, if $\theta_i\in\R$ for some $i\in\mI$, then resource $i$ is in heavy traffic. Our definitions for resources in heavy traffic and light traffic are from \cite{gly90}, see Theorem 12 therein. We let $\mI^L:=\{i\in\mI: \theta_i=-\infty\}$ and $\mI^H:=\{i\in\mI: \theta_i\in\R\}$. Then, $\mI^L$ and $\mI^H$ denote the set of resources that are in light and heavy traffic, respectively. We let $\mJ^H:=\left\{j\in\mJ: \mI_j\cap\mI^H\neq \emptyset \right\}$, $\mJ^L:=\mJ\backslash\mJ^H$, and $\mJ^{HL}:=\left\{j\in\mJ: \mI_j\cap\mI^H\cap\mI^L\neq \emptyset \right\}$. Then, $\mJ^H$ is the set of job types that are processed by resources in heavy traffic, $\mJ^L$ is the set of job types that are processed by only resources in light traffic, and $\mJ^{HL}$ is the set of job types that are processed by both resources in light traffic and resources in light traffic. Finally, we assume that $\mI^H\neq\emptyset$, that is, there is at least one resource in heavy traffic. Otherwise, if all resources are in light traffic, any work conserving policy will perform reasonably well and the scheduling problem will become trivial. The reason is that work-conserving policies exist by Assumption \ref{a_h} (e.g., the GVM policy) and under work-conserving policies, the diffusion limits of the queue length processes associated with resources in light traffic converge to 0 (see for example Theorem 12(c) of \cite{gly90}). Therefore, if all resources are in light traffic, no cost will incur in the diffusion limit under work-conserving policies.

\subsection{Network Dynamics and Scheduling Control}\label{s_dynamics}

Let us fix arbitrary $t\in\R_+$ and $r\in\N_+$. For all $j\in\mJ$, we let $T_j^r(t)$ denote the cumulative amount of time that type $j$ jobs are processed during $[0,t]$. The scheduling control is defined by the process $(T_j^r,j\in\mJ)$. We let
\begin{equation}\label{eq_idle}
I_i^r(t):= t - \sum_{j\in\mJ_i}T_j^r(t)\qquad\forall i\in\mI,
\end{equation}
denote the cumulative time that resource $i$ idles during $[0,t]$.

For all $j\in\mJ$, we let $Q_j^r(t)$ denote the number of type $j$ jobs in the system at time $t$, including the job that is in service. Then, for all $j\in\mJ$,
\begin{equation}\label{eq_queue}
Q_j^r(t)= Q_j^r(0)+E_j^r(t)-S_{j}(T_{j}^r(t)) \geq 0,
\end{equation}
where $S_{j}(T_{j}^r(t))$ denotes the cumulative number of type $j$ jobs processed up to time $t$. We assume that the initial queue length vector, $(Q_j^r(0),j\in\mJ)$, is a random vector independent of all stochastic primitives and takes values in $\N^J$.

For all $j\in\mJ$, we have
\begin{equation}\label{eq_hl}
V_j(S_j(T_j^r(t))) \leq T_j^r(t) < V_j(S_j(T_j^r(t))+1),
\end{equation}
which implies that we consider only head-of-the-line (HL) policies, where jobs are processed in first-in-first-out order within each buffer.

\begin{definition}\label{d_admissible} (\textit{Admissible policies})
In the $r$th system, a scheduling policy $\pi^r:=(T_j^r,j\in\mJ)$ is admissible if the process $(I_i^r,T_j^r, Q_j^r,\;i\in\mI,j\in\mJ)$ satisfy \eqref{eq_idle}, \eqref{eq_queue}, \eqref{eq_hl}, and the following conditions for all $i\in\mI$, $j\in\mJ$, and $t\in\R_+$:
\begin{subequations}\label{eq_admissible}
\begin{align}
&\text{$T_{j}^r(t)$ is $\mF$-measurable (that is, $T_{j}^r(t)\in\mF$}).\label{eq_admissible_1}\\
&\text{Both $T_{j}^r$ and $I_{i}^r$ are continuous and nondecreasing and $T_j^r(0)=I_i^r(0)=0$}.\label{eq_admissible_2}\\
&\text{If $i\in\mI^H$, then $I_{i}^r(t)$ increases if and only if $\sum_{j\in\mJ_i} Q_j^r(t)=0$}.\label{eq_admissible_3}
\end{align}
\end{subequations}
\end{definition}
Condition \eqref{eq_admissible_1} implies that the set of admissible policies includes even the ones that can anticipate the future (recall that the $\sigma$-field $\mF$ is defined in Section \ref{s_not}). Condition \eqref{eq_admissible_3} implies that resources in heavy traffic are utilized in a work-conserving fashion, thus there is no capacity loss in the system. In contrast, the resources in light traffic are allowed to idle even when there are jobs in the system that they can process alone or in collaboration with other resources, because of their excess processing capacity. Condition \eqref{eq_admissible_3} keeps the set of admissible policies within the set of hierarchical greedy control policies introduced by \cite{har14}, under which resource utilization is the primary goal and the cost minimization is the secondary goal.

The set of admissible policies includes the ones which require divisible resources. The reason is that Definition \ref{d_admissible} allows $T_j^r(t)$ and $T_l^r(t)$ to increase simultaneously for some $i\in\mI$, $t\in\R_+$, and $j,l\in\mJ_i$ such that $j\neq l$, that is, resource $i$ can process both type $j$ and type $l$ jobs at the same time. However, in Section \ref{s_policy}, we will propose a policy which does not require resources to be divisible.

Condition \eqref{eq_admissible_3} implicitly forces the admissible policies to be preemptive. To see this, let us consider the network in Figure \ref{f_bc}. Suppose that both resources are in heavy traffic and consider a time epoch $t\in\R_+$ in which $Q_1^r(t)=Q_2^r(t)=0$, $Q_3^r(t)>0$, resource 2 is processing a type 3 job, and there is a type 1 job arrival. Due to Condition \eqref{eq_admissible_3}, resource 2 has to preempt the process of the type 3 job and start processing the newly arrived type 1 job together with resource 1. Because nonpreemptive policies are generally outperformed by the preemptive ones (see Section \ref{s_conc} for details and some numerical examples), the aforementioned restriction imposed by Condition \eqref{eq_admissible_3} is less likely to effect the strength of the asymptotic lower bound that we will propose in Section \ref{s_alb}.

Our objective is to minimize the expected total discounted holding cost. Let $h_j>0$ denote the holding cost rate per a type $j$ job per unit time for all $j\in\mJ$. Let $\delta > 0$ be the discount parameter and $\Pi^r$ denote the set of admissible policies in the $r$th system for all $r\in\N_+$. For notational convenience, we denote the process $X$ in the $r$th system under the admissible policy $\pi^r\in\Pi^r$ by $X^{\pi,r}$. Then, we want to find
\begin{equation}\label{eq_obj}
\argmin_{\pi^r\in\Pi^r} \sum_{j\in\mJ} h_j  \E\left[ \int_0^\infty \e^{-\delta t}Q_j^{\pi,r}(t) \dr t\right].
\end{equation}
Solving \eqref{eq_obj} exactly is very challenging. Therefore, in the next section, we will introduce fluid and diffusion scaled processes and henceforth consider the diffusion scaled processes in the asymptotic regime.

\subsection{Fluid and Diffusion Scaled Processes}\label{s_fluid}

For all $i\in\mI$, $t\in\R_+$, and $r\in\N_+$, we define the workload process associated with resource $i$ as
\begin{equation}\label{eq_workload}
W_i^r(t):=\sum_{j\in\mJ_i} \frac{Q_{j}^r(t)}{\mu_{j}}.
\end{equation}
Observe that $W_i^r(t)$ is the expected time that resource $i$ should spend in order to process all of the associated jobs given that no more jobs arrive in the system and all other resources that are needed for collaboration are available. 

For all $i\in\mI$, $j\in\mJ$, $t\in\R_+$, and $r\in\N_+$, the fluid scaled processes are defined as 
\begin{subequations}\label{eq_fluid}
\begin{align}
& \bar{E}_j^r(t):=r^{-2}E_j^r(r^2t)&& \bar{S}_j^r(t):=r^{-2}S_j(r^2t), \label{eq_fluid_1}\\
&\bar{T}_j^r(t):=r^{-2}T_j^r(r^2t) && \bar{I}_{i}^r(t):=r^{-2}I_{i}^r(r^2t),\label{eq_fluin_2}\\
&\bar{Q}_j^r(t):=r^{-2}Q_j^r(r^2t) &&\bar{W}_i^r(t):=r^{-2}W_i^r(r^2t).\label{eq_fluin_3}
\end{align}
\end{subequations}

For all $i\in\mI$, $j\in\mJ$, $t\in\R_+$, and $r\in\N_+$, the diffusion scaled processes are defined as
\begin{subequations}\label{eq_diffusion}
\begin{align}
& \hat{E}_j^r(t):=r\left(\bar{E}_j^r(t)-\la_j^r t\right)&&\hat{S}_j^r(t):=r\left(\bar{S}_j^r(t)-\mu_jt\right), \label{eq_diffusion_1}\\
& \hat{T}_j^r(t):=r\bar{T}_j^r(t) && \hat{I}_i^r(t):=r \bar{I}_i^r(t),\label{eq_diffusion_2}\\
& \hat{Q}_j^r(t):=r \bar{Q}_j^r(t) && \hat{W}_i^r(t):=r \bar{W}_i^r(t).\label{eq_diffusion_3}
\end{align}
\end{subequations}

For all $i\in\mI$, $t\in\R_+$, and $r\in\N_+$, let
\begin{equation}\label{eq_workload_2}
\hat{X}_i^r(t) := \sum_{j\in\mJ_i} \frac{1}{\mu_{j}}\left(\hat{Q}_j^r(0) + \hat{E}_j^r(t) -  \hat{S}_{j}^r\circ\bar{T}_{j}^r(t)\right) + r\left( \rho_i^r-1\right)t.
\end{equation}
After some algebra, we have
\begin{equation}\label{eq_workload_3}
\hat{W}_i^r=\hat{X}_i^r + \hat{I}_i^r,\qquad \forall i\in\mI,\;r\in\N_+.
\end{equation}
Because the resources in heavy traffic are utilized in a work-conserving fashion and by the minimality property of the reflection map (see Definition 14.2.1 of \cite{whi02}), we have
\begin{align}
&\left(\hat{W}_i^r,\hat{I}_i^r\right) = \left(\Phi,\Psi\right)\left(\hat{X}_i^r\right),\qquad\forall i\in\mI^H,\;r\in\N_+,\label{eq_min}\\
&\left(\hat{W}_i^r,\hat{I}_i^r\right) \geq \left(\Phi,\Psi\right)\left(\hat{X}_i^r\right),\qquad\forall i\in\mI^L,\;r\in\N_+.\nonumber
\end{align}

We assume that $\bar{Q}_j^r(0)\xrightarrow{a.s.} 0$ and $\hat{Q}_j^r(0)\Rightarrow \ot{Q}_j(0)$ for all $j\in\mJ$ such that $\big(\ot{Q}_j(0),\; j\in\mJ\big)$ is a random vector in $R_+^J$ and $\ot{Q}_j(0)= 0$ for all $j\in\mJ^L\cup\mJ^{HL}$. We let $\tilde{W}_i(0):= \sum_{j\in\mJ_i} \tilde{Q}_j(0)/\mu_j$ for all $i\in\mI$. Let $a_1:=\left|\mI^H\right|$ and $R$ be an $(a_1\times a_1)$-dimensional identity matrix. Let us define the $a_1$-dimensional vector $\Theta:=\left(\theta_i,i\in\mI^H\right)$ and the $(a_1\times a_1)$-dimensional matrix $\Sigma$ such that
\begin{equation}\label{eq_cov}
\Sigma_{ik}:= \sum_{j\in\mJ_i\cap\mJ_k}\la_j\left( \beta_j^2 + \sigma_{j}^2\right) \mu_{j}^{-2}, \qquad\forall i,k\in\mI^H.
\end{equation}
We have the following weak convergence result stating that under any sequence of admissible policies, the workload processes associated with the resources in heavy traffic weakly converge to a stochastic process with invariant distribution.

\begin{lemma}\label{l_diffusion}
Let $\pi=\{\pi^r,r\in\N_+\}$ be an arbitrary sequence of admissible policies, that is, $\pi^r\in\Pi^r$ for all $r\in\N_+$. Then,
\begin{equation*}
\left(\hat{W}_i^{\pi,r},i\in\mI^H\right) \Rightarrow \left(\ot{W}_i^*, i\in\mI^H\right),
\end{equation*}
where $\big(\ot{W}_i^*, i\in\mI^H\big)$ is a semimartingale reflected Brownian motion (SRBM) associated with the data $\left(\R_+^{a_1},\Theta, \Sigma, R \right)$. $\R_+^{a_1}$ is the state space of the SRBM; $\Theta$ and $\Sigma$ are the drift vector and the covariance matrix of the underlying Brownian motion of the SRBM, respectively; $R$ is the reflection matrix; and the SRBM starts at $\big(\tilde{W}_i(0),i\in\mI^H\big)$ in $\R^{a_1}_+$.
\end{lemma}
The proof of Lemma \ref{l_diffusion} is presented in E-companion EC.1.1. See Definition 3.1 of \cite{wil98b} for the formal definition of SRBM. The proof of Lemma 1 follows from standard weak convergence arguments. Specifically, the process $\big(\hat{X}_i^{\pi,r},i\in\mI^H\big)$ weakly converges to a Brownian motion by continuous mapping theorem, functional strong law of large numbers, functional central limit theorem etc. Then, the weak convergence of $\big(\hat{W}_i^{\pi,r},i\in\mI^H\big)$ follows from the continuity of the multidimensional reflection map. Lemma \ref{l_diffusion} will facilitate the derivation of the proposed policy and an asymptotic lower bound on the performance of sequences of admissible policies.

Finally, we will consider the following diffusion scaled version of the objective \eqref{eq_obj} henceforth:
\begin{equation}\label{eq_obj_2}
\argmin_{\pi^r\in\Pi^r} \sum_{j\in\mJ} h_j  \E\left[ \int_0^\infty \e^{-\delta t}\hat{Q}_j^{\pi,r}(t) \dr t\right].
\end{equation}

\section{Proposed Policy}\label{s_policy}

We propose a state-dependent, preemptive, and, continuous-review control policy under which an LP and a quadratic program (QP) are solved at discrete time epochs. First, we construct an approximating BCP whose solution will help us to derive a control policy. At this point, let us assume that
\begin{equation}\label{eq_hc}
\left( \hat{Q}_j^r, j\in\mJ\right)\Rightarrow \left( \ot{Q}_j, j\in\mJ\right).
\end{equation}
Then, by \eqref{eq_queue}, \eqref{eq_workload}, \eqref{eq_obj_2}, \eqref{eq_hc}, and Lemma \ref{l_diffusion}, we construct the following BCP: For all $t\in\R_+$ and $\omega\in\Omega$,
\begin{subequations}\label{eq_bcp}
\begin{align}
\min\;& \sum_{j\in\mJ^H} h_j \ot{Q}_j(t,\omega)\label{eq_bcp_1}\\
\text{s.t.}\; & \sum_{j\in\mJ_i} \frac{\ot{Q}_j(t,\omega)}{\mu_{j}} = \ot{W}_i^*(t,\omega)\qquad\forall i\in\mI^{H},\label{eq_bcp_2}\\
& \ot{Q}_j(t,\omega)\geq 0\qquad\forall j\in\mJ^H,\label{eq_bcp_3}
\end{align}
\end{subequations}
where the decision variables are $\big(\ot{Q}_j(t,\omega),j\in\mJ^H\big)$. The objective \eqref{eq_bcp_1} minimizes the total holding cost rate at time $t$ and under sample path $\omega$ associated with the job types processed by resources in heavy traffic. The constraints \eqref{eq_bcp_2} and \eqref{eq_bcp_3} state that we should split the workload of each resource that is in heavy traffic to its associated buffers in order to minimize the instantaneous total holding cost rate. Because the resources in light traffic have excess capacity, we do not expect to have nonzero amount of jobs in the buffers that are processed only by resources in light traffic in the diffusion limit. Hence, we ignore the job types in the set $\mJ^L$ in the BCP \eqref{eq_bcp}. 

For fixed $t\in\R_+$ and $\omega\in\Omega$, the BCP \eqref{eq_bcp} is an LP. Specifically, for given $w:=\big(w_i,i\in\mI^H\big)\in\R_+^{a_1}$, consider the LP
\begin{subequations}\label{eq_lp}
\begin{align}
\min\;& \sum_{j\in\mJ^H} h_j q_j\label{eq_lp_1}\\
\text{s.t. }\; & \sum_{j\in\mJ_i} \frac{q_{j}}{\mu_{j}} = w_i\qquad\forall i\in\mI^H,\label{eq_lp_2}\\
& q_{j}\geq 0\qquad\forall j\in\mJ^H,\label{eq_lp_3}
\end{align}
\end{subequations}
where the decision variables are $\big(q_{j},j\in\mJ^H\big)$. Let $z:\R_+^{a_1}\rightarrow \R_+$ be such that $z(w)$ denotes the optimal objective function value of the LP \eqref{eq_lp} for all $w\in\R_+^{a_1}$. Then, for any given $w\in\R_+^{a_1}$, let us consider the following QP:
\begin{subequations}\label{eq_qp}
\begin{align}
\min\;& \sum_{j\in\mJ^H} q_j^2\label{eq_qp_1}\\
\text{s.t. }\; & \sum_{j\in\mJ_i} \frac{q_{j}}{\mu_{j}} = w_i\qquad\forall i\in\mI^H,\label{eq_qp_2}\\
& \sum_{j\in\mJ^H} h_j q_j\leq z(w), \label{eq_qp_3}\\
& q_{j}\geq 0\qquad\forall j\in\mJ^H,\label{eq_qp_4}
\end{align}
\end{subequations}
where the decision variables are $\big(q_{j},j\in\mJ^H\big)$. Formally, for given $w\in\R_+^{a_1}$, the QP \eqref{eq_qp} finds the optimal solution of the LP \eqref{eq_lp} with the smallest Euclidean norm. However, our main reason for solving the QP \eqref{eq_qp} is a continuity property of its solution (see Section \ref{s_alb} for details). Because the QP \eqref{eq_qp} is convex, it is solvable in polynomial time (see \cite{vav08}).

The proposed policy keeps the numbers of the job types in the set $\mJ^L$ close to 0 and the numbers of the job types in the set $\mJ^H$ close to an optimal LP \eqref{eq_lp} solution under the LP parameter $(W_i^r(t),i\in\mI^H)$ at all times. Moreover, the proposed policy utilizes the resources in heavy traffic in a work-conserving fashion. Specifically, we start by solving LP \eqref{eq_lp} and QP \eqref{eq_qp}, then implement a scheduling rule in order to achieve the optimal buffer levels while utilizing the resources efficiently, next we resolve the LP \eqref{eq_lp} and QP \eqref{eq_qp} and repeat the same procedure.

Let us fix an arbitrary $r\in\N_+$ and a sample path $\omega\in\Omega$ and we will hide the notation $\omega$ for notational convenience in this section. Let $a_2:= \left|\mJ^H\right|$ and $(q_j^{*,r},j\in\mJ^H)\in\D^{a_2}$ denote the optimal solution process associated with the LP \eqref{eq_lp} under the parameter process $(W_i^r,i\in\mI^H)\in\D^{a_1}$ that is computed by the QP \eqref{eq_qp}. By \eqref{eq_workload} and \eqref{eq_lp_2},
\begin{equation}\label{eq_we}
\sum_{j\in\mJ_i} \frac{Q_j^r(t)}{\mu_j}= \sum_{j\in\mJ_i} \frac{q_j^{*,r}(t)}{\mu_j}=W_i^r(t),\qquad\forall i\in\mI^H,\;t\in\R_+.
\end{equation}
For all $i\in\mI^H$ and $t\in\R_+$, let 
\begin{align*}
\mJ^{>,r}(t)&:=\left\{j\in\mJ^H: Q_j^r(t)>\lceil q_j^{*,r}(t) \rceil \right\},\qquad \mJ_i^{>,r}(t) := \mJ_i\cap\mJ^{>,r}(t),\\
\mJ^{\leq,r}(t)&:=\left\{j\in\mJ^H: Q_j^r(t)\leq \lceil q_j^{*,r}(t) \rceil \right\},\qquad \mJ_i^{\leq,r}(t) := \mJ_i\cap\mJ^{\leq,r}(t).
\end{align*}
Then, $\{\mJ^{>,r}(t),\mJ^{\leq,r}(t)\}$ is a disjoint partition of $\mJ^H$ for all $t\in\R_+$ and $\{\mJ_i^{>,r}(t),\mJ_i^{\leq,r}(t)\}$ is a disjoint partition of $\mJ_i$ for all $i\in\mI^H$ and $t\in\R_+$. For all $t\in\R_+$, we let
\begin{align}
L_i^r(t) &:= \max_{j\in\mJ_i^{\leq,r}(t)} \frac{q_j^{*,r}(t)-Q_j^r(t)}{\la_j^r},\qquad\forall i\in\mI^H,\label{eq_review_3}\\
L^r(t)    &:=\max_{i\in\mI^H} L_i^r(t),\label{eq_review_4}
\end{align}
where $L^r(t)$ approximates the expected time required to achieve the optimal workload split among the resources in heavy traffic at time $t\in\R_+$. Roughly speaking, on the interval $[t, t+L^r(t))$, the proposed policy will take action to make the actual buffer levels close to $(q_j^{*,r}(t),j\in\mJ)$, at time $t+L^r(t)$, we will resolve the LP \eqref{eq_lp} and QP \eqref{eq_qp} and repeat the same procedure. We call the time interval between successive LP \eqref{eq_lp} solutions as \textit{review period}.

We will formally introduce the proposed policy in Section \ref{s_p_def}. Then, we present the reasoning behind the definitions in \eqref{eq_review_3} and \eqref{eq_review_4} in Section \ref{s_len}. 

\subsection{Formal Definition of the Proposed Policy}\label{s_p_def}

Let us fix an arbitrary $r\in\N_+$. Let, for all $t\in\R_+$, $x(t):=(x_j(t),j\in\mJ)$ and 
\begin{align*}
\mX(t):=\Bigg\{x(t)\in\R_+^J:\; &x_j(t)\in\{0,1\},\quad\forall j\in\mJ,&&  x_j(t)\leq \I\left(Q_j^r(t)>0\right),\quad\forall j\in\mJ,\\
&\sum_{j\in\mJ_i}x_j(t) \leq 1,\quad\forall i\in\mI^L,&& \sum_{j\in\mJ_i}x_j(t) = \I\left(W_i^r(t)>0\right),\quad\forall i\in\mI^H \Bigg\}.
\end{align*}
Hence, $\mX(t)$ denotes the set of binary resource allocation vectors under which a resource can be assigned to at most a single job type and resources in heavy traffic are utilized in a work-conserving fashion at time $t$. By Assumption \ref{a_h}, $\mX(t)\neq\emptyset$ for all $t\in\R_+$. 

For all $t\in\R_+$, the proposed policy implements a resource allocation vector chosen from the set $\mX(t)$ by using an index rule. For all $s,t\in\R_+$ and $x(s)\in\mX(s)$ such that $s\geq t$, let
\begin{equation}\label{eq_index}
p_1(x(s),t) :=\sum_{j\in\mJ^H} x_j(s)\I\left(Q_j^r(s)> \lceil q_j^{*,r}(t) \rceil \right) \left|\mI_j \cap \mI^H\right|,\qquad p_2(x(s)) :=\sum_{j\in\mJ} x_j(s) |\mI_j|.
\end{equation}
Then, $p_1(x(s),t)$ is the number of resources in heavy traffic processing job types whose number in the system at time $s$ exceeds the desired levels associated with time $t$, and $p_2(x(s))$ is the total number of busy resources at time $s$ under the resource allocation vector $x(s)$. Finally, let us define the partial order $\succeq$ in $\N^2$ such that for all $(m_1,n_1),(m_2,n_2)\in\N^2$, we have $(m_1,n_1)\succeq (m_2,n_2)$ if and only if $m_1>m_2$ or $(\{m_1=m_2\} \cap \{n_1\geq n_2\})$. 

Let $t\in\R_+$ denote the start time of a review period. Then, at time $s\geq t$, the proposed policy will choose a resource allocation vector from the set $\mX(s)$ which achieves the maximum $\left(p_1(x(s),t),p_2(x(s))\right)$ value with respect to the partial order $\succeq$. 

The formal definition of the proposed policy is presented below.

\begin{definition}\label{d_pp}
(\textit{The proposed policy})
\begin{description}[align=left]

\item [Step 1] Let $t\in\R_+$ denote the current time. Solve LP \eqref{eq_lp} and then QP \eqref{eq_qp} at time $t$. If $\mJ^{>,r}(t)=\emptyset$, go to Step 2. Otherwise, go to Step 3.

\item [Step 2] Utilize all of the resources in a work-conserving and admissible fashion. At the first time when an external arrival or a service completion happens, go to Step 1.

\item [Step 3] Let $t\in\R_+$ denote the current time. Compute $L_i^r(t)$ for all $i\in\mI^H$ by \eqref{eq_review_3} and $L^r(t)$ by \eqref{eq_review_4}. Let $\kappa>0$ denote the length of the time interval between time $t$ and the first time when an external arrival or service completion happens after time $t$. Then, Step 3 will last $L^r(t) \vee \kappa$ time units. For all $s\in[t,t+(L^r(t)\vee \kappa))$, implement an arbitrary resource allocation vector from the set
\begin{equation}\label{eq_set}
\left\{x(s)\in\mX(s): \left(p_1(x(s),t),p_2(x(s))\right) \succeq \left(p_1(x'(s),t),p_2(x'(s))\right),\quad\forall x'(s)\in\mX(s) \right\}.
\end{equation}
At time $t+(L^r(t)\vee \kappa)$, go to Step 1.

\end{description}
\end{definition}

We make a few remarks about the proposed policy. Step 1 in Definition \ref{d_pp} is done instantaneously and Step 2 or 3 is a review period. We solve the LP \eqref{eq_lp} and the QP \eqref{eq_qp} only at Step 1. Because the system controller should monitor the system state continuously in a review period to choose resource allocation vectors from the set in \eqref{eq_set} in Step 3 and to determine the ending time of the review period in Steps 2 and 3, the proposed policy has a continuous-review structure. Because Step 2 or 3 can end before the jobs that are in service are processed completely, the proposed policy is preemptive. The proposed policy is work-conserving in Step 2 and Step 2 lasts until the first external arrival or service completion epoch after it starts. Step 3 lasts at least until the first external arrival or service completion epoch after it starts. Therefore, there are finitely many review periods in a finite time interval. The proposed policy chooses resource allocation vectors from the set $\mX(t)$ at all times in Step 3. However, it is enough to make those choices at only service completion or external arrival epochs because the system state does not change in between those events.

Under the proposed policy, resources in heavy traffic are utilized in a work-conserving fashion, that is,
\begin{equation*}
\int_0^\infty W_i^r(t)\dr I_i^r(t)=0,\qquad\forall i\in\mI^H.
\end{equation*}
However, resources in light traffic are not necessarily utilized in a work-conserving fashion. Otherwise, $Q_j$ will be close to 0 for all $j\in\mJ^L\cup\mJ^{HL}$. This implies that resources in heavy traffic should not keep their workload in the buffers associated with the job types in the set $\mJ^{HL}$, which can increase the total cost especially when cost of keeping job types in the set $\mJ^{HL}$ is smaller than the cost of keeping job types in the set $\mJ^H\backslash\mJ^{HL}$. Nevertheless, under the proposed policy, resources in light traffic process the job types in the set $\mJ^L$ in a work-conserving fashion, that is, 
\begin{equation*}
\int_0^\infty \bigg(\sum_{j\in\mJ_i\cap\mJ^L} Q_j^r(t)\bigg) \dr I_i^r(t)=0,\qquad\forall i\in\mI^L.
\end{equation*}

\begin{remark}\label{r_unique}
If there are multiple resource allocation vectors in Step 2 that utilize all resources in a work-conserving fashion or there are multiple resource allocation vectors in Step 3 that maximize the index $(p_1,p_2)$, the proposed policy does not uniquely identify which vector to choose. Instead, it allows to choose any resource allocation vector with the aforementioned properties in Steps 2 and 3.
\end{remark}

We end this section by presenting an algorithm which efficiently finds a resource allocation vector in the set in \eqref{eq_set} in Step 3 of the proposed policy. Let $t\in\R_+$ denote the start time of Step 3 and $s\in[t,t+(L^r(t)\vee \kappa))$. Let
\begin{equation*}
\hat{\mJ}^{r}(t,s) := \left\{j\in\mJ^H: Q_j^r(s)>\lceil q_j^{*,r}(t) \rceil \right\},
\end{equation*}
that is, $\hat{\mJ}^{r}(t,s)$ is the set of job types that are processed by resources in heavy traffic and whose number in system at time $s$ is greater than the optimal LP \eqref{eq_lp} solution at time $t$. The algorithm that can be used at time $s$ in Step 3 of the proposed policy is as follows.

\begin{definition}\label{d_alg} (\textit{The algorithm designed for Step 3 of the proposed policy})
\begin{description}[align=left]

\item [Step 1] Compute $|\mI_j\cap\mI^H|$ for all job types in the set $\hat{\mJ}^{r}(t,s)$ and sort those job types starting from the largest $|\mI_j\cap\mI^H|$ value to the smallest one. Ties can be broken randomly or with respect to a fixed rule. Go to Step 2.

\item [Step 2] Assign the associated resources to the first job type in the ordered set. Go to Step 3.

\item [Step 3] Consider the next job type in the ordered set. Go to Step 4.

\item [Step 4] If none of the resources that are required in the service process of the current job type are assigned to the previous job types, assign those resources to the current job type. Otherwise, the current job type will not be processed. If there are remaining job types in the ordered set, go to Step 3. Otherwise go to Step 5.

\item [Step 5] If  $\mJ\setminus\hat{\mJ}^{r}(t,s)=\emptyset$, go to Step 8. Otherwise, compute $|\mI_j|$ for all job types in the set $\mJ\setminus\hat{\mJ}^{r}(t,s)$ and sort those job types starting from the largest $|\mI_j|$ value to the smallest one. Ties can be broken randomly or with respect to a fixed rule. Call the resulting ordered set as ``the second ordered set''. Consider the first job type in the second ordered set. Go to Step 7.

\item [Step 6] Consider the next job type in the second ordered set. Go to Step 7.

\item [Step 7] If all resources that are required in the service process of the current job type are available, assign those resources to the current job type. Otherwise, the current job type will not be processed. If there are remaining job types in the second ordered set, go to Step 6. Otherwise go to Step 8.

\item [Step 8] End the algorithm.

\end{description}
\end{definition}

The computational complexity of the algorithm in Definition \ref{d_alg} is $O(J\ln J)$, where $O(\cdot)$ denotes the big-$O$ notation. Therefore, the algorithm is computationally efficient. Assumption \ref{a_h} plays a crucial role in the following result whose proof is presented in E-companion EC.1.2.

\begin{lemma}\label{l_alg}
The algorithm in Definition \ref{d_alg} finds a resource allocation vector in the set in \eqref{eq_set} in Step 3 of the proposed policy.
\end{lemma}

\subsection{Length of a Review Period}\label{s_len}

In this section, we present the reasoning behind our choice for the length of a review period for the proposed policy. Let us consider an arbitrary time $t\in\R_+$. Suppose that $\mJ^{>,r}(t)\neq\emptyset$. Then, there exists a $j\in\mJ^H$ such that $Q_j^r(t)>\lceil q_j^{*,r}(t) \rceil\geq q_j^{*,r}(t)$. This implies that there exists an $l\in\mJ^H$ such that $Q_l^r(t)< q_l^{*,r}(t)\leq \ru{q_l^{*,r}(t)}$ by \eqref{eq_we}. Hence, $\mJ^{\leq,r}(t)\neq \emptyset$ by definition. For simplicity, let us assume that $q_j^{*,r}(t)$ is an integer for all $j\in\mJ^H$. We want to utilize the resources in order to decrease $Q_j^r$ to $q_j^{*,r}(t)$ for all $j\in\mJ^{>,r}(t)$, while keeping $Q_j^r$ less than or equal to $q_j^{*,r}(t)$ for all $j\in\mJ^{\leq,r}(t)$. Let $\breve{L}_i^r(t)\in\R_+$ denote the expected length of the time interval that resource $i$ needs to spend in order to achieve the aforementioned goal for all $i\in\mI^H$. If we ignore the collaboration requirements, $\breve{L}_i^r(t)$ should satisfy the following equalities:
\begin{align}
\breve{L}_i^r(t) &= \sum_{j\in\mJ_i^{>,r}(t)} \frac{Q_j^r(t)-q_j^{*,r}(t)}{\mu_j} + \breve{L}_i^r(t)\sum_{j\in\mJ_i^{>,r}(t)} \frac{\la_j^r}{\mu_j} + \sum_{j\in\mJ_i^{\leq,r}(t)} \frac{\left(\la_j^r \breve{L}_i^r(t) - q_j^{*,r}(t)+Q_j^r(t)\right)^+}{\mu_j}\label{eq_review_1}\\
&= \sum_{j\in\mJ_i}\frac{\left(\la_j^r \breve{L}_i^r(t) -q_j^{*,r}(t)+Q_j^r(t)\right)^+}{\mu_j}.\label{eq_review_2}
\end{align}
Observe that \eqref{eq_review_2} is a compact version of the right-hand-side (RHS) of \eqref{eq_review_1}. The first term in the RHS of \eqref{eq_review_1} denotes the average time that resource $i$ should spend to deplete the excess jobs associated with the job types in the set $\mJ_i^{>,r}(t)$. In the mean time, there will be external type $j$ job arrivals for all $j\in\mJ_i^{>,r}(t)$. Hence, the second term in the RHS of \eqref{eq_review_1} denotes the average time that resource $i$ should spend to process the excess jobs due to external job arrivals associated with the jobs in the set $\mJ_i^{>,r}(t)$. Lastly, the third term in the RHS of \eqref{eq_review_1} denotes the average time that resource $i$ should spend to process type $j$ jobs for all $j\in\mJ_i^{\leq,r}(t)$ if the average number of external type $j$ job arrivals is greater than $q_j^{*,r}(t) - Q_j^r(t)$. Then, we have the following result.

\begin{lemma}\label{l_review}
If $\la_j^r = \la_j$ for all $j\in\mJ$, that is, if the arrival rates are equal to the limiting ones, then, for all $i\in\mI^H$, $\breve{L}_i^r(t)\in\R_+$ is a solution of the equality \eqref{eq_review_1} if and only if
\begin{equation}\label{eq_review}
\breve{L}_i^r(t) \geq \max_{j\in\mJ_i^{\leq,r}(t)} \frac{q_j^{*,r}(t)- Q_j^r(t)}{\la_j}.
\end{equation}
\end{lemma}
The proof of Lemma \ref{l_review} is presented in E-companion EC.1.3. Lemma \ref{l_review} provides a lower bound on the expected amount of time that resource $i$ should spend in order to split its workload to the associated buffers according to an optimal LP \eqref{eq_lp} solution under the limiting arrival rates. By Lemma \ref{l_review} and considering Assumption \ref{a_regime} Part 1, we define $L_i^r(t)$ and $L^r(t)$ as in \eqref{eq_review_3} and \eqref{eq_review_4}, respectively.

Next, we present our main theoretical results.

\section{The Main Theoretical Results}\label{s_th}

We present an asymptotic lower bound in Section \ref{s_alb}. Then, we will prove that the proposed policy achieves that asymptotic lower bound and thus is asymptotically optimal in Section \ref{s_ao}.

\subsection{Asymptotic Lower Bound}\label{s_alb}

We present an asymptotic lower bound on the performance of any sequence of admissible policies. The following result presents some useful properties of the LP \eqref{eq_lp} and the QP \eqref{eq_qp}.

\begin{lemma}\label{l_lp}
\begin{enumerate}

\item For any given $w^{(1)}:=\big(w_i^{(1)},i\in\mI^H\big)\in\R_+^{a_1}$ and $w^{(2)}:=\big(w_i^{(2)},i\in\mI^H\big)\in\R_+^{a_1}$,
\begin{equation*}
\left| z\big(w^{(1)}\big)-z\big(w^{(2)}\big)\right| \leq C_1 \big| w^{(1)}-w^{(2)}\big|_{\infty},
\end{equation*}
where $C_1>0$ is a constant dependent only on the objective function coefficients and the left-hand-side (LHS) parameters of the constraints.

\item For each $w\in\R_+^{a_1}$, there exists a unique optimal solution of the QP \eqref{eq_qp}. For any given $w^{(1)}\in\R_+^{a_1}$ and $w^{(2)}\in\R_+^{a_1}$, let $\big(q_j^{(1)},j\in\mJ^H\big)$ and $\big(q_j^{(2)},j\in\mJ^H\big)$ be the optimal solutions of the QP \eqref{eq_qp} under $w^{(1)}$ and $w^{(2)}$, respectively. Then,
\begin{equation*}
\max_{j\in\mJ^H} \left| q_j^{(1)}-q_j^{(2)}\right| \leq C_2 \big| w^{(1)}-w^{(2)} \big|_{\infty},
\end{equation*}
where $C_2>0$ is a constant dependent only on the objective function coefficients and the LHS parameters of the constraints. 
\end{enumerate}
\end{lemma}
The first part of Lemma \ref{l_lp} states that the optimal objective function value of the LP \eqref{eq_lp} is Lipschitz continuous in the RHS parameter $w\in\R_+^{a_1}$. For a proof, see the derivation of the equation (10.22) of \cite{sch98}, which states the aforementioned continuity result for general LPs. Because we will solve LP \eqref{eq_lp} regularly over time (at discrete time epochs) and LP \eqref{eq_lp} may have multiple optimal solutions at some time epochs, we need to choose an optimal solution among the set of optimal solutions at those time epochs such that the optimal solutions that we will use over time will not fluctuate a lot. The second part of Lemma \ref{l_lp} presents a method to achieve the aforementioned goal. The second part of Lemma \ref{l_lp} states that the optimal solution of the QP \eqref{eq_qp} is unique and Lipschitz continuous in the RHS parameter $w\in\R_+^{a_1}$ and this result follows directly from Proposition 4.1.d of \cite{han12}. A direct consequence of the second part of Lemma \ref{l_lp} is the following Lipschitz continuity result.

\begin{lemma}\label{l_reg}
For any given nonnegative parameter process $(\bm{w}(t),t\in\R_+)\in\D^{a_1}$, let $\big(q_j^*(t),j\in\mJ^H,t\in\R_+)\in\D^{a_2}$ denote the optimal solution process associated with the LP \eqref{eq_lp} derived by the QP \eqref{eq_qp}. For all $s,t\in\R_+$, we have
\begin{equation*}
\max_{j\in\mJ^H} \left| q_j^*(s)-q_j^*(t)\right| \leq C_2 \big| \bm{w}(s)-\bm{w}(t) \big|_{\infty}.
\end{equation*}
\end{lemma}

\begin{remark}\label{r_sol}
Any optimal solution process associated with the LP \eqref{eq_lp} that satisfy the Lipchitz continuity result in Lemma \ref{l_reg} is enough for the policy that we will propose to be asymptotically optimal. Solving the QP \eqref{eq_qp} is just one way to find the optimal LP \eqref{eq_lp} solutions with the desired Lipschitz property.
\end{remark}

Finally, we prove that the optimal objective function value of the LP \eqref{eq_lp} provides an asymptotic lower bound on the performance of admissible policies with respect to the objective \eqref{eq_obj_2}.

\begin{theorem}\label{t_lb}
Let $\pi=\{\pi^r,r\in\N_+\}$ be an arbitrary sequence of admissible policies. Then,
\begin{equation*}
\liminf_{r\rightarrow\infty} \sum_{j\in\mJ} h_j  \E\left[ \int_0^\infty \e^{-\delta t} \hat{Q}_j^{\pi,r}(t) \dr t\right] \geq \E\left[ \int_0^\infty \e^{-\delta t} z\big(\ot{W}^*(t)\big) \dr t\right],
\end{equation*}
where $\ot{W}^*:=\big(\ot{W}_i^*,i\in\mI^H\big)$ is defined in Lemma \ref{l_diffusion}.
\end{theorem}
The proof of Theorem \ref{t_lb} is presented in E-companion EC.2 and follows from the continuity of the optimal objective function value of the LP \eqref{eq_lp} (see Lemma \ref{l_lp} Part 1) and Fatou's lemma. Any sequence of admissible policies that achieves the lower bound in Theorem \ref{t_lb} is asymptotically optimal.

\subsection{Asymptotic Optimality of the Proposed Policy}\label{s_ao}

The following theorem states that the proposed policy achieves the asymptotic lower bound in Theorem \ref{t_lb}. 

\begin{theorem}\label{t_ao}
Consider the following conditions about the initial system state.
\begin{enumerate}
\item There exists an $r_0\in\N_+$ such that
\begin{equation}\label{eq_ini1}
\sup_{r\geq r_0} \max_{j\in\mJ}\E\left[ \left(\hat{Q}_j^r(0)\right)^2\right]=:C_0 <\infty.
\end{equation}

\item For all $\ep>0$,
\begin{equation}\label{eq_ini2}
\pr\left(\max_{j\in\mJ^H} \left| Q_j^{r}(0) - q_j^{*,r}(0)\right| > \ep r \right) \rightarrow 0.
\end{equation}
\end{enumerate}

If the conditions above hold, then, under the proposed policy,
\begin{equation*}
\lim_{r\rightarrow\infty} \sum_{j\in\mJ} h_j  \E\left[ \int_0^\infty \e^{-\delta t} \hat{Q}_j^{r}(t) \dr t\right] = \E\left[ \int_0^\infty \e^{-\delta t} z\big(\ot{W}^*(t)\big) \dr t\right],
\end{equation*}
where $\ot{W}^*=\big(\ot{W}_i^*,i\in\mI^H\big)$ is defined in Lemma \ref{l_diffusion}.
\end{theorem}
The proof of Theorem \ref{t_ao} is presented in E-companion EC.3 and follows from the continuity of the optimal LP \eqref{eq_lp} solutions (see Lemma \ref{l_reg}) and large deviations theory. Specifically, by large deviations theory, we show that under the proposed policy, the queue length processes track the lagged optimal LP \eqref{eq_lp} solutions sufficiently closely and the time lag is sufficiently small. Because the optimal LP \eqref{eq_lp} solutions are Lipschitz continuous, they do not change significantly within short time intervals, and thus the queue length processes track the non-lagged optimal LP \eqref{eq_lp} solutions sufficiently closely over time under the proposed policy. Consequently, we obtain Theorem \ref{t_ao}. Theorems \ref{t_lb} and \ref{t_ao} imply that the proposed policy is asymptotically optimal with respect to the objective \eqref{eq_obj_2}.

Next, we will present applications of the proposed policy in two different networks.

\section{Examples}\label{s_ex}

We demonstrate the proposed policy in two different examples. The first example shows that if there is a single resource in the system, then the proposed policy becomes equivalent to the $c\mu$-rule. The second example demonstrates the superior performance of the proposed policy against policies from the literature.

\subsection{Example 1: The $c\mu$-rule}\label{s_ex_1}

Suppose that all job types require the same resources in their service processes, that is, $\mI_j=\mI_l$ for all $j,l\in\mJ$. Because the coordination of the resources is trivial in this case, without loss of generality, we can assume that there is a single resource in the network, that is, $I=1$. Therefore, the network is equivalent to a multi-class $GI/GI/1$ queueing system. Let us assume that the single resource is in heavy traffic. Then, in this case, the LP \eqref{eq_lp} is:
\begin{equation*}
\min\left\{\sum_{j\in\mJ} h_j q_j:  \sum_{j\in\mJ} \frac{q_{j}}{\mu_{j}} = w,\quad q_{j}\geq 0,\;\forall j\in\mJ.\right\}
\end{equation*}
Without loss of generality, let us assume that $h_1\mu_1\geq h_2\mu_2\geq\ldots\geq h_J\mu_J$. Then, an optimal solution of the LP \eqref{eq_lp} satisfying Lemma \ref{l_reg} (recall Remark \ref{r_sol}) is
\begin{equation}\label{eq_cmu_s}
q_1^*=q_2^*=\ldots=q_{J-1}^*=0,\quad q_J^*=\mu_J w.
\end{equation}
Therefore, all of the workload is kept in the cheapest buffer. By recalling Remark \ref{r_unique}, the proposed policy associated with the solution in \eqref{eq_cmu_s} includes the preemptive $c\mu$-rule.

\subsection{Example 2: The network in Figure \ref{f_bc}}\label{s_ex_2}

Consider the network in Figure \ref{f_bc}. Suppose that both resources are in heavy traffic. In this case, the LP \eqref{eq_lp} is:
\begin{equation*}
\min\left\{\sum_{j=1}^3 h_j q_j:  \frac{q_1}{\mu_1}+ \frac{q_2}{\mu_2}= w_1,\quad\frac{q_1}{\mu_1}+ \frac{q_3}{\mu_3}= w_2,\quad q_{j}\geq 0,\;\forall j\in\{1,2,3\}.\right\}
\end{equation*}
There are two cases to consider. 
\begin{enumerate}
\item Suppose that $h_1\mu_1 \geq h_2\mu_2 + h_3\mu_3$. Then, an optimal solution of the LP \eqref{eq_lp} that satisfy Lemma \ref{l_reg} (recall Remark \ref{r_sol}) is 
\begin{equation}\label{eq_bc_s1}
q_1^*=0,\quad q_2^*=\mu_2 w_1,\quad q_3^*=\mu_3 w_2.
\end{equation}
The proposed policy associated with the solution in \eqref{eq_bc_s1} gives type 1 jobs preemptive priority over the type 2 and 3 jobs, and includes (recall Remark \ref{r_unique}) the GVM policy.

\item Suppose that $h_1\mu_1 \leq h_2\mu_2 + h_3\mu_3$. Then, an optimal solution of the LP \eqref{eq_lp} that satisfy Lemma \ref{l_reg} (recall Remark \ref{r_sol}) is 
\begin{equation*}
q_1^*=\mu_1(w_1\wedge w_2),\quad q_2^*=\mu_2(w_1-w_2)^+,\quad q_3^*=\mu_3(w_2-w_1)^+.
\end{equation*}
Therefore, the proposed policy keeps as much workload as possible in buffer 1 while utilizing the resources in a work-conserving fashion.
\end{enumerate}

We present a simulation experiment associated with the second case above. Suppose that the inter-arrival and service times are exponentially distributed, $\mu_j=2$, and $\la_j=\la<1$ for all $j\in\{1,2,3\}$. Suppose that $h_j=1$ for all $j\in\{1,2,3\}$ and thus let us consider the objective of minimizing the long-run average number of jobs in the system, which is a natural objective in simulation experiments. We compare the performance of the proposed policy with the ones of two static prioritization policies, namely the GVM policy and the PIA policy, the performance of the \textit{preemptive hierarchical threshold} (PHT) policy, and the performance of the \textit{maximum pressure} (MP) policy. 

Under the PIA policy, the resources give preemptive priority to the job types that they process individually. For example, if a type 2 job arrives when there is not any type 2 and 3 job in the system and the resources are processing a type 1 job, then R1 immediately starts the process of the type 2 job while R2 idles even though there is a type 1 job in the system (because there is no type 3 job in the system). The PHT policy is suggested by \cite{gur18} (see Section 7 therein). Under the PHT policy, if the number of type 1 jobs reaches a threshold level, the resources give preemptive priority to type 1 jobs until the first time there is no type 1 job in the system. At that moment, the resources start giving preemptive priority to the job types that they individually process until the number of type 1 jobs reach the threshold value again. Observe that, the PHT policy acts like the GVM policy or the PIA policy depending on the system state. If the threshold value is equal to 1, then the PHT policy becomes equivalent to the GVM policy. Finally, we test the performance of the MP policy proposed by \cite{dai05}. Under the MP policy, if the number of type 1 jobs is greater than or equal to the summation of the numbers of type 2 and 3 jobs, the resources give preemptive priority to type 1 jobs. Otherwise, the resources give preemptive priority to the type 2 and 3 jobs. For example, under the MP policy, if there is no type 3 job in the system and the number of type 1 jobs is less than the number of type 2 jobs, resource 1 processes type 2 jobs and resource 3 idles. Due to the violation of Condition \eqref{eq_admissible_3}, the PIA, the PHT, and the MP policies are not admissible. In contrast, the GVM policy is admissible.

\begin{figure}[ht]
  \centering
  \subfloat[$\# X$ where  $X\in\{\text{Pro., GVM, MP, PHT, PIA}\}$.]{\label{f1}\includegraphics[width=0.5\textwidth]{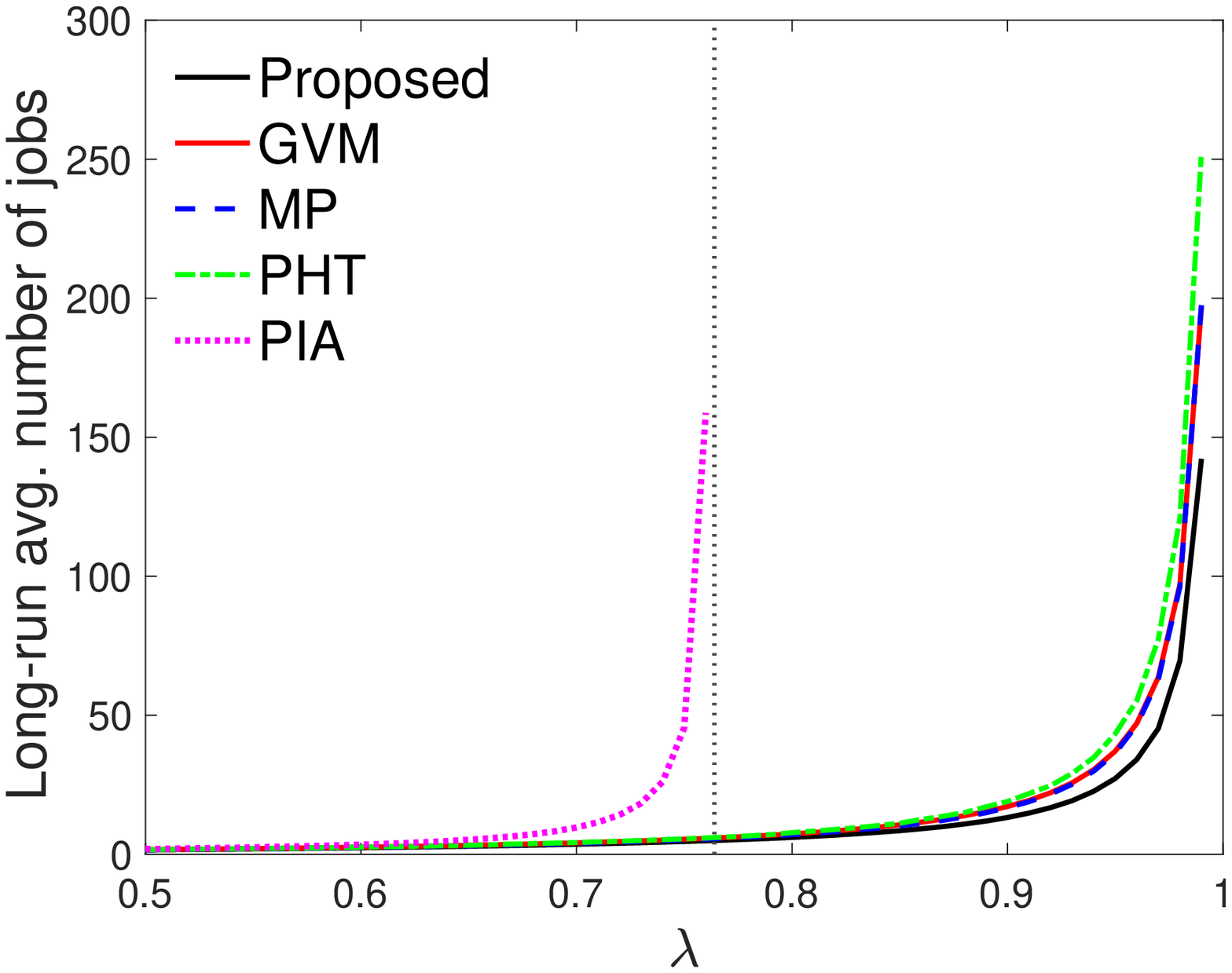}}
  \subfloat[$100\times (\# X - \# Proposed)/ \# Proposed$ where $X\in\{\text{GVM, MP, PHT}\}$.]{\label{f2}\includegraphics[width=0.5\textwidth]{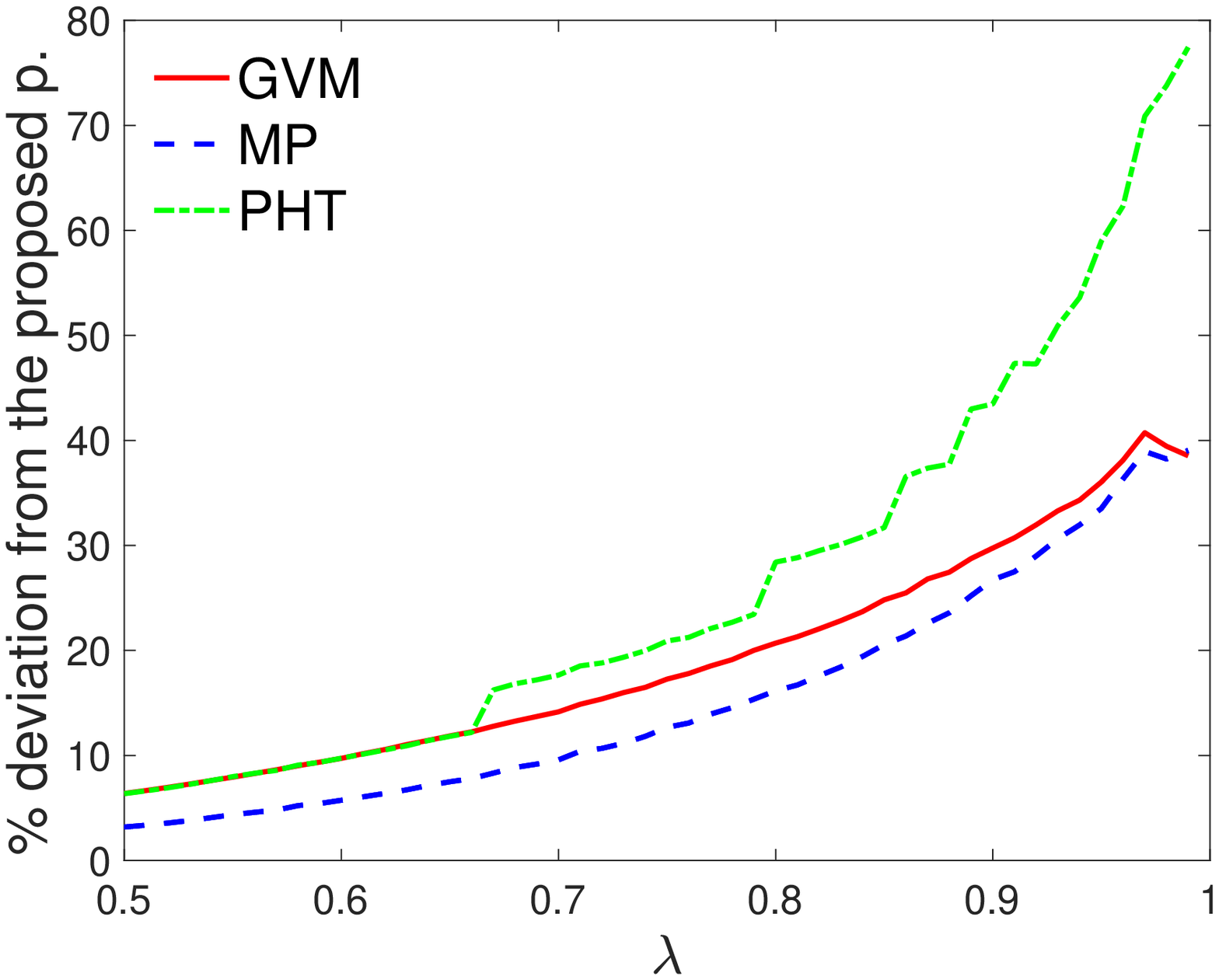}}\\
  \caption{(Color online) Performances of the proposed, the GVM, the MP, the PHT, and the PIA policies.}
  \label{f0}
\end{figure}

Let $\# X$ denote the long-run average number of jobs in the system under policy X. Figure \ref{f1} shows $\# Proposed$, $\# GVM$, $\# MP$, $\# PHT$, and $\# PIA$ for $\la\in\{0.5,0.51,\ldots,0.99\}$. For each $\la$ value, we use 4 different threshold values for the PHT policy and present the result associated with the best performing one. By considering the results of \cite{gur18} (see Section 7 therein), we choose those threshold values in the order of $(1-\la)^{-1}$. Specifically, for given $\la$, the 4 threshold values are given by the formula
\begin{equation*}
\left[ \frac{\ga}{1-\la}\right],\quad \text{where } \ga\in\{0.5,1,1.5, 2\}.
\end{equation*}

According to the results, the proposed policy outperforms the other four policies. Figure \ref{f2} shows the relative performances of the GVM, the MP, and the PHT policies with respect to the performance of the proposed policy. Let us consider the states in which all three job types are in the system. Under the proposed policy, the resources often give priority to the type 2 and 3 jobs and so the throughput is often 4. In contrast, under the GVM policy, the resources give priority to the type 1 jobs and so the throughput is only 2. Under the MP policy, if the number of type 1 jobs is greater than the summation of the numbers of type 2 and 3 jobs, the resources give priority to the type 1 jobs and thus the throughput is only 2 in those cases. Consequently, the proposed policy outperforms both the GVM and the MP policies. 

Figure \ref{f1} shows that the system becomes unstable when $\la\geq 0.77$ under the PIA policy. Therefore, there is capacity loss under the PIA policy, that is, the system capacity is strictly less than the bottleneck capacity. If there is a type 1 and a type 2 job in the system and there is no type 3 job in the system, resource 2 idles even though there is a type 1 job in the system because resource 1 processes a type 2 job. Similarly, resource 1 can idle when there is a type 1 job in the system. Hence, the resources are not fully utilized under the PIA policy which results in a capacity loss (see Section 3.1 of \cite{gur18} for a theoretical explanation of the aforementioned capacity loss). In contrast, both the proposed and the GVM policies are work-conserving. Although the MP policy is not work-conserving, it is path-wise stable (see Theorem 2 of \cite{dai05}) and thus there is no capacity loss under the MP policy as seen in Figure \ref{f1}.

\begin{figure}[htb]
\begin{center}
\includegraphics[width=0.5\textwidth]{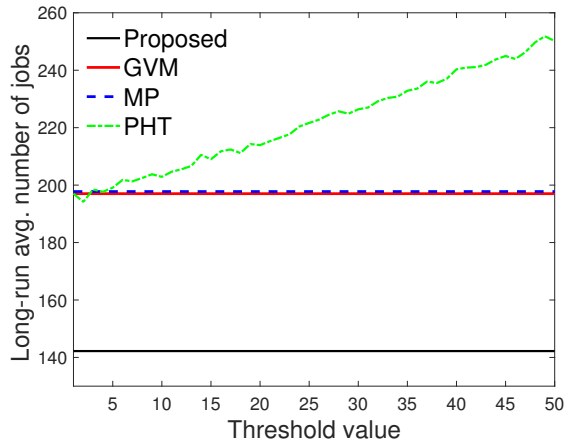}
\caption{(Color online) $\# X$ where $X\in\{\text{Pro., GVM, MP, PHT}\}$ for $\la=0.99$ and threshold values $\{1,2,\ldots,50\}$.}\label{f4}
\end{center}
\end{figure}

According to Figure \ref{f2}, the PHT policy performs worse than the GVM policy does. Recall that the PHT policy acts like the GVM policy or the PIA policy depending on the system state. Because the resources are not fully utilized under the PIA policy, the PHT policy does not benefit by acting like the PIA policy and thus it does not outperform the GVM policy. A potential reason for the underperformance of the PHT policy is that the associated threshold value is not optimally chosen. Because the threshold value associated with $\ga=0.5$ results in the best performing PHT policy among the threshold values associated with $\ga\in\{0.5,1,1.5,2\}$ for all $\la\in\{0.5,0.51,\ldots,0.99\}$, Figure \ref{f4} shows $\# PHT$ for $\la = 0.99$ and the threshold values in the set $\{1,2,\ldots,50\}$ (equivalently $\ga\in\{0.01,0.02,\ldots,0.5\}$). According to Figure \ref{f4}, the PHT policy never outperforms the proposed policy and slightly outperforms the GVM and the MP policies only when the threshold value is equal to 2.

\section{Networks with Capacity Loss}\label{s_ext}

The proposed policy is designed for and asymptotically optimal in parallel networks with hierarchical collaboration architecture (see Definition \ref{d_hca} and Assumption \ref{a_h}). Recall that there is no capacity loss in those networks. In this section, we will show that the proposed policy is implementable and asymptotically optimal in a class of parallel networks with capacity loss. First, we will illustrate the idea behind this extension in the simple network in Figure \ref{f_bc2}. Then, we will generalize and formalize the extension.

\begin{figure}[htb]
\begin{center}
\includegraphics[width=0.4\textwidth]{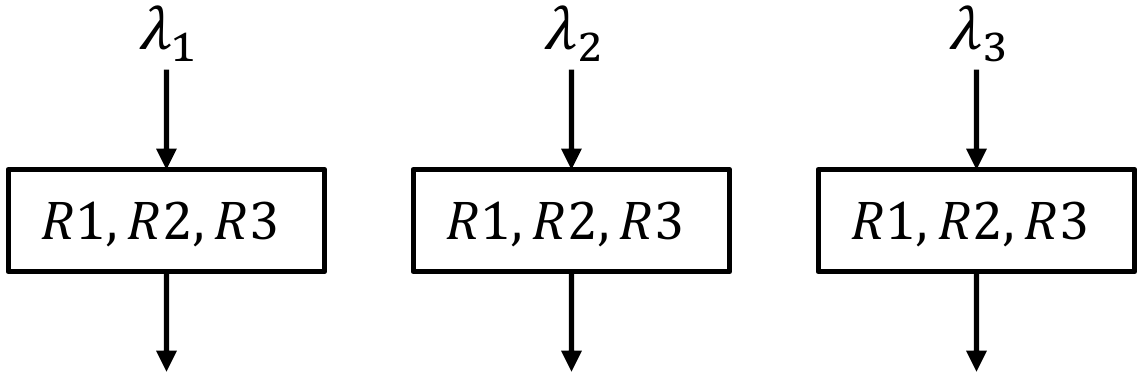}
\caption{A parallel network with three job types and three multitasking resources (R1, R2, R3).}\label{f_bc3}
\end{center}
\end{figure}

Let us consider the parallel network in Figure \ref{f_bc2}. There is capacity loss in that network because whenever a job is processed, a resource is forced to idle due to collaboration requirements. For example, resource 3 is forced to idle whenever a type 1 job is processed. Next, let us consider the parallel network in Figure \ref{f_bc3}, in which no resource is forced to idle during a service of a job. The only difference between the networks in Figures \ref{f_bc2} and \ref{f_bc3} is that the resource that is forced to idle in the former one helps the ongoing service process in the latter network. Observe that any scheduling policy implemented in the latter (former) network can also be implemented in the former (latter) network. Therefore, the networks in Figures \ref{f_bc2} and \ref{f_bc3} are path-wise identical to each other with respect to the total discounted holding cost. Specifically, for fixed stochastic primitives and scheduling policy, the holding costs incurred in those networks are the same. Because the network in Figure \ref{f_bc3} has hierarchical collaboration architecture, the proposed policy is asymptotically optimal in that network. In fact, the network in Figure \ref{f_bc3} is equivalent to a multi-class $GI/GI/1$ queuing system (see Example \ref{s_ex_1}). Consequently, the proposed policy associated with the network in Figure \ref{f_bc3} is also asymptotically optimal in the network in Figure \ref{f_bc2}. 

\begin{figure}[ht]
  \centering
  \subfloat[A parallel network with 4 job types and 3 multitasking resources.]{\label{f_cl2}\includegraphics[width=0.43\textwidth]{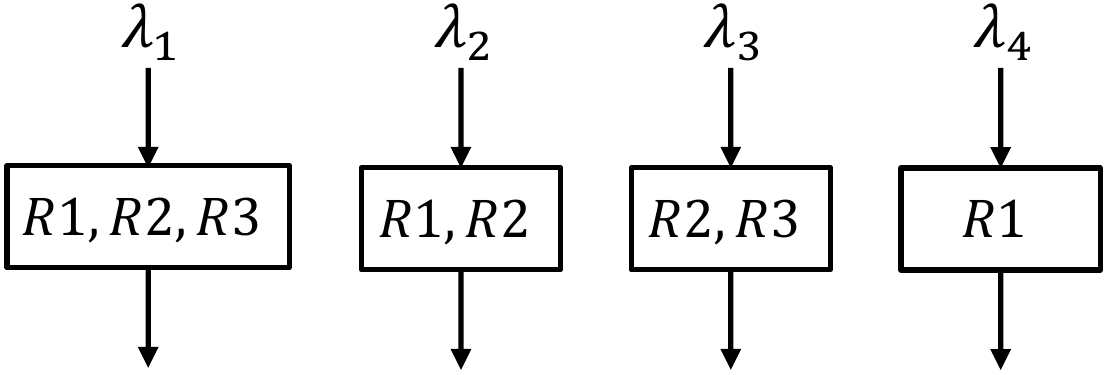}}\hspace{2cm}
  \subfloat[The graph associated with Figure \ref{f_cl2}.]{\label{f_cl2g}\includegraphics[width=0.2\textwidth]{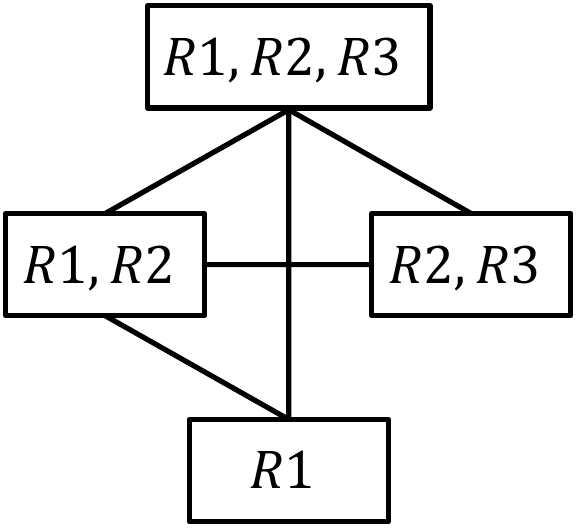}}\\
    \subfloat[The associated hypothetical network.]{\label{f_cl3}\includegraphics[width=0.45\textwidth]{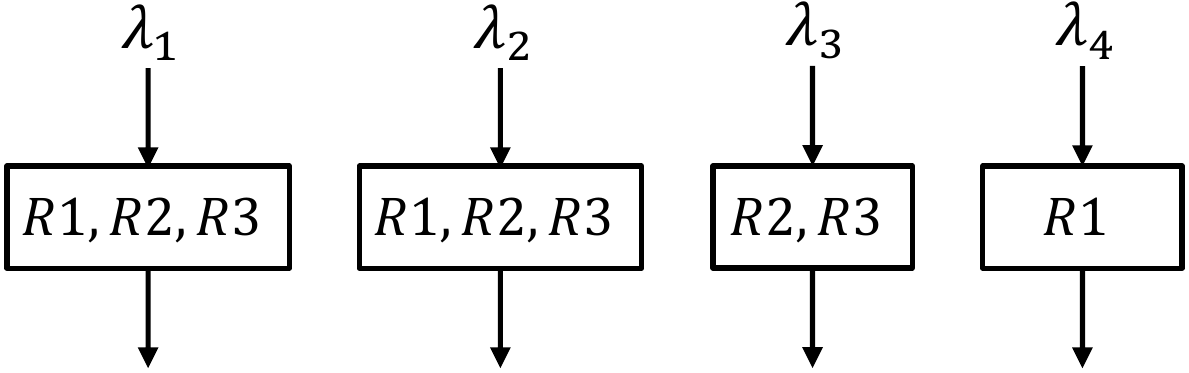}}\hspace{2cm}
  \subfloat[The graph associated with Figure \ref{f_cl3}.]{\label{f_cl3g}\includegraphics[width=0.24\textwidth]{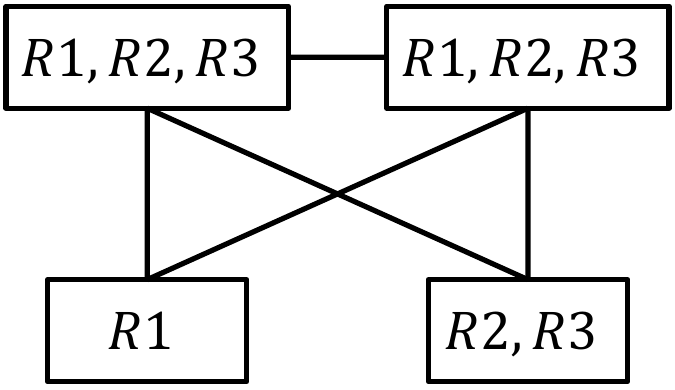}}\\
  \caption{A parallel network whose hypothetical network has hierarchical collaboration architecture.}
  \label{f_network_2}
\end{figure}

We can generalize the extension idea illustrated above for general parallel networks. The key idea is to construct a hypothetical parallel network such that the resources that are forced to idle during a service process in the original network help the associated service process in the hypothetical network. The formal construction of the hypothetical network is as follows. Let $\mP_j:=\{l\in\mJ: \mI_j\cap\mI_l = \emptyset\}$ for all $j\in\mJ$. Then, $\mP_j$ denotes the set of job types that can be processed in parallel with type $j$ jobs. Let 
\begin{equation*}
\mZ_j:=\left\{i\in\mI: i\notin\mI_j\text{ and $i\notin\mI_l$ for all $l\in\mP_j$}\right\},\qquad\forall j\in\mJ.
\end{equation*}
Then, $\mZ_j$ denotes the set of resources that are forced to idle during the service of a type $j$ job. Let us represent the hypothetical network with the $'$ symbol and construct it such that
\begin{equation*}
\mI'_1 := \mI_1\cup\mZ_1,\qquad \mI'_j := \mI_j\cup\left(\mZ_j\backslash \left( \mZ_{j-1}\cup\mZ_{j-2}\cup\ldots\cup\mZ_1\right) \right),\quad\forall j>1.
\end{equation*}
Hence, the resources that are forced to idle in the original network are no longer forced to idle in the hypothetical network. By construction, the only difference between the original and hypothetical networks is the resource-job type incidence matrix $A$. Furthermore, any scheduling policy that is admissible in the original (hypothetical) network is also admissible in the hypothetical (original) network with path-wise the same total discounted holding cost. Therefore, if the hypothetical network has hierarchical collaboration architecture, then the proposed policy associated with the hypothetical network is admissible and asymptotically optimal in the original network too. Figure \ref{f_network_2} presents such an example. Figure \ref{f_cl2} shows the original network that has a capacity loss (R3 is forced to idle whenever a type 2 job is processed). Figure \ref{f_cl3} shows the associated hypothetical network, which has hierarchical collaboration architecture (see Figure \ref{f_cl3g}).

\begin{figure}[ht]
  \centering
  \subfloat[A parallel network with 6 job types and 4 multitasking resources.]{\label{f_cl}\includegraphics[width=0.56\textwidth]{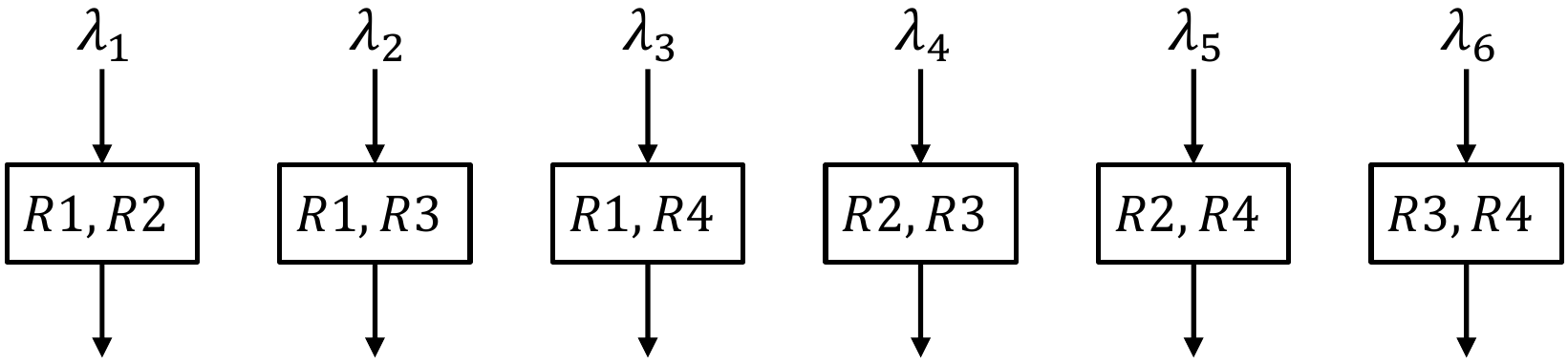}}\hspace{1cm}
  \subfloat[The graph associated with Figure \ref{f_cl}.]{\label{f_clg}\includegraphics[width=0.3\textwidth]{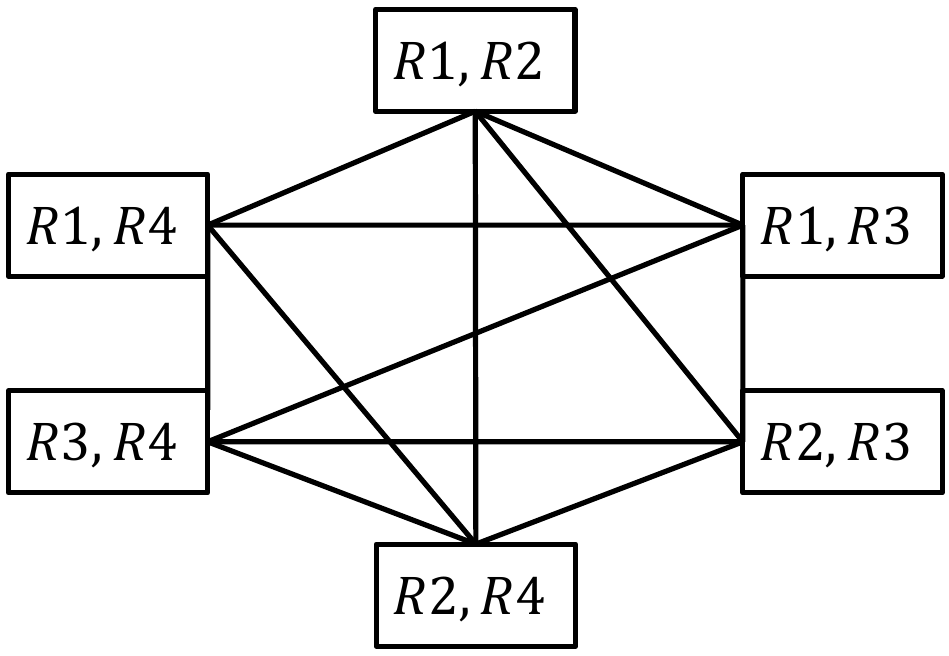}}
  \caption{An identical original and hypothetical network pair with capacity loss.}
  \label{f_network_3}
\end{figure}

Observe that if the original network has hierarchical collaboration architecture, then the associated hypothetical network is the same of the original one because $\bigcup_{j\in\mJ}\mZ_j=\emptyset$. However, the fact that $\bigcup_{j\in\mJ}\mZ_j=\emptyset$ in the original network does not necessarily imply that there is no capacity loss in the original network. For example, consider the network in Figure \ref{f_cl} in which $\bigcup_{j\in\mJ}\mZ_j=\emptyset$ (and so the original and the associated hypothetical networks are the same) but there is capacity loss. In order to intuitively understand the capacity loss, consider the state in which there are only type 1, type 2, and type 3 jobs in the system. In that state, exactly two resources will be forced to idle. Finally, Figure \ref{f_clg} shows that the network in Figure \ref{f_cl} does not have hierarchical collaboration architecture.

In summary, our results can be extended to a class of parallel networks with capacity loss. Specifically, if the hypothetical network associated with a parallel network has hierarchical collaboration architecture, then the proposed policy is admissible and asymptotically optimal in the original network.

\section{Concluding Remarks}\label{s_conc}

We consider scheduling control of parallel networks with resource collaboration and multitasking. Specifically, we consider parallel networks with hierarchical collaboration architecture, which do not have any capacity loss. We propose a dynamic prioritization policy which is asymptotically optimal in diffusion scale in the conventional heavy-traffic regime. Under the proposed policy, an LP is solved at discrete time epochs and the resource allocation decisions are done with respect to an index rule that depends on the optimal LP solution. Finally, we extend our results to a class of parallel networks with capacity loss.

\begin{figure}[htb]
\begin{center}
\includegraphics[width=0.5\textwidth]{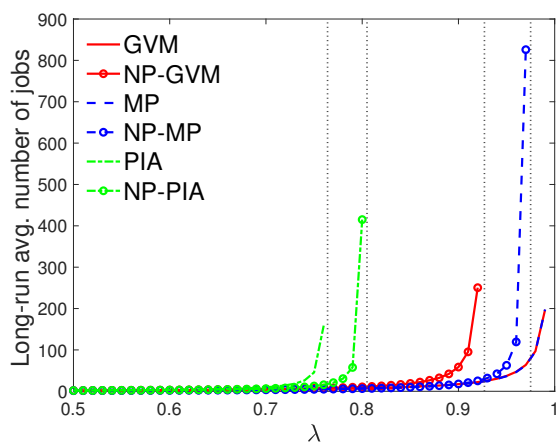}
\caption{(Color online) Comparison of the performances of the GVM policy, the MP policy, the PIA policy, and their nonpreemptive versions, namely NP-GVM, NP-MP, and NP-PIA, respectively.}\label{f3}
\end{center}
\end{figure}

Because nonpreemptive policies have more practical appeal than the preemptive ones have, an interesting future research topic is to propose well performing nonpreemptive control policies. It is not even clear how to construct a non-preemptive version of the proposed policy. In conventional queueing networks, the performance gap between a preemptive control policy and its nonpreemptive counterpart is negligible in the heavy-traffic regime. However, in networks with resource collaboration and multitasking, the aforementioned gap is not negligible in the heavy-traffic regime (see Section 3.2 of \cite{gur18} for details). For example, Figure \ref{f3} shows that there is a capacity loss under the nonpreemptive GVM (NP-GVM) policy because the NP-GVM policy is not work-conserving. Figure \ref{f3} shows that there is a capacity loss under the nonpreemptive MP (NP-MP) policy. If preemption is not allowed, the feasible set of resource allocations is reduced which, in turn, decreases the capacity (see Section 8 of \cite{dai05} for details). Interestingly, Figure \ref{f3} shows that the system capacity under the nonpreemptive PIA (NP-PIA) policy is higher than the one under the PIA policy. This is because the nonpreemptive structure of the NP-PIA policy prevents the blow-up of the buffer 1 to some extent. Consequently, deriving non-preemptive control policies that perform well is not trivial and requires further research.

Parallel networks whose associated hypothetical network not having hierarchical collaboration architecture may have capacity loss (see for example the network in Figure \ref{f_cl}). It is not clear how to derive control policies that perform well for those networks. Because some resources cannot be fully utilized in those networks, the BCP \eqref{eq_bcp} no longer provides an achievable lower bound on the performance of admissible policies. Nevertheless, the aforementioned problem is valuable for further research. Similarly, control of networks with general topological structure, that is, control of networks beyond the parallel ones, is a relevant area for future research. However, even the first-order analysis of those networks is very challenging. For example, \cite{bo19} prove that not only the exact computation of the system capacity of general networks is a strongly NP-hard problem but also approximating the system capacity accurately is an NP-hard problem.


%
%
%



\bibliographystyle{informs2014} 
\bibliography{collab}

\newpage
\pagenumbering{arabic}
\numberwithin{theorem}{section}
\numberwithin{lemma}{section}
\numberwithin{proposition}{section}
\numberwithin{assumption}{section}
\numberwithin{remark}{section}
\numberwithin{table}{section}
\numberwithin{figure}{section}
\setcounter{section}{0}
\renewcommand{\thesection}{EC.\arabic{section}}
\setcounter{equation}{0}
\renewcommand{\theequation}{EC.\arabic{equation}}

\begin{center}
\textbf{ELECTRONIC COMPANION}
\end{center}

This electronic companion is associated with the manuscript titled ``On the Optimal Control of Parallel Processing Networks with Resource Collaboration and Multitasking''. Section \ref{l_proof} presents the proofs of Lemmas \ref{l_diffusion}, \ref{l_alg}, and \ref{l_review}. Section \ref{t_lb_proof} presents the proof of Theorem \ref{t_lb}. Section \ref{t_ao_proof} presents the proof of Theorem \ref{t_ao}. Section \ref{p_wc_proof} presents the proof of Proposition \ref{p_wc}, which is used in the proof of Theorem \ref{t_ao}. Section \ref{l_good_set_proof} presents the proof of Lemma \ref{l_good_set}, which is used in the proof of Proposition \ref{p_wc}. Finally, Section \ref{l_good_4_proof} presents some proofs used in the proof of Lemma \ref{l_good_set}.

\section{Lemma Proofs}\label{l_proof}

In this section we prove Lemmas \ref{l_diffusion}, \ref{l_alg}, and \ref{l_review} in Sections \ref{l_diffusion_proof}, \ref{l_alg_proof}, and \ref{l_review_proof}, respectively.

\subsection{Proof of Lemma \ref{l_diffusion}}\label{l_diffusion_proof}

The proof is pretty standard and follows from \eqref{eq_min}, functional strong law of large numbers (FSLLN, see Theorem 5.10 of \cite{che01}), functional central limit theorem (FCLT, see Theorems 4.3.5 and 7.3.2 of \cite{whi02}), joint convergence when one limit is deterministic (see Theorem 11.4.5 of \cite{whi02}), continuous mapping theorem (see Theorem 3.4.4 of \cite{whi02}), and the continuity of the multidimensional reflection map (see Theorem 14.2.5 of \cite{whi02}). For a similar proof, see the proof of Proposition 2 of \cite{ozk19}.

Next, we will prove that the covariance matrix $\Sigma$ defined in \eqref{eq_cov} is positive definite, which is a necessary condition for the limiting process in Lemma \ref{l_diffusion} to be an SRBM (see Definition 3.1 of \cite{wil98b} for details). Let $z=\left( z_i,i\in\mI^H\right)$ be an arbitrary nonzero vector in $\R^{a_1}$, the superscript $'$ denote the transpose of a matrix or vector, and $c_j := \la_j\left( \beta_j^2 + \sigma_{j}^2\right) /\mu_{j}^2$. After some algebra and because $c_j>0$ for all $j\in\mJ^H$, we have
\begin{equation*}
z'\Sigma z = \sum_{i\in\mI^H} \sum_{k\in\mI^H} \sum_{j\in\mJ_i\cap\mJ_k} c_j z_iz_k = \sum_{j\in\mJ^H} c_j \left(\sum_{i\in\mI_j\cap \mI^H} z_i \right)^2 \geq 0, 
\end{equation*}
which implies that $\Sigma$ is positive semi-definite. 

By Assumption \ref{a_h}, if there is a job type that is processed by a single resource, then that job type is a leaf node in the graph associated with the network structure. Recall that, without loss of generality, we assume that for all $i,k\in\mI$ such that $i\neq k$, both $\mJ_i\not\subset\mJ_k$ and $\mJ_k\not\subset\mJ_i$. This assumption together with Assumption \ref{a_h} implies that each leaf node is processed by a single resource. Otherwise, if there exists a leaf node that is processed in collaboration with multiple resources, then those resources process exactly the same job types. Therefore, a job type is a leaf node in the graph associated with the network structure if and only if it is processed by a single resource. Furthermore, the aforementioned result together with Assumption \ref{a_h} imply that for each resource, there exists a job type that is processed by only that resource. In other words, for each $i\in\mI$, there exists a $j\in\mJ_i$ such that $j\notin\mJ_k$ for all $k\in\mI\backslash\{i\}$. Otherwise, if there exists an $i\in\mI$ which does not have a job type that is processed by only resource $i$, then resource $i$ does not process a leaf node. This further implies that there exists a resource $k\in\mI$ such that $\mI_i\subset\mI_k$, which is a contradiction. Let $\ga_i\in\mJ^H$ denote the job type that is processed by only resource $i$ for all $i\in\mI^H$. Then,
\begin{equation*}
z'\Sigma z = \sum_{i\in\mI^H} \sum_{k\in\mI^H} \sum_{j\in\mJ_i\cap\mJ_k} c_j z_iz_k = \sum_{j\in\mJ^H} c_j \left(\sum_{i\in\mI_j\cap \mI^H} z_i \right)^2 \geq \sum_{i\in\mI^H} c_{\ga_i} z_i^2 >0, 
\end{equation*}
which implies that $\Sigma$ is positive definite.

\subsection{Proof of Lemma \ref{l_alg}}\label{l_alg_proof}

Because $\mX(s)$ has finite number of elements by definition and $\mX(s)\neq \emptyset$ by Assumption \ref{a_h}, the set in \eqref{eq_set} has also finite number of elements. Let $x(s)$ be an arbitrary resource allocation vector from the set in \eqref{eq_set} and let $x^A(s)$ denote the resource allocation vector given by the algorithm in Definition \ref{d_alg}. Let
\begin{equation*}
\mR^r(x(s),t):=\left\{i\in\mI^H: \exists j\in\mJ^H \text{ s.t. } x_j(s)=1, Q_j^r(s)> \lceil q_j^{*,r}(t) \rceil, i\in\mI_j \right\},
\end{equation*}
that is, $\mR^r(x(s),t)$ is the set of resources in heavy traffic that are processing job types whose number in system at time $s$ is greater than the optimal LP \eqref{eq_lp} solution associated with time $t$ under $x(s)$. Suppose that there exists a $k\in\mR^r(x(s),t)\setminus\mR^r\big(x^A(s),t \big)$. There must exist an $l\in\mJ^H$ such that $x_l(s)=1$, $k\in\mI_l$, $Q_l^r(s)> \lceil q_l^{*,r}(t) \rceil$, and $x^A_l(s)=0$. Let $\mS_l:=\{j\in\mJ: \mI_l\subseteq\mI_j\}$. By Assumption \ref{a_h} and the construction of the Steps 1-4 of the algorithm in Definition \ref{d_alg}, it must be the case that $Q_j^r(s)\leq \lceil q_j^{*,r}(t) \rceil$ for all $j\in\mS_l$, otherwise $x^A_j(s)$ would be equal to 1 for some $j\in \mS_l$ such that $Q_j^r(s)> \lceil q_j^{*,r}(t) \rceil$ implying $k\in\mI_l\subseteq\mI_j\in\mR^r\big(x^A(s),t \big)$. However, this is a contradiction because $Q_l^r(s)> \lceil q_l^{*,r}(t) \rceil$ by assumption. Therefore, $\mR^r(x(s),t)\subseteq\mR^r\big(x^A(s),t \big)$ implying $p_1\big(x^A(s),t\big)\geq p_1(x(s),t)$.

Second, by considering the definition of the partial order $\succeq$, suppose that $p_1\big(x^A(s),t\big) = p_1(x(s),t)$ and $p_2\big(x^A(s)\big)<p_2(x(s))$. Because $p_1\big(x^A(s),t\big) = p_1(x(s),t)$ and $\mR^r(x(s),t)\subseteq\mR^r\big(x^A(s),t \big)$, we have $\mR^r(x(s),t)=\mR^r\big(x^A(s),t \big)$. Furthermore, because $p_2\big(x^A(s)\big)<p_2(x(s))$, there must exist an $k\in\mI$ such that $x_l(s)=1$ for some $l\in\mJ_k$ but $x^A_j(s)=0$ for all $j\in\mJ_k$. By Assumption \ref{a_h} and the construction of the Steps 1-4 of the algorithm in Definition \ref{d_alg} and because $x^A_j(s)=0$ for all $j\in\mJ_k$, it must be the case that $Q_j^r(s)\leq \lceil q_j^{*,r}(t) \rceil$ for all $j\in\mJ_k$. Because $\mR^r(x(s),t)=\mR^r\big(x^A(s),t \big)$ and $x_l(s)=1$ for some $l\in\mJ_k$, by Assumption \ref{a_h}, and by the construction of the Steps 5-9 of the algorithm in Definition \ref{d_alg}, it must be the case that $x^A_j(s)=1$ for some $j\in\mS_l$, which is a contradiction because $x^A_j(s)=0$ for all $j\in\mJ_k$. Therefore, $p_2\big(x^A(s)\big)\geq p_2(x(s))$, which completes the proof.

Finally, due to sorting, the computational complexity of the algorithm in Definition \ref{d_alg} is $O(J\ln J)$.

\subsection{Proof of Lemma \ref{l_review}}\label{l_review_proof}

Consider an arbitrary $i\in\mI^H$. By \eqref{eq_we}, we have
\begin{equation}\label{eq_review_p_1}
\sum_{j\in\mJ_i^{>,r}(t)} \frac{Q_j^r(t)-q_j^{*,r}(t)}{\mu_j}= \sum_{j\in\mJ_i^{\leq,r}(t)} \frac{q_j^{*,r}(t)-Q_j^r(t)}{\mu_j},\qquad\forall i\in\mI^H,\;t\in\R_+.
\end{equation}
By Assumption \ref{a_regime} Part 2, we have $\sum_{j\in\mJ_i}\la_j/\mu_j =1$ for all $i\in\mI^H$. Then, by \eqref{eq_review_p_1}, the RHS of \eqref{eq_review_1} is equal to
\begin{align}
&\sum_{j\in\mJ_i^{\leq,r}} \frac{q_j^{*,r}(t)-Q_j^r(t)}{\mu_j} +\breve{L}_i^r(t)\left(1-\sum_{j\in\mJ_i^{\leq,r}(t)} \frac{\la_j}{\mu_j}\right)  + \sum_{j\in\mJ_i^{\leq,r}(t)} \frac{\left(\la_j \breve{L}_i^r(t) - q_j^{*,r}(t)+Q_j^r(t)\right)^+}{\mu_j}\nonumber\\
&\hspace{3cm}= \breve{L}_i^r(t)+ \sum_{j\in\mJ_i^{\leq,r}(t)} \frac{q_j^{*,r}(t)-Q_j^r(t)-\la_j\breve{L}_i^r(t)+ \left(\la_j \breve{L}_i^r(t) - q_j^{*,r}(t)+Q_j^r(t)\right)^+}{\mu_j}\nonumber\\
&\hspace{3cm}= \breve{L}_i^r(t)+ \sum_{j\in\mJ_i^{\leq,r}(t)} \frac{\left(q_j^{*,r}(t)-Q_j^r(t)-\la_j\breve{L}_i^r(t)\right)^+}{\mu_j}.\label{eq_l_review}
\end{align}
The term in \eqref{eq_l_review} is equal to the LHS of \eqref{eq_review_1} if and only if  $q_j^{*,r}(t)-Q_j^r(t)-\la_j\breve{L}_i^r(t)\leq 0$ for all $j\in\mJ_i^{\leq,r}(t)$ which holds if and only if \eqref{eq_review} holds.

\section{Proof of Theorem \ref{t_lb}}\label{t_lb_proof}

Let us fix an arbitrary sequence of admissible policies denoted by $\pi=\{\pi^r,r\in\N_+\}$. By \eqref{eq_queue}, \eqref{eq_workload}, and \eqref{eq_diffusion},
\begin{equation*}
\sum_{j\in\mJ_i} \frac{\hat{Q}_j^{\pi,r}(t)}{\mu_j} = \hat{W}_i^{\pi,r}(t)\quad\forall i\in\mI^H,\hspace{2cm}\hat{Q}_j^{\pi,r}(t)\geq 0\quad\forall j\in\mJ^H.
\end{equation*}
Therefore, $\big(\hat{Q}_j^{\pi,r}(t), j\in\mJ^H\big)$ is a feasible point of the LP \eqref{eq_lp} under the RHS parameter $\big(\hat{W}_i^{\pi,r}(t),i\in\mI^H\big)$ for all $t\in\R_+$. Under all sample paths and for all $t\in\R_+$, we have
\begin{equation}\label{eq_lb_proof_0}
\sum_{j\in\mJ^H} h_j \hat{Q}_j^{\pi,r}(t) \geq  z\Big(\hat{W}_i^{\pi,r}(t),i\in\mI^H\Big).
\end{equation}

Let $Z$ be a mapping from $\D^{a_1}$ such that $Z(f)(t):=z(f(t))$ for all $f\in\D^{a_1}$ and $t\in\R_+$. Then, $Z$ is the process version of $z$. Because $z$ is Lipschitz continuous (see Lemma \ref{l_lp} Part 1), $Z$ maps the functions from $\D^{a_1}$ to $\D$, that is, $Z:\D^{a_1}\rightarrow\D$. Let $d(\cdot)$ denote the Skorokhod distance (see Equation (12.13) of \cite{bil99}). For arbitrary $X,Y\in\D^{a_1}$, because $z$ is Lipschitz continuous, one can see that $d(Z(X),Z(Y))\leq (C_1\vee 1)d(X,Y)$. Therefore, $Z$ is Lipschitz continuous. Then, by Lemma \ref{l_diffusion} and continuous mapping theorem, we have
\begin{equation}\label{eq_lb_proof_1}
Z\left(\hat{W}_i^{\pi,r},i\in\mI^H\right)\Rightarrow Z\big(\ot{W}^*\big).
\end{equation}
Next, we invoke the Skorokhod representation theorem (see Theorem 3.2.2 of \cite{whi02}) on the processes in \eqref{eq_lb_proof_1} such that 
\begin{equation}\label{eq_lb_proof_2}
\breve{Z}^{\pi,r}  \overset{d}{=} Z\left(\hat{W}_i^{\pi,r},i\in\mI^H\right),\;\;\forall r\in\N_+,\qquad \breve{Z}  \overset{d}{=} Z\big(\ot{W}^*\big),\qquad\breve{Z}^{\pi,r} \xrightarrow{a.s.} \breve{Z}\;\;\text{u.o.c.}
\end{equation}

Finally,
\begin{align}
\liminf_{r\rightarrow\infty} \sum_{j\in\mJ} h_j  \E\left[ \int_0^\infty \e^{-\delta t} \hat{Q}_j^{\pi,r}(t) \dr t\right] &= \liminf_{r\rightarrow\infty} \E\left[ \int_0^\infty \e^{-\delta t} \left(\sum_{j\in\mJ} h_j\hat{Q}_j^{\pi,r}(t)\right) \dr t\right]\nonumber\\
&\geq \liminf_{r\rightarrow\infty} \E\left[ \int_0^\infty \e^{-\delta t} z\Big(\hat{W}_i^{\pi,r}(t),i\in\mI^H\Big) \dr t\right]\label{eq_lb_proof_4}\\
&= \liminf_{r\rightarrow\infty} \E\left[ \int_0^\infty \e^{-\delta t} Z\Big(\hat{W}_i^{\pi,r},i\in\mI^H\Big)(t) \dr t\right]\label{eq_lb_proof_5}\\
&= \liminf_{r\rightarrow\infty} \E\left[ \int_0^\infty \e^{-\delta t} \breve{Z}^{\pi,r}(t) \dr t\right]\label{eq_lb_proof_6}\\
&\geq  \E\left[ \int_0^\infty \e^{-\delta t} \liminf_{r\rightarrow\infty}\breve{Z}^{\pi,r}(t) \dr t\right]\label{eq_lb_proof_7}\\
& = \E\left[ \int_0^\infty \e^{-\delta t} \breve{Z}(t) \dr t\right]\label{eq_lb_proof_8} \\
&=\E\left[ \int_0^\infty \e^{-\delta t} z\big(\ot{W}^*(t)\big) \dr t\right],\label{eq_lb_proof_9}
\end{align}
where \eqref{eq_lb_proof_4} is by \eqref{eq_lb_proof_0}, \eqref{eq_lb_proof_5} is by the definition of the mapping $Z$, \eqref{eq_lb_proof_6} is by \eqref{eq_lb_proof_2}, \eqref{eq_lb_proof_7} is by Fatou's lemma, \eqref{eq_lb_proof_8} is by \eqref{eq_lb_proof_2}, and \eqref{eq_lb_proof_9} is by \eqref{eq_lb_proof_2}.

\section{Proof of Theorem \ref{t_ao}}\label{t_ao_proof}

The following proposition, whose proof is presented in E-companion \ref{p_wc_proof}, will be very helpful in the proof of Theorem \ref{t_ao}.

\begin{proposition}\label{p_wc}
Under the proposed policy (see Definition \ref{d_pp}), we have
\begin{equation*}
\sum_{j\in\mJ} h_j\hat{Q}_j^r\Rightarrow Z\big(\ot{W}^*\big).
\end{equation*}
\end{proposition}
We will prove Theorem \ref{t_ao} by Proposition \ref{p_wc} and proving a uniform integrability result (see for example the proof of Theorem 5.3 in \cite{bel01}). 

We invoke the Skorokhod representation theorem to obtain equivalent distributional representations of the processes in Proposition \ref{p_wc} (for which we use the same symbols and call ``Skorokhod represented versions'') such that the Skorokhod represented versions of the processes are defined in the same probability space and satisfy
\begin{equation*}
\sum_{j\in\mJ} h_j\hat{Q}_j^r\xrightarrow{a.s.} Z\big(\ot{W}^*\big)\quad\text{u.o.c.}
\end{equation*}
Then,
\begin{equation*}
\sum_{j\in\mJ} h_j\hat{Q}_j^r\rightarrow Z\big(\ot{W}^*\big)\qquad\text{$(m\times\pr)-$almost everywhere},
\end{equation*}
where $\dr m := \delta\e^{-\delta t}\dr t$ on $(\R_+,\mB(\R_+))$ and $\mB(\R_+)$ is the Borel $\sigma$-algebra on $\R_+$. Because $(\R_+\times\Omega,\mB(\R_+)\times\mF, m\times\pr)$ is a probability space, the condition
\begin{equation}\label{ui_cond_0}
\limsup_{r\rightarrow \infty} \E\left[\int_0^{\infty} \e^{-\delta t} \left(\sum_{j\in\mJ} h_j\hat{Q}_j^r(t)\right)^2 \dr t \right] <\infty,
\end{equation}
which implies the uniform integrability of $\sum_{j\in\mJ} h_j\hat{Q}_j^r$ with respect to expectation under $m\times\pr$, gives us the desired result. In order to prove \eqref{ui_cond_0}, it is enough to show that
\begin{equation}\label{ui_cond_1}
\limsup_{r\rightarrow \infty}\int_0^{\infty} \e^{-\delta t} \E\left[ \left(\hat{Q}_j^{r}(t)\right)^2\right]  \dr t <\infty,\qquad\forall j\in\mJ.
\end{equation}
Hence, we will prove \eqref{ui_cond_1} henceforth. 

By Assumption \ref{a_regime}, there exists an $r_1\geq r_0$ such that if $r\geq r_1$, then
\begin{align}
&0.5\la_j\leq \la_j^r \leq 1+\la_j,\quad\forall j\in\mJ,\label{eq_rdef_1}\\
&\rho_i^r \leq 1,\quad\forall i\in\mI^L,\qquad r\left| 1- \rho_i^r\right| \leq 1+|\theta_i|,\quad\forall i\in\mI^H,\label{eq_rdef_2}
\end{align}
where $r_0$ is defined in \eqref{eq_ini1}. We will use the following result in the rest of the proof.
\begin{lemma}\label{ui_renewal}
\textbf{(Equation (172) of \cite{bel01})} If $r\geq r_1$, then for all $j\in\mJ$ and $t\in\R_+$, we have
\begin{align*}
&\E\left[ \sup_{0\leq s\leq t} \left( \hat{E}_j^r(s)\right)^2 \right] \leq 18\left(1+ 4\beta_j^2 \left((1+\la_j) t + \beta_j^2 +2\right)\right),\\
&\E\left[ \sup_{0\leq s\leq t} \left( \hat{S}_j^r(s)\right)^2 \right] \leq 18\left(1+ 4\sigma_j^2 \left(\mu_j t + \sigma_j^2 +2\right)\right).
\end{align*}
\end{lemma}
Notice that we present a slightly different version of (172) of \cite{bel01} such that there is no $r\geq r_1$  condition and there is $\la_j^rt$ instead of the term $(1+\la_j)t$ in (172) of \cite{bel01}.

First, let us consider an arbitrary $j\in\mJ^H$. Let $i\in\mI_j$ be an arbitrary resource in heavy traffic that process type $j$ jobs. By \eqref{eq_workload}, \eqref{eq_workload_2}, and \eqref{eq_workload_3}, if $r\geq r_1$, then for all $t\in\R_+$, we have
\begin{align}
\frac{\hat{Q}_j^r(t)}{\mu_j}\leq \hat{W}^r_i(t) &= \hat{X}_i^r(t) + \hat{I}_i^r(t)= \hat{X}_i^r(t) + \Psi(\hat{X}_i^r)(t)\label{eq_ui_proof_1}\\
&=  \sum_{l\in\mJ_i} \frac{1}{\mu_{l}}\left(\hat{Q}_l^r(0)  + \hat{E}_l^r(t) -  \hat{S}_{l}^r\circ\bar{T}_{l}^r(t)\right) + r\left(\rho_i^r-1\right)t\nonumber\\
&\hspace{1cm} + \sup_{0\leq s\leq t} \left( \sum_{l\in\mJ_i} \frac{1}{\mu_{l}}\left(-\hat{Q}_l^r(0) -\hat{E}_l^r(s) +\hat{S}_{l}^r\circ\bar{T}_{l}^r(s)\right) + r\left( 1-\rho_i^r\right)s \right)^+ \nonumber\\
&\leq \sum_{l\in\mJ_i} \left[\frac{2}{\mu_l}\left(\hat{Q}_l^r(0) + \sup_{0\leq s\leq t} \left|\hat{E}_l^r(s) \right|+ \sup_{0\leq s\leq t} \left|\hat{S}_l^r(s) \right|\right)\right]  + r\left| 1- \rho_i^r\right|t \label{eq_ui_proof_2}\\
&\leq \sum_{l\in\mJ_i} \left[\frac{2}{\mu_l}\left(\hat{Q}_l^r(0) + \sup_{0\leq s\leq t} \left|\hat{E}_l^r(s) \right|+ \sup_{0\leq s\leq t} \left|\hat{S}_l^r(s) \right|\right)\right] +  (1+| \theta_i |)t, \label{eq_ui_proof_3}
\end{align}
where the last equality in \eqref{eq_ui_proof_1} is by \eqref{eq_min}, \eqref{eq_ui_proof_2} is by the fact that  $\bar{T}_l^r(t) \leq t$ for all $l\in\mJ$ and $t\in\R_+$, and \eqref{eq_ui_proof_3} holds because $r\geq r_1$ (see \eqref{eq_rdef_2}). Therefore, if $r\geq r_1$, for all $t\in\R_+$, we have
\begin{equation}\label{eq_ui_proof_4}
\sup_{0\leq s\leq t} \left(\hat{Q}^r_j(s)\right)^2 \leq 4(3J+1)\mu_j^2\left[ \sum_{l\in\mJ_i}\frac{1}{\mu_l^2}\left( \left(\hat{Q}_l^r(0)\right)^2 + \sup_{0\leq s\leq t} \left|\hat{E}_l^r(s) \right|^2 + \sup_{0\leq s\leq t} \left|\hat{S}_l^r(s) \right|^2\right) + (1+|\theta_i|) ^2 t^2 \right].
\end{equation}
Finally, \eqref{eq_ini1}, Lemma \ref{ui_renewal}, and \eqref{eq_ui_proof_4} imply that
\begin{equation}\label{eq_ui_proof_5}
\E\left[  \sup_{0\leq s\leq t} \left(\hat{Q}_j^{r}(s)\right)^2\right] \leq C_3+C_4t +C_5t^2,\qquad\forall t\in\R_+, r\geq r_1,
\end{equation}
where
\begin{align*}
&\bar{\la}:= \max_{j\in\mJ}\la_j,\qquad \underline{\la}:=\min_{j\in\mJ}\la_j,\qquad\bar{\mu}:= \max_{j\in\mJ}\mu_j,\qquad \underline{\mu}:=\min_{j\in\mJ}\mu_j,\\
&C_3 :=  2^43^2J(3J+1) \frac{\bar{\mu}^2}{\underline{\mu}^2}\left(1+4\left(\max_{j\in\mJ} \beta_j^4 \vee \max_{j\in\mJ} \sigma_j^4 \right) +8\left(\max_{j\in\mJ} \beta_j^2 \vee \max_{j\in\mJ} \sigma_j^2 \right) + C_0 \right),\\
&C_4:= 2^63^2J(3J+1)\frac{\bar{\mu}^2}{\underline{\mu}^2}\left(\max_{j\in\mJ} \beta_j^2 \vee \max_{j\in\mJ} \sigma_j^2 \right) (1+\bar{\la}\vee\bar{\mu}),\\
&C_5:= 4(3J+1)\bar{\mu}^2\left(1+\max_{i\in\mI^H}|\theta_i|\right) ^2.
\end{align*}

Observe that \eqref{eq_ui_proof_5} implies that \eqref{ui_cond_1} holds for all $j\in\mJ^H$.

Next, let us consider an arbitrary $j\in\mJ^L$. Observe that all resources in the set $\mI_j$ are in light traffic by definition of $\mJ^L$. Because the resources in the set $\mI_j$ are not necessarily utilized in a work-conserving fashion under the proposed policy, we do not necessarily have the equality $\hat{I}_i^r = \Psi(\hat{X}_i^r)$. Hence, the last equality in \eqref{eq_ui_proof_1} does not necessarily hold and so the proof technique that we use above to prove the uniform integrability of the number of jobs in the buffers processed by resources in heavy traffic does not work in this case. Instead, we will use coupling arguments to derive a uniformly integrable upper bound on $Q_j^r$, which will imply uniform integrability of $Q_j^r$.

We construct a hypothetical system, denoted by the superscript $^{(1)}$, in which the stochastic primitives are the same of the ones in the original system. Let us consider an arbitrary $i\in\mI_j$. Suppose that, in the hypothetical system, resource $i$ operates in a work-conserving fashion; type $j$ jobs are processed at time $t\in\R_+$ if and only if $Q_j^{(1),r}(t)>0$ and $Q_l^{(1),r}(t)=0$ for all $l\in\mJ_{i}\backslash\{j\}$, that is, type $j$ jobs receive the least amount of priority from resource $i$; and whenever resource $i$ needs to collaborate with some other resources to process a job, then those resources are assigned to collaborate with resource $i$, that is, there is no collaboration constraint for resource $i$. Observe that, by Assumption \ref{a_h} and because the proposed policy utilizes maximum number of resources in light traffic (recall the index $p_2(x(s))$ in \eqref{eq_index}), in the original system and under the proposed policy, type $j$ jobs are processed at time $t\in\R_+$ if $Q_j^r(t)>0$ and $Q_l^r(t)=0$ for all $l\in\mJ_{i}\backslash\{j\}$. Therefore, we have the following path-wise upper bound: 
\begin{equation}\label{eq_coupling}
Q_j^r \leq Q_j^{(1),r}.
\end{equation}
By \eqref{eq_coupling} and because resource $i$ is utilized in a work-conserving fashion in the hypothetical system, similar to how we derive \eqref{eq_ui_proof_5}, we can prove that
\begin{equation}\label{eq_ui_proof_6}
\E\left[  \sup_{0\leq s\leq t} \left(\hat{Q}_j^{(1),r}(s)\right)^2\right] \leq C_3+C_4t,\qquad\forall t\in\R_+,\;r\geq r_1.
\end{equation}
The only difference between the derivations of \eqref{eq_ui_proof_5} and \eqref{eq_ui_proof_6} is that the term $r\left|1-\rho_i^r\right|t$ that we see in \eqref{eq_ui_proof_2} does not appear in the derivation of \eqref{eq_ui_proof_6} when $r\geq r_1$ (see \eqref{eq_rdef_2}). Observe that \eqref{eq_coupling} and \eqref{eq_ui_proof_6} imply that \eqref{ui_cond_1} holds for all $j\in\mJ^L$.

\section{Proof of Proposition \ref{p_wc}}\label{p_wc_proof}

By Lemma \ref{l_diffusion} and continuous mapping theorem, we have
\begin{equation}\label{eq_pwc_6}
Z\left(\hat{W}_i^r,i\in\mI^H\right) \Rightarrow Z\big(\ot{W}^*\big).
\end{equation}
By \eqref{eq_pwc_6} and convergence-together theorem (see Theorem 11.4.7 of \cite{whi02}), in order to prove Proposition  \ref{p_wc}, it is enough to prove the following convergence result:
\begin{equation}\label{eq_pwc_7}
\lim_{r\rightarrow\infty} \pr\left( \bigg\Vert \sum_{j\in\mJ} h_j\hat{Q}_j^r - Z\left(\hat{W}_i^r,i\in\mI^H\right) \bigg\Vert_T >\ep \right) = 0, \qquad\forall \ep,T>0.
\end{equation}
We will prove \eqref{eq_pwc_7} henceforth. 

Let us fix arbitrary $\ep,T>0$. Let $\hat{q}_j^{*,r}(t):= r^{-1}q_j^{*,r} (r^2t)$ denote the diffusion scaled version of the optimal solution process associated with the LP \eqref{eq_lp} for all $t\in\R_+$, $j\in\mJ^H$, and $r\in\N_+$. Then, the probability in \eqref{eq_pwc_7} is equal to
\begin{align}
&\pr\left( \bigg\Vert \sum_{j\in\mJ^H} h_j\left(\hat{Q}_j^r - \hat{q}_j^{*,r}\right) + \sum_{j\in\mJ^L} h_j\hat{Q}_j^r \bigg\Vert_T >\ep \right)\nonumber\\
&\hspace{2cm}\leq \sum_{j\in\mJ^L}\pr\left( \left\Vert \hat{Q}_j^r \right\Vert_T >\frac{\ep}{J\max_{j\in\mJ} h_j} \right) + \sum_{j\in\mJ^H}\pr\left( \left\Vert \hat{Q}_j^r - \hat{q}_j^{*,r}\right\Vert_T >\frac{\ep}{J\max_{j\in\mJ} h_j} \right)\nonumber\\
&\hspace{2cm}= \sum_{j\in\mJ^L}\pr\left( \left\Vert \hat{Q}_j^r \right\Vert_T > \ep_1 \right) + \sum_{j\in\mJ^H}\pr\left( \left\Vert \hat{Q}_j^r - \hat{q}_j^{*,r}\right\Vert_T > \ep_1 \right),\label{eq_pwc_8}
\end{align}
where $\ep_1:=\ep/(J\max_{j\in\mJ}h_j)$. Therefore, in order to prove \eqref{eq_pwc_7}, it is enough to prove that the sum in \eqref{eq_pwc_8} converges to 0. In other words, we need to show that the diffusion scaled number of jobs that are processed by only resources in light traffic should be close to 0 and the diffusion scaled number of jobs that are processed by resources in heavy traffic should be close to the diffusion scaled optimal solution process associated with the LP \eqref{eq_lp}.

Let us first consider the first sum in \eqref{eq_pwc_8}. Let us fix arbitrary $j\in\mJ^L$ and $i\in\mI_j$. We will use coupling arguments again. Specifically, we will use the hypothetical system defined in Section \ref{t_ao_proof} and the path-wise inequality in \eqref{eq_coupling}. Then,
\begin{equation*}
\pr\left(\left\Vert \hat{Q}_j^r \right\Vert_T > \ep_1\right) \leq \pr\left(\left\Vert \hat{Q}_j^{(1),r} \right\Vert_T > \ep_1\right)\leq \pr\left( \left\Vert \hat{W}_i^{(1),r} \right\Vert_T > \frac{\ep_1}{\bar{\mu}}\right) \rightarrow 0,
\end{equation*}
where the first inequality is by \eqref{eq_coupling}, the second inequality is by \eqref{eq_workload}, and the convergence is because resource $i$ is in light traffic and utilized in a work-conserving fashion in the hypothetical system. Because we choose $j\in\mJ^L$ arbitrarily, the first sum in \eqref{eq_pwc_8} converges to 0.

Next, let us consider the second sum in \eqref{eq_pwc_8}. Let $\tau_{n}^r:\Omega\rightarrow\R_+\cup\{\infty\}$ denote the start time of the $n$th review period (Step 2 or 3) under the proposed policy for all $n,r\in\N_+$. For completeness, if $\tau_{n}^r(\omega)=\infty$ for some $n,r\in\N_+$ and $\omega\in\Omega$, then $\tau_{m}^r(\omega):=\infty$ for all $m>n$. Then, $\tau_{1}^r(\omega)=0$ and $\tau_{n+1}^r(\omega) \geq\tau_{n}^r(\omega)$ for all $n,r\in\N_+$ and $\omega\in\Omega$. 

By construction of the proposed policy (see Definition \ref{d_pp}), there exists an external arrival or service completion in each review period. Let $N^r:=1+\ru{2J(1+\bar{\la}\vee\bar{\mu})r^2T}$ for all $r\in\N_+$. Then,
\begin{align}
\pr\left( \tau_{N^r}^r \leq r^2T\right) &\leq \pr\left( \sum_{j\in\mJ} \left(E^r_j(r^2T)+ S_j(r^2T)\right) \geq N^r -1\right)\nonumber \\
& \leq  \sum_{j\in\mJ} \left[\pr\left( E^r_j(r^2T)\geq (1+\bar{\la}\vee\bar{\mu})r^2T\right) + \pr\left(S_j(r^2T)\geq (1+\bar{\la}\vee\bar{\mu})r^2T\right) \right]\nonumber\\
&= \sum_{j\in\mJ} \left[\pr\left( \bar{E}^r_j(T)\geq (1+\bar{\la}\vee\bar{\mu})T\right) + \pr\left(\bar{S}_j^r(T)\geq (1+\bar{\la}\vee\bar{\mu})T\right) \right]\rightarrow 0,\label{eq_pwc_9}
\end{align}
where \eqref{eq_pwc_9} is by FSLLN. The convergence result in \eqref{eq_pwc_9} implies that there are at most $O(r^2)$ review periods in the interval $[0,r^2T]$ with a high probability when $r$ is sufficiently large.

With the convention that $\infty - \infty := \infty$, let us define the following sets for all $\ep_2>0$ and $n,r\in\N_+$:
\begin{subequations}\label{eq_good_set}
\begin{align}
&\mA_{n}^{(1),r} (\ep_2):= \left\{ \tau_n^r > r^2T \right\}\cup \left\{ \tau_{n+1}^r - \tau_n^r \leq \ep_2 r \right\},\label{eq_good_set_1}\\
&\mA_{n}^{(2),r} (\ep_2):= \left\{ \tau_n^r > r^2T \right\}\cup\bigcap_{j\in\mJ^H} \left\{ \sup_{\tau_n^r\leq t\leq \tau_{n+1}^r}\left|Q_j^r(t) -Q_j^r(\tau_n^r)  \right| \leq C_6\ep_2 r \right\},\label{eq_good_set_2}\\
&\mA_{n}^{(3),r} (\ep_2):= \left\{ \tau_n^r > r^2T \right\}\cup \bigcap_{j\in\mJ^H} \left\{ \sup_{\tau_n^r\leq t\leq \tau_{n+1}^r}\left|q_j^{*,r}(t) -q_j^{*,r}(\tau_n^r)  \right| \leq C_7\ep_2 r \right\},\label{eq_good_set_3}\\
&\mA_{n}^{(4),r} (\ep_2):= \left\{ \tau_n^r > r^2T \right\}\cup \bigcap_{j\in\mJ^H} \left\{\left|Q_j^r(\tau_{n+1}^r) - q_j^{*,r}(\tau_{n+1}^r)  \right| \leq C_8\ep_2 r \right\},\label{eq_good_set_4}\\
&\mA_{n}^r (\ep_2):=\bigcap_{m=1}^4 \mA_{n}^{(m),r}(\ep_2),\label{eq_good_set_5}
\end{align}
\end{subequations}
where $C_6$, $C_7$, and $C_8$ are arbitrary strictly positive constants such that
\begin{equation}\label{eq_constants}
C_6 > \max\left\{2(1+\bar{\lambda}), 2\bar{\mu}\right\},\qquad 2C_7 < C_8<0.5\underline{\lambda}.
\end{equation}
We let $\mA_0^r(\ep_2):=\Omega$ for all $\ep_2>0$ and $r\in\N_+$ for completeness. 

The event in \eqref{eq_good_set_1} implies that the length of a review period is short in $[0,r^2T]$. The event in \eqref{eq_good_set_2} implies that the queue length processes in the buffers do not change a lot during a review period in $[0,r^2T]$. The event in \eqref{eq_good_set_3} implies that the optimal LP \eqref{eq_lp} solution does not change a lot during a review period in $[0,r^2T]$. The event in \eqref{eq_good_set_4} implies that the queue lengths in the buffers processed by resources in heavy traffic do not deviate a lot from the desired levels at the end of a review period which starts in $[0,r^2T]$. Finally, the following result states that the aforementioned events are realized jointly in the review periods $\{1,2,\ldots,N^r\}$ with high probability when $r$ is large.

\begin{lemma}\label{l_good_set}
For all $\ep_2>0$, we have
\begin{equation*}
\pr\left(\bigcap_{n=1}^{N^r} \mA_n^r(\ep_2)\right) \rightarrow 1.
\end{equation*}

\end{lemma}
The proof of Lemma \ref{l_good_set} is presented in E-companion \ref{l_good_set_proof}.

Let $\ep_2>0$ be such that $(C_6+C_7+C_8)\ep_2\leq \ep_1$. Then, the second sum in \eqref{eq_pwc_8} is less than or equal to
\begin{align}
&\sum_{j\in\mJ^H}\pr\left(\left\Vert\hat{Q}_{j}^r - \hat{q}_{j}^{*,r}\right\Vert_{T}>\ep_1,\; \tau_{N^r}^r > r^2T,\;  \bigcap_{n=1}^{N^r} \mA_n^r(\ep_2) \right)\label{eq_pwc_10}\\
&\hspace{6cm} + J\pr\left( \tau_{N^r}^r\leq r^2T\right) + J\pr\left( \left(\bigcap_{n=1}^{N^r} \mA_n^r(\ep_2)\right)^c\; \right),\label{eq_pwc_12}
\end{align}
where the superscript $c$ denotes the complement of the associated set. The sum in \eqref{eq_pwc_12} converges to 0 by \eqref{eq_pwc_9} and Lemma \ref{l_good_set}, respectively. The sum in \eqref{eq_pwc_10} is equal to
\begin{align}
&\sum_{j\in\mJ^H}\pr\left(\sup_{0\leq t\leq r^2T} \left|Q_{j}^r(t) - q_{j}^{*,r}(t)\right|>\ep_1r,\; \tau_{N^r}^r > r^2T,\;  \bigcap_{n=1}^{N^r} \mA_n^r(\ep_2) \right)\nonumber\\
&\hspace{1cm}\leq \sum_{j\in\mJ^H}\pr\Bigg( \bigcup_{n=1}^{N^r} \Bigg\{\sup_{\tau_n^r\leq t \leq \tau_{n+1}^r}\left| Q_j^r(t) - q_j^{*,r}(t)\right|>\ep_1r,\;\tau_n^r\leq r^2T \Bigg\},\; \bigcap_{n=1}^{N^r} \mA_n^r(\ep_2) \Bigg)\nonumber\\
&\hspace{1cm}\leq \sum_{j\in\mJ^H}\pr\Bigg( \bigcup_{n=1}^{N^r} \Bigg\{ \bigg( \sup_{\tau_n^r\leq t \leq \tau_{n+1}^r}\left| Q_j^r(t) - Q_j^r(\tau_n^r) \right| + \left| Q_j^r(\tau_n^r) - q_j^{*,r}(\tau_n^r) \right|\nonumber\\ 
&\hspace{4cm} + \sup_{\tau_n^r\leq t \leq \tau_{n+1}^r}\left| q_j^{*,r}(t) - q_j^{*,r}(\tau_n^r) \right|\bigg)>\ep_1r,\;\tau_n^r\leq r^2T \Bigg\},\; \bigcap_{n=1}^{N^r} \mA_n^r(\ep_2) \Bigg).\label{eq_pwc_13}
\end{align}
Let us focus on the event inside the probability in \eqref{eq_pwc_13}. In the set $\{\tau_n^r\leq r^2T\}\cap \mA_{n-1}^r(\ep_2) \cap \mA_n^r(\ep_2)$, 
\begin{align}
&\sup_{\tau_n^r\leq t \leq \tau_{n+1}^r}\left| Q_j^r(t) - Q_j^r(\tau_n^r) \right| + \left| Q_j^r(\tau_n^r) - q_j^{*,r}(\tau_n^r) \right| + \sup_{\tau_n^r\leq t \leq \tau_{n+1}^r}\left| q_j^{*,r}(t) - q_j^{*,r}(\tau_n^r) \right| \nonumber\\
&\hspace{10cm} \leq (C_6+C_7+C_8)\ep_2 r \leq \ep_1 r\label{eq_pwc_14}
\end{align}
for all $n\in\{1,2,\ldots,N^r\}$ and $r\in\N_+$ by \eqref{eq_good_set}. Hence, the event inside the probability in \eqref{eq_pwc_13} is equal to $\emptyset$ by \eqref{eq_pwc_14}. Therefore, the sum in \eqref{eq_pwc_10} is equal to 0.

\section{Proof of Lemma \ref{l_good_set}}\label{l_good_set_proof}

Let us define
\begin{equation}\label{eq_gsp_0}
B^r:=\left\{ \max_{j\in\mJ} \left(E_j^r(r^2T) \vee S_j(r^2T)\right)  \leq \rd{C_9r^2T} - 1,\;  \max_{j\in\mJ^H} \left| Q_j^{r}(0) - q_j^{*,r}(0)\right| \leq C_8\ep_2 r\right\},
\end{equation}
where $C_9 := 1+ \bar{\lambda} \vee \bar{\mu}$. By FSLLN (see for example \eqref{eq_pwc_9}) and \eqref{eq_ini2}, we have $\pr\left(B^r\right)\rightarrow 1$.

The following lemma from \cite{ozk20} provides an exponentially decaying upper bound on renewal processes and will be useful later.

\begin{lemma}\label{l_tail_1}
(Lemma 5 of \cite{ozk20}) Let us fix arbitrary $a>0$ and $b>0$. There exists an $r_2\in\N_+$ such that if $r\geq r_2$, then for all $j\in\mJ$ and $n\in\{1,2,\ldots,N^r\}$, we have
\begin{subequations}\label{eq_t1}
\begin{align*}
&\pr\left( \sup_{0\leq t\leq ar} \left| E_j^{r}(\tau_{n}^r + t) -  E_j^{r}(\tau_{n}^r) -\lambda_j^r t \right| >br,\;\tau_{n}^r \leq r^2T,\;B^r\right) \leq C_{10}r^{2}\e^{-C_{11}r},\\
&\pr\left( \sup_{0\leq t\leq ar} \left| S_j(T_j^r(\tau_{n}^r) + t) -  S_j(T_j^r(\tau_{n}^r)) -\mu_jt \right| >br,\;\tau_n^r \leq r^2T,\;B^r\right) \leq C_{10}r^{2}\e^{-C_{11}r},
\end{align*}
\end{subequations}
where $C_{10}$ and $C_{11}$ are strictly positive constants independent of $j$, $n$, and $r$.
\end{lemma}

The following lemma from \cite{ozk20} states that the workloads associated with resources in heavy traffic do not fluctuate a lot within a time interval with length $O(r)$ with high probability when $r$ is large.

\begin{lemma}\label{l_tail_2}
(Lemma 8 of \cite{ozk20}) Let us fix arbitrary $a>0$ and $b>0$. There exists an $r_3\in\N_+$ such that if $r\geq r_3$, then for all $n\in\{1,2,\ldots,N^r\}$, we have
\begin{equation*}
\pr\Bigg(\sup_{0\leq t\leq ar } \max_{i\in\mI^H} \left|W_i^{r}(\tau_n^r+t)- W_i^{r}(\tau_n^r) \right|>b r,\; \tau_n^r \leq r^2T,\;B^r\Bigg)\leq C_{12}r^5\e^{-C_{13}r},
\end{equation*}
where $C_{12}$ and $C_{13}$ are strictly positive constants independent of $n$ and $r$.
\end{lemma}

Because $\pr\left(B^r\right)\rightarrow 1$, proving Lemma \ref{l_good_set} is equivalent to proving 
\begin{equation*}
\pr\left( \bigcup_{n=1}^{N^r}\left(\mA_n^r(\ep_2)\right)^c,\;B^r\right)\rightarrow 0,\qquad\forall \ep_2>0.
\end{equation*}

Let us fix an arbitrary $\ep_2>0$. Let $\{A_n,n\in\N_+\}$ be an arbitrary sequence of sets. Then
\begin{equation}\label{eq_gsp_1}
\bigcup_{n=1}^N A_n=A_1 \cup \bigcup_{n=2}^N \left(A_n\cap A_{n-1}^c \cap A_{n-2}^c \cap \ldots \cap A_1^c\right),\qquad\forall N\in\N_+.
\end{equation}
Therefore,
\begin{equation*}
\pr\left( \bigcup_{n=1}^{N^r}\left(\mA_n^r\right)^c,B^r\right)\leq \sum_{n=1}^{N^r} \pr\left(\left(\mA_n^r\right)^c \cap \mA_{n-1}^r \cap B^r\right).
\end{equation*}

Let us fix an arbitrary $n\in\{1,2,\ldots,N^r\}$. By \eqref{eq_good_set_5} and \eqref{eq_gsp_1}, we have
\begin{subequations}\label{eq_gsp_2}
\begin{align}
\pr\left(\left(\mA_n^r\right)^c \cap \mA_{n-1}^r \cap B^r\right) &\leq \pr\left(\left(\mA_n^{(1),r}\right)^c \cap \mA_{n-1}^r \cap B^r\right)\label{eq_gsp_21}\\
&\hspace{1cm} +\pr\left(\left(\mA_n^{(2),r}\right)^c \cap \mA_n^{(1),r} \cap B^r\right)\label{eq_gsp_22}\\
&\hspace{1cm} +\pr\left(\left(\mA_n^{(3),r}\right)^c \cap \mA_n^{(1),r} \cap B^r\right)\label{eq_gsp_23}\\
&\hspace{1cm} +\pr\left(\left(\mA_n^{(4),r}\right)^c \cap \mA_n^{(1),r} \cap \mA_{n-1}^r \cap B^r\right).\label{eq_gsp_24}
\end{align}
\end{subequations}
We will derive exponentially decaying upper bounds on each of the probabilities in the RHS of \eqref{eq_gsp_2}. We will use the definition of the proposed policy to derive the desired bound on the probability in \eqref{eq_gsp_21}. We will use Lemma \ref{l_tail_1} to derive desired bound on the probability in \eqref{eq_gsp_22} and will use Lemmas \ref{l_reg} and \ref{l_tail_2} to derive desired bound on the probability in \eqref{eq_gsp_23}. Finally, we will use the structure of the proposed policy to derive the desired bound on the probability in \eqref{eq_gsp_24}.

\paragraph{The probability in the RHS of \eqref{eq_gsp_21}:} By \eqref{eq_good_set_1}, it is equal to
\begin{equation}\label{eq_a_1}
\pr\left( \tau_{n+1}^r - \tau_{n}^r > \ep_2 r,\; \tau_n^r\leq r^2T,\;\mA_{n-1}^r,\; B^r \right).
\end{equation}
We will consider \eqref{eq_a_1} case by case. First, suppose that the $n$th review period is Step 2 in Definition \ref{d_pp}. Then, the length of the $n$th review period is at most minimum of residual service times and residual inter-arrival times at time $\tau_n^r$. Therefore, by \eqref{eq_gsp_0}, if $r\geq r_1$ (see \eqref{eq_rdef_1}), \eqref{eq_a_1} is less than or equal to
\begin{align}
& \pr\left( \max_{m\in\left\{1,2,\ldots,\rd{C_9r^2T}\right\}}\max_{j\in\mJ} \left(u_{jm}^r\vee v_{jm}\right) > \ep_2 r \right)\nonumber\\
&\hspace{2cm}\leq \sum_{m=1}^{\rd{C_9r^2T}} \sum_{j\in\mJ} \left[\pr\left(u_{jm}^r> \ep_2 r \right)+\pr\left(v_{jm}> \ep_2 r \right)\right]\nonumber\\
&\hspace{2cm}\leq C_9r^2T \sum_{j\in\mJ} \left[\pr\left(u_{j1}^r> \ep_2 r \right)+\pr\left(v_{j1}> \ep_2 r \right)\right]\nonumber\\
&\hspace{2cm}= C_9r^2T \sum_{j\in\mJ} \left[\pr\left(\frac{\bar{u}_{j1}}{\la_j^r}> \ep_2 r \right)+\pr\left(v_{j1}> \ep_2 r \right)\right]\nonumber\\
&\hspace{2cm}\leq C_9r^2T \sum_{j\in\mJ} \left[\pr\left(\bar{u}_{j1}> 0.5\underline{\la}\ep_2 r \right)+\pr\left(v_{j1}> \ep_2 r \right)\right]\label{eq_a_2}\\
&\hspace{2cm}= C_9r^2T \sum_{j\in\mJ} \left[\pr\left(\e^{0.5\bar{\al}\bar{u}_{j1}}> \e^{0.25\bar{\al}\underline{\la}\ep_2 r} \right)+\pr\left(\e^{0.5\bar{\al}v_{j1}}> \e^{0.5\bar{\al}\ep_2 r} \right)\right]\label{eq_a_3}\\
&\hspace{2cm}\leq C_9r^2T \sum_{j\in\mJ} \left(\E\left[\e^{0.5\bar{\al}\bar{u}_{j1}}\right] \e^{-0.25\bar{\al}\underline{\la}\ep_2 r} +\E\left[\e^{0.5\bar{\al}v_{j1}}\right]\e^{-0.5\bar{\al}\ep_2 r}\right)\label{eq_a_4}\\
&\hspace{2cm}\leq C_{14}r^2\e^{-C_{15} r}, \label{eq_a_5}
\end{align}
where \eqref{eq_a_2} is by the fact that $r\geq r_1$ (see \eqref{eq_rdef_1}), $\bar{\al}$ used in \eqref{eq_a_3} is defined in Assumption \ref{a_moment}, \eqref{eq_a_4} is by Markov's inequality, and
\begin{equation*}
C_{14}:=C_9T J\max_{j\in\mJ} \left\{\E\left[\e^{0.5\bar{\al}\bar{u}_{j1}}\right]+ \E\left[\e^{0.5\bar{\al}v_{j1}}\right]\right\},\qquad C_{15}:= \left(0.5\wedge \left(0.25\underline{\lambda}\right)\right)\bar{\al}\ep_2.
\end{equation*}
Notice that both $C_{14}$ and $C_{15}$ are strictly positive and finite constants by Assumption \ref{a_moment} and independent of $n$ and $r$. 

Next, suppose that the $n$th review period is Step 3 in Definition \ref{d_pp}. Let $\kappa_n^r$ denote the length of the time interval between $\tau_n^r$ and the first time when an external arrival or service completion happens after $\tau_n^r$. Then, the length of the $n$th review period is equal to $L^r(\tau_n^r)\vee\kappa_n^r$ by definition of Step 3. Hence, \eqref{eq_a_1} is equal to
\begin{align}
& \pr\left( L^r(\tau_n^r)\vee\kappa_n^r > \ep_2 r,\; \tau_n^r\leq r^2T,\;\mA_{n-1}^r,\; B^r \right)\nonumber\\
&\hspace{3cm} \leq \pr\left( \kappa_n^r > \ep_2 r,\; \tau_n^r\leq r^2T,\; B^r \right) + \pr\left( L^r(\tau_n^r)> \ep_2 r,\; \tau_n^r\leq r^2T,\;\mA_{n-1}^r,\; B^r  \right).\label{eq_a_6}
\end{align}
Observe that if $r\geq r_1$ (see \eqref{eq_rdef_1}), the first probability in \eqref{eq_a_6} is less than or equal to the term in \eqref{eq_a_5}. By \eqref{eq_review_3} and \eqref{eq_review_4}, and if $r\geq r_1$, the second probability in \eqref{eq_a_6} is equal to
\begin{equation}\label{eq_a_7}
 \pr\left( \max_{i\in\mI^H} \max_{j\in\mJ_i^{\leq,r}(\tau_n^r)} \frac{q_j^{*,r}(\tau_n^r)-Q_j^r(\tau_n^r)}{\la_j^r} > \ep_2 r,\; \tau_n^r\leq r^2T,\;\mA_{n-1}^r,\; B^r  \right)\leq \I\left( \frac{C_8\ep_2r}{0.5\underline{\la}} > \ep_2 r\right)=0,
\end{equation}
where the inequality in \eqref{eq_a_7} is by \eqref{eq_good_set_4}, \eqref{eq_gsp_0}, and the fact that $r\geq r_1$ (see \eqref{eq_rdef_1}) and the equality in \eqref{eq_a_7} is by \eqref{eq_constants}.

Consequently, if $r\geq r_1$ (see \eqref{eq_rdef_1}), the probability in the RHS of \eqref{eq_gsp_21} is less than or equal to the term in \eqref{eq_a_5}.

\paragraph{The probability in \eqref{eq_gsp_22}:} By \eqref{eq_good_set_2}, it is less than or equal to
\begin{align}
&\sum_{j\in\mJ^H}\pr\left( \sup_{\tau_n^r\leq t\leq \tau_{n+1}^r}\left|Q_j^r(t) -Q_j^r(\tau_n^r)  \right| > C_6\ep_2 r,\;\tau_n^r \leq r^2T,\; \mA_n^{(1),r},\;B^r\right)\nonumber\\
&\hspace{0cm}\leq \sum_{j\in\mJ^H}\pr\left( \sup_{0\leq t\leq \ep_2 r}\left|Q_j^r(\tau_n^r +t) -Q_j^r(\tau_n^r)  \right| > C_6\ep_2 r,\;\tau_n^r \leq r^2T,\;B^r \right)\nonumber\\
&\hspace{0cm}= \sum_{j\in\mJ^H}\pr\bigg( \sup_{0\leq t\leq \ep_2 r}\left|E_{j}^r(\tau_n^r+t)-S_j(T_j^r(\tau_n^r+t)) - E_{j}^r(\tau_n^r)+S_j(T_j^r(\tau_n^r))\right| > C_6\ep_2 r,\; \tau_n^r \leq r^2T,\;B^r\bigg)\nonumber\\
&\hspace{0cm}\leq \sum_{j\in\mJ^H}\pr\left( S_j(T_j^r(\tau_n^r)+\ep_2 r)-S_j(T_j^r(\tau_n^r)) > 0.5C_6\ep_2 r,\;\tau_n^r \leq r^2T,\;B^r \right)\nonumber\\
&\hspace{4cm}+\sum_{j\in\mJ^H}\pr\left( E_{j}^r(\tau_n^r+\ep_2 r) - E_{j}^r(\tau_n^r) > 0.5C_6\ep_2 r,\;\tau_n^r \leq r^2T,\;B^r \right)\label{eq_b_1}\\
&\hspace{0cm}\leq\sum_{j\in\mJ^H}\pr\left( S_j(T_j^r(\tau_n^r)+\ep_2 r)-S_j(T_j^r(\tau_n^r)) -\mu_j\ep_2 r > \left(0.5C_6-\bar{\mu}\right)\ep_2 r,\;\tau_n^r \leq r^2T,\;B^r \right)\nonumber\\
&\hspace{1cm}+\sum_{j\in\mJ^H}\pr\left( E_{j}^r(\tau_n^r+\ep_2 r) - E_{j}^r(\tau_n^r) -\lambda_j^r \ep_2 r > \left(0.5C_6-1-\bar{\lambda}\right)\ep_2 r,\;\tau_n^r \leq r^2T,\;B^r \right),\label{eq_b_2}
\end{align}
where \eqref{eq_b_1} is by triangular inequality and the fact that $T_j^r(\tau_n^r+t) \leq T_j^r(\tau_n^r)+t$ for all $t\in\R_+$ and \eqref{eq_b_2} holds if $r\geq r_1$ (see \eqref{eq_rdef_1}). By \eqref{eq_constants} and Lemma \ref{l_tail_1}, there exists an $r_4\in\N_+$ independent of $n$ such that if $r\geq r_4$, the sum in \eqref{eq_b_2} is less than or equal to 
\begin{equation}\label{eq_b_3}
C_{16}r^2\e^{-C_{17}r},
\end{equation}
where $C_{16}$ and $C_{17}$ are strictly positive constants independent of $n$ and $r$. Finally, if $r\geq r_1\vee r_4$, the probability in \eqref{eq_gsp_22} is less than or equal to the term in \eqref{eq_b_3}.

\paragraph{The probability in \eqref{eq_gsp_23}:} By \eqref{eq_good_set_3}, it is less than or equal to
\begin{align}
&\sum_{j\in\mJ^H}\pr\left( \sup_{\tau_n^r\leq t\leq \tau_{n+1}^r}\left|q_j^{*,r}(t)- q_j^{*,r}(\tau_n^r) \right| > C_7\ep_2 r,\;\tau_n^r \leq r^2T,\; \mA_n^{(1),r},\;B^r\right)\label{eq_c_1}\\
&\hspace{0.5cm}\leq \sum_{j\in\mJ^H}\pr\left(C_2\sup_{\tau_n^r\leq t\leq \tau_{n+1}^r} \max_{i\in\mI^H} \left|W_i^{r}(t)- W_i^{r}(\tau_n^r) \right|>C_7\ep_2 r,\; \tau_n^r \leq r^2T,\; \mA_n^{(1),r},\;B^r\right)\label{eq_c_2}\\
&\hspace{0.5cm}\leq \sum_{j\in\mJ^H}\pr\left(\sup_{0\leq t\leq \ep_2 r} \max_{i\in\mI^H} \left|W_i^{r}(\tau_n^r+t)- W_i^{r}(\tau_n^r) \right| >\frac{C_7}{C_2}\ep_2 r,\; \tau_n^r \leq r^2T,\;B^r\right),\label{eq_c_3}
\end{align}
where \eqref{eq_c_2} is by Lemma \ref{l_reg}. By Lemma \ref{l_tail_2}, there exists an $r_5\in\N_+$ independent of $n$ such that if $r\geq r_5$, the sum in \eqref{eq_c_3} is less than or equal to
\begin{equation}\label{eq_c_4}
C_{18}r^5 \e^{-C_{19}r},
\end{equation}
where $C_{18}$ and $C_{19}$ are strictly positive constants independent of $n$ and $r$.

\paragraph{The probability in \eqref{eq_gsp_24}:} Because deriving an exponentially decaying upper bound on the probability in \eqref{eq_gsp_24} requires exploitation of the structure of the proposed policy, it is long and involved. Hence, we present the following result whose proof is presented in E-companion \ref{l_good_4_proof}: There exists an $r_6\in\N_+$ independent of $n$ such that if $r\geq r_6$,
\begin{equation}\label{eq_d_1}
\pr\left(\left(\mA_n^{(4),r}\right)^c \cap \mA_n^{(1),r} \cap \mA_{n-1}^r \cap B^r\right) \leq C_{20}r^5 \e^{-C_{21}r},
\end{equation}
where $C_{20}$ and $C_{21}$ are strictly positive constants independent of $n$ and $r$.

Finally, let $r_7:=r_1\vee r_4\vee r_5\vee r_6$. Then, $r_7$ is independent of $n$. By \eqref{eq_gsp_2}, \eqref{eq_a_5}, \eqref{eq_b_3}, \eqref{eq_c_4}, and \eqref{eq_d_1}, if $r\geq r_7$ 
\begin{equation*}
\pr\left(\left(\mA_n^r\right)^c \cap \mA_{n-1}^r \cap B^r\right)\leq C_{22}r^5 \e^{-C_{23}r},
\end{equation*}
where $C_{22}$ and $C_{23}$ are strictly positive constants independent of $n$ and $r$. Consequently, if $r\geq r_7$, 
\begin{align*}
\pr\left( \bigcup_{n=1}^{N^r}\left(\mA_n^r\right)^c,B^r\right)&\leq \sum_{n=1}^{N^r} \pr\left(\left(\mA_n^r\right)^c \cap \mA_{n-1}^r \cap B^r\right)\\
&\leq N^rC_{22} r^{5}\e^{-C_{23}r} \leq \left(2+2\left(1+\bar{\la}\vee\bar{\mu}\right)JT\right) C_{22} r^{7}\e^{-C_{23}r},
\end{align*}
which converges to $0$ exponentially fast and this completes the proof.

\section{Derivation of \eqref{eq_d_1}}\label{l_good_4_proof}

By \eqref{eq_good_set_4}, the probability in the LHS of \eqref{eq_d_1} is less than or equal to
\begin{align}
&\sum_{j\in\mJ^H}\pr\left( \left|Q_j^r(\tau_{n+1}^r)- q_j^{*,r}(\tau_{n+1}^r) \right| > C_8\ep_2 r,\;\tau_n^r \leq r^2T,\; \mA_n^{(1),r},\;\mA_{n-1}^r,\;B^r\right)\nonumber\\
&\hspace{1cm}\leq\sum_{j\in\mJ^H}\pr\left( \left|q_j^{*,r}(\tau_{n+1}^r)- q_j^{*,r}(\tau_n^r) \right| > 0.5C_8\ep_2 r,\;\tau_n^r \leq r^2T,\; \mA_n^{(1),r},\;B^r\right)\label{eq_d_4}\\
&\hspace{2cm} +\sum_{j\in\mJ^H}\pr\left( \left|Q_j^r(\tau_{n+1}^r)- q_j^{*,r}(\tau_n^r)\right| > 0.5C_8\ep_2 r,\;\tau_n^r \leq r^2T,\; \mA_n^{(1),r},\;\mA_{n-1}^r,\;B^r\right).\label{eq_d_5}
\end{align}
By \eqref{eq_constants}, the sum in \eqref{eq_d_4} is less than or equal to the sum in \eqref{eq_c_1}. Therefore, if $r\geq r_5$, the sum in \eqref{eq_d_4} is less than or equal to the term in \eqref{eq_c_4}. Next, let us consider the sum in \eqref{eq_d_5} case by case. There are two cases.

\paragraph{Case 1 associated with \eqref{eq_d_5}} Suppose that the $n$th review period is Step 2 in Definition \ref{d_pp}. Let us fix an arbitrary $r\in\N_+$. Recall that Step 2 ends at the first time when an external arrival or a service completion happens after $\tau_n^r$. Therefore,
\begin{equation}\label{eq_d_6}
\left|Q_j^r(\tau_{n+1}^r) -Q_j^r(\tau_n^r)\right| \leq 1,\qquad \forall j\in\mJ^H.
\end{equation}
By \eqref{eq_d_6} and because $\mJ^{>,r}(\tau_n^r)=\emptyset$ and so $Q_j^r(\tau_n^r)\leq \lceil q_j^{*,r}(\tau_n^r) \rceil$ for all $j\in\mJ^H$ by definition of Step 2, we have $Q_j^r(\tau_{n+1}^r) -\lceil q_j^{*,r}(\tau_n^r) \rceil \leq 1$ for all $j\in\mJ^H$, which implies
\begin{equation}\label{eq_d_7}
Q_j^r(\tau_{n+1}^r) - q_j^{*,r}(\tau_n^r) \leq 2,\qquad\forall j\in\mJ^H.
\end{equation}
Let us consider an arbitrary $j\in\mJ^H$. Let $i\in\mI_j\cap\mI^H$ be an arbitrary resource in heavy traffic processing type $j$ jobs. By \eqref{eq_we} and the fact that $Q_l^r(\tau_n^r)\leq \lceil q_l^{*,r}(\tau_n^r) \rceil$ for all $l\in\mJ_i$, we have
\begin{equation*}
q_j^{*,r}(\tau_n^r) - Q_j^r(\tau_n^r) = \mu_j\sum_{l\in\mJ_i\backslash\{j\}} \frac{Q_l^r(\tau_n^r)-q_l^{*,r}(\tau_n^r) }{\mu_l} \leq (J-1)\frac{\bar{\mu}}{\underline{\mu}}.
\end{equation*}
Then, by \eqref{eq_d_6}, we have
\begin{equation}\label{eq_d_8}
q_j^{*,r}(\tau_n^r) - Q_j^r(\tau_{n+1}^r)  \leq 1+ (J-1)\frac{\bar{\mu}}{\underline{\mu}},\qquad \forall j\in\mJ^H.
\end{equation}
Therefore, by \eqref{eq_d_7} and \eqref{eq_d_8}, we have
\begin{equation}\label{eq_d_9}
\left| Q_j^r(\tau_{n+1}^r)-q_j^{*,r}(\tau_n^r) \right| \leq 2+ (J-1)\frac{\bar{\mu}}{\underline{\mu}},\qquad \forall j\in\mJ^H.
\end{equation}
By \eqref{eq_d_9}, there exists an $r_8\in\N_+$ independent of $n$ such that if the $n$th review period is Step 2 in Definition \ref{d_pp}, then the sum in \eqref{eq_d_5} is less than or equal to
\begin{equation}\label{eq_d_10}
\I\left( 2+ (J-1)\frac{\bar{\mu}}{\underline{\mu}} > 0.5C_8\ep_2 r \right) =0,\qquad\forall  r\geq r_8.
\end{equation}

\paragraph{Case 2 associated with \eqref{eq_d_5}} Suppose that the $n$th review period is Step 3 in Definition \ref{d_pp}. Then, the sum in \eqref{eq_d_5} is less than or equal to
\begin{align}
&\sum_{j\in\mJ^H}\pr\left( \left|Q_j^r(\tau_{n+1}^r)- q_j^{*,r}(\tau_n^r)\right| > 0.5C_8\ep_2 r,\;\kappa_n^r> L^r(\tau_n^r),\;\tau_n^r \leq r^2T,\;\mA_{n-1}^r,\;B^r\right)\label{eq_d_11}\\
&\hspace{0cm} +\sum_{j\in\mJ^H}\pr\left( \left|Q_j^r(\tau_{n+1}^r)- q_j^{*,r}(\tau_n^r)\right| > 0.5C_8\ep_2 r,\;\kappa_n^r \leq L^r(\tau_n^r),\;\tau_n^r \leq r^2T,\; \mA_n^{(1),r},\;\mA_{n-1}^r,\;B^r\right).\label{eq_d_12}
\end{align}
On the one hand, the sum in \eqref{eq_d_11} considers the case in which Step 3 ends at the first time when an external arrival or a service completion happens after $\tau_n^r$, that is, $\kappa_n^r > L^r(\tau_n^r)$. On the other hand, the sum in \eqref{eq_d_12} considers the case in which Step 3 ends at time $L^r(\tau_n^r)$. 

\paragraph{The sum in \eqref{eq_d_11}} First, let us consider the sum in \eqref{eq_d_11}. In the set $\{\tau_n^r \leq r^2T\}\cap\mA_{n-1}^r \cap B^r$ (see \eqref{eq_good_set_4} and \eqref{eq_gsp_0}), we have
\begin{equation*}
\left|Q_j^r(\tau_n^r) - q_j^{*,r}(\tau_n^r)  \right| \leq C_8\ep_2 r,\qquad\forall j\in\mJ^H.
\end{equation*}
Because Step 3 ends at the first time when an external arrival or a service completion happens after $\tau_n^r$, then \eqref{eq_d_6} holds. Let us consider an arbitrary $j\in\mJ^H$ and let $i\in\mI_j\cap\mI^H$ be an arbitrary resource in heavy traffic processing type $j$ jobs. There are two sub-cases. 

\paragraph{Case 1 associated with \eqref{eq_d_11}} Suppose that 
\begin{equation*}
\left|Q_j^r(\tau_n^r) - q_j^{*,r}(\tau_n^r)  \right| < 0.5C_8\ep_2 r -1.
\end{equation*}
Then, by \eqref{eq_d_6}, we have
\begin{equation*}
\left|Q_j^r(\tau_{n+1}^r) - q_j^{*,r}(\tau_n^r)  \right| < 0.5C_8\ep_2 r.
\end{equation*}
Hence,
\begin{align}
&\pr\Big( \left|Q_j^r(\tau_{n+1}^r)- q_j^{*,r}(\tau_n^r)\right| > 0.5C_8\ep_2 r,\;\left|Q_j^r(\tau_n^r) - q_j^{*,r}(\tau_n^r)  \right| < 0.5C_8\ep_2 r,\nonumber\\
&\hspace{8cm}\kappa_n^r> L^r(\tau_n^r),\;\tau_n^r \leq r^2T,\;\mA_{n-1}^r,\;B^r\Big) =0.\label{eq_d_13}
\end{align}

\paragraph{Case 2 associated with \eqref{eq_d_11}} Suppose that 
\begin{equation*}
\left|Q_j^r(\tau_n^r) - q_j^{*,r}(\tau_n^r) \right| \geq 0.5C_8\ep_2 r -1.
\end{equation*}

\paragraph{Case 2.1 associated with \eqref{eq_d_11}} Suppose that $Q_j^r(\tau_n^r) - q_j^{*,r}(\tau_n^r) \geq 0.5C_8\ep_2 r-1$. Then, by \eqref{eq_we}, 
\begin{equation*}
\sum_{l\in\mJ_i\backslash\{j\}}\frac{q_l^{*,r}(\tau_n^r)-Q_l^r(\tau_n^r)}{\mu_l} = \frac{Q_j^r(\tau_n^r) - q_j^{*,r}(\tau_n^r)}{\mu_j} \geq \frac{0.5C_8\ep_2 r-1}{\mu_j},
\end{equation*}
which implies that there exists $l\in\mJ_i\backslash\{j\}$ such that
\begin{equation*}
\frac{q_l^{*,r}(\tau_n^r)-Q_l^r(\tau_n^r)}{\mu_l} \geq\frac{0.5C_8\ep_2 r-1}{(J-1)\mu_j}. 
\end{equation*}
Therefore, by \eqref{eq_review_3} and \eqref{eq_review_4}, if $r\geq r_1$ (see \eqref{eq_rdef_1}), it must be the case that
\begin{equation}\label{eq_d_14}
L^r(\tau_n^r)\geq L_i^r(\tau_n^r)\geq\frac{\underline{\mu}(0.5C_8\ep_2 r-1)}{(J-1)\bar{\mu}(1+\bar{\la})}.
\end{equation}

\paragraph{Case 2.2 associated with \eqref{eq_d_11}} Suppose that $q_j^{*,r}(\tau_n^r) - Q_j^r(\tau_n^r)\geq 0.5C_8\ep_2 r-1$. Then, by \eqref{eq_review_3} and \eqref{eq_review_4}, if $r\geq r_1$, again \eqref{eq_d_14} holds. 

Therefore, in both Case 2.1 and Case 2.2, \eqref{eq_d_14} holds. Consequently, by \eqref{eq_d_14} and because $j\in\mJ^H$ is arbitrarily chosen,
\begin{align}
&\sum_{j\in\mJ^H}\pr\Big( \left|Q_j^r(\tau_{n+1}^r)- q_j^{*,r}(\tau_n^r)\right| > 0.5C_8\ep_2 r,\;\left|Q_j^r(\tau_n^r) - q_j^{*,r}(\tau_n^r)  \right| \geq 0.5C_8\ep_2 r-1,\nonumber\\
&\hspace{9cm}\kappa_n^r> L^r(\tau_n^r),\;\tau_n^r \leq r^2T,\;\mA_{n-1}^r,\;B^r\Big)\nonumber\\
&\hspace{1cm}\leq \sum_{j\in\mJ^H}\pr\left( \kappa_n^r > \frac{\underline{\mu}(0.5C_8\ep_2 r-1)}{(J-1)\bar{\mu}(1+\bar{\la})},\;\tau_n^r \leq r^2T,\;\mA_{n-1}^r,\;B^r\right)\nonumber\\
&\hspace{1cm}\leq J\pr\left( \kappa_n^r> \frac{\underline{\mu}(0.5C_8\ep_2 r-1)}{(J-1)\bar{\mu}(1+\bar{\la})},\;\tau_n^r \leq r^2T,\;\mA_{n-1}^r,\;B^r\right)\label{eq_d_15}.
\end{align}
Similar to how we derive \eqref{eq_a_5}, we can prove that there exists $r_9\in\N_+$ independent of $n$ such that if $r>r_9$, the term in \eqref{eq_d_15} is less than or equal to 
\begin{equation}\label{eq_d_16}
C_{24}r^2\e^{-C_{25}r},
\end{equation}
where $C_{24}$ and $C_{25}$ are strictly positive constants independent of $n$ and $r$. Finally, by also considering \eqref{eq_d_13}, if $r\geq r_9$, the sum in \eqref{eq_d_11} is less than or equal to the term in \eqref{eq_d_16}. 

\subsection{The sum in \eqref{eq_d_12}} 

Let us define the set
\begin{align}
\mE_n^r:=&\bigg\{ \sup_{0\leq t\leq \ep_2r} \left| E_j^{r}(\tau_{n}^r + t) -  E_j^{r}(\tau_{n}^r) -\lambda_j^r t \right|\leq\ep_3 r, \nonumber\\
&\hspace{3cm} \sup_{0\leq t\leq \ep_2r} \left| S_j(T_j^r(\tau_{n}^r) + t) -  S_j(T_j^r(\tau_{n}^r)) -\mu_jt \right| \leq \ep_3 r,\quad\forall j\in\mJ \bigg\},\label{eq_d_17}
\end{align}
where
\begin{equation}\label{eq_d_17_1}
\ep_3:= \frac{C_8}{4\left(1+ 2J^2\frac{\bar{\mu}}{\underline{\mu}}\right)} \ep_2.
\end{equation}
Then, the sum in \eqref{eq_d_12} is less than or equal to
\begin{align}
&\sum_{j\in\mJ^H}\pr\left( \left(\mE_n^r\right)^c,\; \tau_n^r \leq r^2T,\;B^r \right)\label{eq_d_19}\\
&\hspace{0.5cm}+\sum_{j\in\mJ^H}\pr\left( \left|Q_j^r(\tau_{n+1}^r)- q_j^{*,r}(\tau_n^r)\right| > 0.5C_8\ep_2 r,\; \mE_n^r,\; \kappa_n^r \leq L^r(\tau_n^r),\;\tau_n^r \leq r^2T,\; \mA_n^{(1),r},\;\mA_{n-1}^r\right).\label{eq_d_18}
\end{align}
By Lemma \ref{l_tail_1}, there exists an $r_{10}\in\N_+$ such that $r_{10}$ is independent of $n$ and if $r\geq r_{10}$, the sum in \eqref{eq_d_19} is less than or equal to
\begin{equation}\label{eq_d_20}
C_{26}r^{2}\e^{-C_{27}r},
\end{equation}
where $C_{26}$ and $C_{27}$ are strictly positive constants independent of $n$ and $r$.

Let us define the following shifted processes for all $t\in\R_+$, $i\in\mI^H$, $j\in\mJ^H$, $n\in\{1,2,\ldots,N^r\}$, and $r\in\N_+$,
\begin{equation*}
I_i^{n,r}(t):= I_i^r(\tau_n^r + t) - I_i^r(\tau_n^r),\hspace{2cm} T_j^{n,r}(t):= T_j^r(\tau_n^r + t) - T_j^r(\tau_n^r).
\end{equation*}
Observe that
\begin{equation}\label{eq_d_21}
I_i^{n,r}(t) + \sum_{j\in\mJ_i} T_j^{n,r}(t) =t,\qquad\forall t\in\R_+,i\in\mI^H.
\end{equation}

The following result, which states that the shifted cumulative busy time processes associated with the job types in the set $\mJ^H$ do not deviate a lot from their expected values, will be very helpful in order to prove that the sum in \eqref{eq_d_18} is equal to 0.

\begin{lemma}\label{l_busy}
In Step 3 and in the set $\left\{\mE_n^r,\; \kappa_n^r \leq L^r(\tau_n^r),\;\tau_n^r\leq r^2T,\;\mA_n^{(1),r}(\ep_2),\;\mA_{n-1}^r(\ep_2)\right\}$, there exists an $r_{11}\in\N_+$ independent of $n$ such that if $r\geq r_{11}$, then
\begin{equation*}
\max_{j\in\mJ^H} \left| T_j^{n,r}(L^r(\tau_n^r)) - \frac{1}{\mu_j}\left(Q_j^r(\tau_n^r)-q_j^{*,r}(\tau_n^r) + \la_j^rL^r(\tau_n^r)\right)\right| \leq \frac{4J^2}{\underline{\mu}}\ep_3 r.
\end{equation*}
\end{lemma}
The proof of Lemma \ref{l_busy} is presented in E-companion \ref{l_busy_proof}.

Next, let us consider the sum in \eqref{eq_d_18}. Let us fix arbitrary $r\geq r_{11}$, $j\in\mJ^H$, and sample path in the set $\left\{\mE_n^r,\; \kappa_n^r \leq L^r(\tau_n^r),\;\tau_n^r\leq r^2T,\;\mA_n^{(1),r},\;\mA_{n-1}^r\right\}$. Then,
\begin{align}
Q_j^r(\tau_{n+1}^r) &= Q_j^r(\tau_n^r) + E_j^r(\tau_{n+1}^r) - E_j^r(\tau_n^r) - S_j(T_j^r(\tau_{n+1}^r))+S_j(T_j^r(\tau_n^r))\nonumber\\
& = Q_j^r(\tau_n^r) + E_j^r(\tau_n^r+L^r(\tau_n^r)) - E_j^r(\tau_n^r) - \left(S_j\left(T_j^r(\tau_n^r) +T_j^{n,r}(L^r(\tau_n^r))\right)-S_j(T_j^r(\tau_n^r))\right)\nonumber\\
& \leq Q_j^r(\tau_n^r) + \la_j^r L^r(\tau_n^r) +\ep_3r - \mu_jT_j^{n,r}(L^r(\tau_n^r))+\ep_3r \label{eq_d_22}\\
& \leq Q_j^r(\tau_n^r) + \la_j^r L^r(\tau_n^r) +2\ep_3r - \left(Q_j^r(\tau_n^r)-q_j^{*,r}(\tau_n^r) + \la_j^rL^r(\tau_n^r)\right) + \frac{4J^2\bar{\mu}}{\underline{\mu}}\ep_3 r\label{eq_d_23}\\
& =  q_j^{*,r}(\tau_n^r)+\left(2+ \frac{4\bar{\mu}J^2}{\underline{\mu}}\right)\ep_3 r,\label{eq_d_24}
\end{align}
where \eqref{eq_d_22} is by \eqref{eq_d_17} and the fact that $L^r(\tau_n^r) = \tau_{n+1}^r - \tau_n^r \leq \ep_2r$ and \eqref{eq_d_23} is by Lemma \ref{l_busy}. Similarly,
\begin{align}
Q_j^r(\tau_{n+1}^r) &= Q_j^r(\tau_n^r) + E_j^r(\tau_{n+1}^r) - E_j^r(\tau_n^r) - S_j(T_j^r(\tau_{n+1}^r))+S_j(T_j^r(\tau_n^r))\nonumber\\
& \geq Q_j^r(\tau_n^r) + \la_j^r L^r(\tau_n^r) -\ep_3r - \mu_jT_j^{n,r}(L^r(\tau_n^r))-\ep_3r \nonumber\\
& \geq Q_j^r(\tau_n^r) + \la_j^r L^r(\tau_n^r) -2\ep_3r - \left(Q_j^r(\tau_n^r)-q_j^{*,r}(\tau_n^r) + \la_j^rL^r(\tau_n^r)\right) - \frac{4J^2\bar{\mu}}{\underline{\mu}}\ep_3 r\nonumber\\
& =  q_j^{*,r}(\tau_n^r)-\left(2+ \frac{4J^2\bar{\mu}}{\underline{\mu}}\right)\ep_3 r.\label{eq_d_25}
\end{align}
By \eqref{eq_d_17_1}, \eqref{eq_d_24}, and \eqref{eq_d_25}, we have
\begin{equation}\label{eq_d_26}
\left|Q_j^r(\tau_{n+1}^r) - q_j^{*,r}(\tau_n^r)\right| \leq \left(2+ \frac{4J^2\bar{\mu}}{\underline{\mu}}\right)\ep_3 r = 0.5C_8\ep_2 r.
\end{equation}
By \eqref{eq_d_26}, if $r\geq r_{11}$, then the sum in \eqref{eq_d_18} is equal to 0. Therefore, if $r\geq r_{10}\vee r_{11}$, the sum in \eqref{eq_d_12} is less than or equal to the term in \eqref{eq_d_20}.

Therefore, if $r\geq r_1\vee r_8\vee r_9 \vee r_{10}\vee r_{11}$, then the sum in \eqref{eq_d_5} is less than or equal to
\begin{equation}\label{eq_d_27}
C_{28}r^{2}\e^{-C_{29}r},
\end{equation}
where $C_{28}:=C_{24}+C_{26}$ and $C_{29}:=C_{25}\wedge C_{27}$. 

Finally, we complete the derivation of \eqref{eq_d_1} by letting $r_6:= r_1\vee r_5\vee r_8\vee r_9 \vee r_{10}\vee r_{11}$, $C_{20}:= C_{18}+C_{28}$, and $C_{21}:= C_{19}\wedge C_{29}$.

\subsubsection{Proof of Lemma \ref{l_busy}}\label{l_busy_proof}

First, we will derive an upper bound on the shifted cumulative idle time at the end of the review period. Observe that in the set $\left\{\mE_n^r,\; \kappa_n^r \leq L^r(\tau_n^r),\;\tau_n^r\leq r^2T,\;\mA_n^{(1),r},\;\mA_{n-1}^r\right\}$, we have $L^r(\tau_n^r) \leq \ep_2r$ by \eqref{eq_good_set_1}. Let us fix an arbitrary $i\in\mI^H$. Let
\begin{equation*}
\eta^{n,r}:=\begin{cases}
\sup\left\{t\in[0,L^r(\tau_n^r)]: \dr I_i^{n,r}(t)>0 \right\},&\mbox{if there exists $t\in[0,L^r(\tau_n^r)]$ such that $\dr I_i^{n,r}(t)>0$,}\\
0,&\mbox{otherwise}.
\end{cases}
\end{equation*}
Because $\tau_n^r + \eta^{n,r}$ is the last time resource $i$ idles during Step 3, $I_i^{n,r}(\eta^{n,r}) = I_i^{n,r}(L^r(\tau_n^r))$. If $\eta^{n,r}=0$, then $I_i^{n,r}(L^r(\tau_n^r))=0$ and so we have a trivial upper bound on $I_i^{n,r}(L^r(\tau_n^r))$. 

Suppose that $\eta^{n,r}>0$. Let $r_{12}\geq r_1$ (see \eqref{eq_rdef_2}) be such that if $r\geq r_{12}$, then
\begin{equation}\label{eq_rdef_3}
\ep_2 \max_{i\in\mI^H}(1+|\theta_i|) + \frac{J}{\underline{\mu}}\leq \frac{J\ep_3r}{\underline{\mu}}.
\end{equation}
Then, $r_{12}$ is independent of $n$. Because resource $i$ is utilized in work conserving fashion, by the definition of $\eta^{n,r}$, because $\eta^{n,r}>0$, and by \eqref{eq_workload}, we have
\begin{equation*}
0=W_i^r(\tau_n^r + \eta^{n,r}) = W_i^r(\tau_n^r)+ \sum_{j\in\mJ_i}\frac{E_j^r(\tau_n^r + \eta^{n,r})- E_j^r(\tau_n^r)- \left(S_j(T_j^r(\tau_n^r + \eta^{n,r}))-S_j(T_j^r(\tau_n^r))\right)}{\mu_j},
\end{equation*}
which implies that if $r\geq r_{12}$, then
\begin{align}
& W_i^r(\tau_n^r)+ \sum_{j\in\mJ_i}\frac{E_j^r(\tau_n^r + \eta^{n,r})- E_j^r(\tau_n^r)}{\mu_j} = \sum_{j\in\mJ_i}\frac{S_j(T_j^r(\tau_n^r + \eta^{n,r}))-S_j(T_j^r(\tau_n^r))}{\mu_j},\nonumber\\
\implies&  W_i^r(\tau_n^r)+ \sum_{j\in\mJ_i}\frac{E_j^r(\tau_n^r + \eta^{n,r})- E_j^r(\tau_n^r)}{\mu_j}=\sum_{j\in\mJ_i}\frac{S_j\left(T_j^r(\tau_n^r )+ T_j^{n,r}(\eta^{n,r})\right)-S_j(T_j^r(\tau_n^r))}{\mu_j},\nonumber\\
\implies& W_i^r(\tau_n^r)+ \sum_{j\in\mJ_i}\frac{\la_j^r \eta^{n,r} -\ep_3r}{\mu_j}\leq \sum_{j\in\mJ_i} \left( T_j^{n,r}(\eta^{n,r}) + \frac{\ep_3r}{\mu_j}\right),\label{eq_idle_p_1}\\
\implies&W_i^r(\tau_n^r)+  \eta^{n,r}\rho_i^r  -\sum_{j\in\mJ_i}\frac{\ep_3r}{\mu_j} \leq \eta^{n,r}-I_i^{n,r}(\eta^{n,r}) + \sum_{j\in\mJ_i}\frac{\ep_3r}{\mu_j},\label{eq_idle_p_2}\\
\implies&I_i^{n,r}(\eta^{n,r})\leq  \eta^{n,r} \left|1- \rho_i^r\right| + \frac{2J\ep_3r}{\underline{\mu}} ,\nonumber\\
\implies& I_i^{n,r}(L^r(\tau_n^r))\leq \ep_2 r \left|1-\rho_i^r\right| + \frac{2J\ep_3r}{\underline{\mu}},\label{eq_idle_p_3}\\
\implies&I_i^{n,r}(L^r(\tau_n^r))\leq  \ep_2 \max_{i\in\mI^H}(1+|\theta_i|) + \frac{2J\ep_3r}{\underline{\mu}} ,\label{eq_idle_p_4}\\
\implies& I_i^{n,r}(L^r(\tau_n^r))\leq \frac{3J\ep_3r}{\underline{\mu}} ,\label{eq_idle_p_5}
\end{align}
where \eqref{eq_idle_p_1} is by \eqref{eq_d_17} and the fact that $T_j^{n,r}(\eta^{n,r})\leq \eta^{n,r}\leq L^r(\tau_n^r)\leq \ep_2r$, \eqref{eq_idle_p_2} is by \eqref{eq_d_21}, \eqref{eq_idle_p_3} is by the fact that $I_i^{n,r}(\eta^{n,r}) = I_i^{n,r}(L^r(\tau_n^r))$, \eqref{eq_idle_p_4} holds by Assumption \ref{a_regime} Part 2 and the fact that $r\geq r_{12}\geq r_1$ (see \eqref{eq_rdef_2} and \eqref{eq_rdef_3}), and \eqref{eq_idle_p_5} holds because $r\geq r_{12}$ (see \eqref{eq_rdef_3}). Because $i\in\mI^H$ is chosen arbitrarily, in the set $\left\{\mE_n^r,\; \kappa_n^r \leq L^r(\tau_n^r),\;\tau_n^r\leq r^2T,\;\mA_n^{(1),r},\;\mA_{n-1}^r\right\}$, we have
\begin{equation}\label{eq_idle2}
I_i^{n,r}(L^r(\tau_n^r))\leq \frac{3J\ep_3r}{\underline{\mu}},\qquad\forall i\in\mI^H.
\end{equation}

Next, we will consider the shifted cumulative busy time processes. Let us fix an arbitrary sample path in the set $\left\{\mE_n^r,\; \kappa_n^r \leq L^r(\tau_n^r),\;\tau_n^r\leq r^2T,\;\mA_n^{(1),r},\;\mA_{n-1}^r\right\}$. On the one hand, if $Q_j^r(\tau_n^r+ t) \geq q_j^{*,r}(\tau_n^r)$ for some $t\in[0,L^r(\tau_n^r)]$ and $j\in\mJ^H$, then
\begin{align}
& q_j^{*,r}(\tau_n^r) \leq Q_j^r(\tau_n^r) + E_j^r(\tau_n^r + t)- E_j^r(\tau_n^r) - S_j(T_j^r(\tau_n^r + t)) +S_j(T_j^r(\tau_n^r))\nonumber\\
\implies &S_j(T_j^r(\tau_n^r + t))-S_j(T_j^r(\tau_n^r)) \leq Q_j^r(\tau_n^r) -q_j^{*,r}(\tau_n^r) + E_j^r(\tau_n^r + t)- E_j^r(\tau_n^r)\nonumber\\
\implies & S_j\left(T_j^r(\tau_n^r)+T_j^{n,r}(t)\right)-S_j(T_j^r(\tau_n^r))  \leq Q_j^r(\tau_n^r) -q_j^{*,r}(\tau_n^r) + E_j^r(\tau_n^r + t)- E_j^r(\tau_n^r)\nonumber\\
\implies & \mu_jT_j^{n,r}(t) -\ep_3r  \leq Q_j^r(\tau_n^r) -q_j^{*,r}(\tau_n^r) + \la_j^r t + \ep_3 r \label{eq_busy_p_1}\\
\implies & T_j^{n,r}(t) \leq \frac{1}{\mu_j}\left(Q_j^r(\tau_n^r) -q_j^{*,r}(\tau_n^r) + \la_j^r t \right)+ \frac{2\ep_3 r}{\underline{\mu}}, \label{eq_busy_p_2}
\end{align}
where \eqref{eq_busy_p_1} is by \eqref{eq_d_17} and the fact that $T_j^{n,r}(t)\leq t \leq L^r(\tau_n^r)\leq \ep_2r$.

On the other hand, if $Q_j^r(\tau_n^r+ t) \leq \lceil q_j^{*,r}(\tau_n^r) \rceil$ for some $t\in[0,L^r(\tau_n^r)]$ and $j\in\mJ^H$, then
\begin{align}
& \lceil q_j^{*,r}(\tau_n^r) \rceil \geq Q_j^r(\tau_n^r) + E_j^r(\tau_n^r + t)- E_j^r(\tau_n^r) - S_j(T_j^r(\tau_n^r + t)) +S_j(T_j^r(\tau_n^r))\nonumber\\
\implies &S_j(T_j^r(\tau_n^r + t))-S_j(T_j^r(\tau_n^r)) \geq Q_j^r(\tau_n^r) -\lceil q_j^{*,r}(\tau_n^r) \rceil + E_j^r(\tau_n^r + t)- E_j^r(\tau_n^r)\nonumber\\
\implies & S_j\left(T_j^r(\tau_n^r)+T_j^{n,r}(t)\right)-S_j(T_j^r(\tau_n^r))  \geq Q_j^r(\tau_n^r) -\lceil q_j^{*,r}(\tau_n^r) \rceil + E_j^r(\tau_n^r + t)- E_j^r(\tau_n^r)\nonumber\\
\implies & \mu_jT_j^{n,r}(t) +\ep_3r  \geq Q_j^r(\tau_n^r) -\lceil q_j^{*,r}(\tau_n^r) \rceil + \la_j^r t - \ep_3 r \nonumber\\
\implies & T_j^{n,r}(t) \geq \frac{1}{\mu_j}\left(Q_j^r(\tau_n^r) -\lceil q_j^{*,r}(\tau_n^r) \rceil + \la_j^r t \right)- \frac{2\ep_3 r}{\underline{\mu}}\nonumber\\
\implies & T_j^{n,r}(t) \geq \frac{1}{\mu_j}\left(Q_j^r(\tau_n^r) - q_j^{*,r}(\tau_n^r) + \la_j^r t \right)- \frac{2\ep_3 r+1}{\underline{\mu}}. \label{eq_busy_p_3}
\end{align}

Let us fix an arbitrary $j\in\mJ^H$ and consider the following inequality:
\begin{equation}\label{eq_busy_p_4}
T_j^{n,r}(L^r(\tau_n^r)) > \frac{1}{\mu_j}\left(Q_j^r(\tau_n^r) -q_j^{*,r}(\tau_n^r) + \la_j^r L^r(\tau_n^r) \right) + \frac{4J\ep_3 r}{\underline{\mu}}.
\end{equation}
We will show that the inequality in \eqref{eq_busy_p_4} cannot hold when $r$ is sufficiently large by contradiction. Suppose that \eqref{eq_busy_p_4} holds. Then, it must be the case that $Q_j^r(\tau_n^r+L^r(\tau_n^r)) < q_j^{*,r}(\tau_n^r)$ by \eqref{eq_busy_p_2}. Let
\begin{equation*}
\eta^{n,r}_1:=\begin{cases}
\sup\left\{t\in[0,L^r(\tau_n^r)]: Q_j^r(\tau_n^r+t) \geq q_j^{*,r}(\tau_n^r) \right\},&\mbox{if such a $t$ exists,}\\
0,&\mbox{otherwise}.
\end{cases}
\end{equation*}
On the one hand, if $\eta^{n,r}_1=0$, then $T_j^{n,r}(\eta^{n,r}_1)=0$. On the other hand, if $\eta^{n,r}_1>0$, then by \eqref{eq_busy_p_2},
\begin{equation*}
T_j^{n,r}(\eta^{n,r}_1) \leq \frac{1}{\mu_j}\left(Q_j^r(\tau_n^r) -q_j^{*,r}(\tau_n^r) + \la_j^r \eta^{n,r}_1 \right) + \frac{2\ep_3 r}{\underline{\mu}}.
\end{equation*}
Therefore,
\begin{equation}\label{eq_busy_p_5}
T_j^{n,r}(\eta^{n,r}_1) \leq \left(\frac{1}{\mu_j}\left(Q_j^r(\tau_n^r) -q_j^{*,r}(\tau_n^r) + \la_j^r \eta^{n,r}_1 \right) + \frac{2\ep_3 r}{\underline{\mu}}\right)^+.
\end{equation}
Because $T_j^{n,r}$ is continuous, nondecreasing, almost everywhere differentiable under the proposed policy (specifically, $\dr T_j^{n,r}(t)\in\{0,1\}$ for almost every $t\in[0,r^2T]$), and by \eqref{eq_busy_p_4} and \eqref{eq_busy_p_5}, there exists an $\eta^{n,r}_2\in(\eta^{n,r}_1,L^r(\tau_n^r))$ such that 
\begin{equation}\label{eq_busy_p_6}
T_j^{n,r}(\eta^{n,r}_2) > \frac{1}{\mu_j}\left(Q_j^r(\tau_n^r) -q_j^{*,r}(\tau_n^r) + \la_j^r L^r(\tau_n^r) \right) + \frac{3J\ep_3 r}{\underline{\mu}}\quad\text{and}\quad \dr T_j^{n,r}(\eta^{n,r}_2) >0.
\end{equation}
By \eqref{eq_busy_p_2} and \eqref{eq_busy_p_6}, it must be the case that $Q_j^r(\tau_n^r+\eta^{n,r}_2) < q_j^{*,r}(\tau_n^r)$. Although the number of type $j$ jobs is less than the desired level $q_j^{*,r}(\tau_n^r)$ at time $\tau_n^r+\eta^{n,r}_2$, those jobs are still being processed at time $\eta^{n,r}_2$ by \eqref{eq_busy_p_6}. 

We make a couple of observations at this point. First, there cannot exist an $l\in\mJ^H$ such that $\mI_j\subset\mI_l$ and $Q_l^r(\tau_n^r+\eta^{n,r}_2) > \lceil q_l^{*,r}(\tau_n^r) \rceil$. Otherwise, the proposed policy would not process type $j$ jobs at time $\tau_n^r+\eta^{n,r}_2$. This is because the maximum $p_1(x(\tau_n^r+\eta^{n,r}_2),\tau_n^r)$ value under which type $l$ jobs are processed at time $\tau_n^r+\eta^{n,r}_2$ is at least $|\mI_j\cap \mI^H|$ greater than the same value under which type $j$ jobs are processed at that time (see \eqref{eq_index}) by Assumption \ref{a_h} and the proposed policy implements a resource allocation vector with the maximum $p_1(x(\tau_n^r+\eta^{n,r}_2),\tau_n^r)$ value at time $\tau_n^r+\eta^{n,r}_2$.

Second, suppose that there does not exist an $l\in\mJ^H$ such that $\mI_l\subset\mI_j$, $\mI_l\neq\mI_j$, and $Q_l^r(\tau_n^r+\eta^{n,r}_2) > \lceil q_l^{*,r}(\tau_n^r) \rceil$. By Assumption \ref{a_h}, we must have $Q_m^r(\tau_n^r+\eta^{n,r}_2) \leq \lceil q_m^{*,r}(\tau_n^r) \rceil$ for all $m\in\mJ_i$ and $i\in\mI_j\cap\mI^H$.

Third, suppose that there exist an $l\in\mJ^H$ such that $\mI_l\subset\mI_j$, $\mI_l\neq\mI_j$, and $Q_l^r(\tau_n^r+\eta^{n,r}_2) > \lceil q_l^{*,r}(\tau_n^r) \rceil$. Because type $j$ jobs are being processed instead of type $l$ jobs at time $\tau_n^r+\eta^{n,r}_2$ under the proposed policy, it must be the case that under any resource allocation vector that processes type $l$ jobs at time $\tau_n^r+\eta^{n,r}_2$, some resources in the set $\mI_j\cap\mI^H$ idle and because the proposed policy utilizes the resources in heavy traffic in a work-conserving fashion, those resource allocation vectors are not implemented at time $\tau_n^r+\eta^{n,r}_2$. Therefore, by Assumption \ref{a_h}, there must exist an $i\in\mI_j\cap\mI^H$ such that $Q_m^r(\tau_n^r+\eta^{n,r}_2) \leq \lceil q_m^{*,r}(\tau_n^r) \rceil$ for all $m\in\mJ_i$.

Consequently, there exists an $i\in\mI_j\cap\mI^H$ such that $Q_l^r(\tau_n^r+\eta^{n,r}_2) \leq \lceil q_l^{*,r}(\tau_n^r) \rceil$ for all $l\in\mJ_i$. In other words, there exists a resource in heavy traffic that process type $j$ jobs such that numbers of all job types that are processed by that resource are less than or equal to their associated desired levels at time $\tau_n^r+\eta^{n,r}_2$. Therefore, by \eqref{eq_busy_p_3},
\begin{equation}\label{eq_busy_p_7}
T_l^{n,r}(\eta^{n,r}_2) \geq \frac{1}{\mu_l}\left(Q_l^r(\tau_n^r) -q_l^{*,r}(\tau_n^r) + \la_l^r \eta^{n,r}_2 \right) - \frac{2\ep_3 r+1}{\underline{\mu}},\qquad\forall l\in\mJ_i.
\end{equation}
By \eqref{eq_d_21}, we have
\begin{align}
&I_i^{n,r}(\eta^{n,r}_2) + T_j^{n,r}(\eta^{n,r}_2) +\sum_{l\in\mJ_i\backslash\{j\}} T_l^{n,r}(\eta^{n,r}_2) = \eta^{n,r}_2 \nonumber\\
\implies & I_i^{n,r}(\eta^{n,r}_2) +  \frac{1}{\mu_j}\left(Q_j^r(\tau_n^r) -q_j^{*,r}(\tau_n^r) + \la_j^r L^r(\tau_n^r) \right) + \frac{3J\ep_3 r}{\underline{\mu}}\nonumber \\
&\hspace{3cm}+ \sum_{l\in\mJ_i\backslash\{j\}}\left(\frac{1}{\mu_l}\left(Q_l^r(\tau_n^r) -q_l^{*,r}(\tau_n^r) + \la_l^r \eta^{n,r}_2 \right) - \frac{2\ep_3 r+1}{\underline{\mu}}\right) < \eta^{n,r}_2 \label{eq_busy_p_8}\\
\implies & \frac{1}{\mu_j}\left(Q_j^r(\tau_n^r) -q_j^{*,r}(\tau_n^r) + \la_j^r \eta^{n,r}_2 \right) + \frac{3J\ep_3 r}{\underline{\mu}} \nonumber\\
&\hspace{3.5cm}+ \sum_{l\in\mJ_i\backslash\{j\}}\frac{1}{\mu_l}\left(Q_l^r(\tau_n^r) -q_l^{*,r}(\tau_n^r) + \la_l^r \eta^{n,r}_2 \right) - \frac{2J\ep_3 r+J}{\underline{\mu}} < \eta^{n,r}_2\label{eq_busy_p_80}\\
\implies &  \sum_{l\in\mJ_i}\frac{1}{\mu_l}\left(Q_l^r(\tau_n^r) -q_l^{*,r}(\tau_n^r) + \la_l^r \eta^{n,r}_2 \right) + \frac{J(\ep_3 r-1)}{\underline{\mu}} < \eta^{n,r}_2\nonumber\\
\implies &  \sum_{l\in\mJ_i}\frac{Q_l^r(\tau_n^r) -q_l^{*,r}(\tau_n^r) }{\mu_l} +\eta^{n,r}_2\rho_i^r + \frac{J(\ep_3 r-1)}{\underline{\mu}} < \eta^{n,r}_2\nonumber\\
\implies &  \frac{J(\ep_3 r-1)}{\underline{\mu}} < \eta^{n,r}_2 \left(1-\rho_i^r \right)\label{eq_busy_p_9}\\
\implies &  \frac{J(\ep_3 r-1)}{\underline{\mu}} < \ep_2 r \left|1- \rho_i^r \right|,\label{eq_busy_p_10}
\end{align}
where \eqref{eq_busy_p_8} follows from \eqref{eq_busy_p_6} and \eqref{eq_busy_p_7}, \eqref{eq_busy_p_80} is by the fact that $\eta^{n,r}_2 \leq L^r(\tau_n^r)$, \eqref{eq_busy_p_9} follows from \eqref{eq_we}, and \eqref{eq_busy_p_10} follows from the fact that $\eta^{n,r}_2\leq L^r(\tau_n^r)\leq \ep_2r$. If $r \geq r_1$ (see \eqref{eq_rdef_2}), \eqref{eq_busy_p_10} implies that
\begin{equation}\label{eq_busy_p_11}
\frac{J(\ep_3 r-1)}{\underline{\mu}} < \ep_2 \max_{i\in\mI^H}(1+|\theta_i|).
\end{equation}
However, if $r\geq r_{12} \geq r_1$ (see \eqref{eq_rdef_3}), then \eqref{eq_busy_p_11} cannot hold. Hence, we have a contradiction and so \eqref{eq_busy_p_4} cannot hold if $r\geq r_{12}$. Because  \eqref{eq_busy_p_4} is incorrect and $j\in\mJ^H$ is arbitrarily chosen, if $r\geq r_{12}$, we have
\begin{equation}\label{eq_busy_p_12}
T_j^{n,r}(L^r(\tau_n^r)) \leq \frac{1}{\mu_j}\left(Q_j^r(\tau_n^r) -q_j^{*,r}(\tau_n^r) + \la_j^r L^r(\tau_n^r) \right) + \frac{4J\ep_3 r}{\underline{\mu}},\qquad\forall j\in\mJ^H.
\end{equation}

The inequality \eqref{eq_busy_p_12} provides an upper bound on the shifted cumulative busy time processes. The following lemma, whose proof is presented in Section \ref{l_length_proof}, will be very useful to derive a lower bound on the shifted cumulative busy time processes.

\begin{lemma}\label{l_length}
In Step 3 and in the set $\left\{\mE_n^r,\; \kappa_n^r \leq L^r(\tau_n^r),\;\tau_n^r\leq r^2T,\;\mA_n^{(1),r}(\ep_2),\;\mA_{n-1}^r(\ep_2)\right\}$, there exists an $r_{13}\in\N_+$ independent of $n$ such that if $r\geq r_{13}$, we have
\begin{equation*}
\left| L^r(\tau_n^r)-\sum_{j\in\mJ_i} \frac{Q_j^r(\tau_n^r) -q_j^{*,r}(\tau_n^r) + \la_j^r L^r(\tau_n^r)}{\mu_j}\right| \leq \frac{J\ep_3r}{\underline{\mu}},\qquad\forall i\in\mI^H.
\end{equation*}
\end{lemma}

Let us consider arbitrary $j\in\mJ^H$, $i\in\mI_j\cap\mI^H$, and $r\geq r_{12}\vee r_{13}$. By \eqref{eq_d_21}, we have
\begin{align}
&I_i^{n,r}(L^r(\tau_n^r))+ T_j^{n,r}(L^r(\tau_n^r)) +\sum_{l\in\mJ_i\backslash\{j\}} T_l^{n,r}(L^r(\tau_n^r)) = L^r(\tau_n^r)\nonumber\\
\implies& T_j^{n,r}(L^r(\tau_n^r)) + \frac{3J\ep_3r }{\underline{\mu}}+\sum_{l\in\mJ_i\backslash\{j\}} \left(\frac{Q_l^r(\tau_n^r) -q_l^{*,r}(\tau_n^r) + \la_l^r L^r(\tau_n^r) }{\mu_l} + \frac{4J\ep_3 r}{\underline{\mu}}\right) \geq L^r(\tau_n^r)\label{eq_busy_p_13}\\
\implies& T_j^{n,r}(L^r(\tau_n^r))\geq L^r(\tau_n^r)-\sum_{l\in\mJ_i\backslash\{j\}} \frac{Q_l^r(\tau_n^r) -q_l^{*,r}(\tau_n^r) + \la_l^r L^r(\tau_n^r)}{\mu_l} -\frac{4J(J-1)\ep_3 r}{\underline{\mu}}  - \frac{3J\ep_3r}{\underline{\mu}}\nonumber\\
\implies& T_j^{n,r}(L^r(\tau_n^r))\geq \frac{Q_j^r(\tau_n^r) -q_j^{*,r}(\tau_n^r) + \la_j^r L^r(\tau_n^r)}{\mu_j}- \frac{J\ep_3r}{\underline{\mu}} -\frac{4J(J-1)\ep_3 r}{\underline{\mu}}  - \frac{3J\ep_3r}{\underline{\mu}}\label{eq_busy_p_14}\\
\implies& T_j^{n,r}(L^r(\tau_n^r))\geq  \frac{1}{\mu_j}\left(Q_j^r(\tau_n^r) -q_j^{*,r}(\tau_n^r) + \la_j^r L^r(\tau_n^r) \right)-\frac{4J^2\ep_3 r}{\underline{\mu}} ,\label{eq_busy_p_15}
\end{align}
where \eqref{eq_busy_p_13} is by \eqref{eq_idle2} and \eqref{eq_busy_p_12}, and \eqref{eq_busy_p_14} is by Lemma \ref{l_length}. Let $r_{11}:=r_{12}\vee r_{13}$. Then, the desired result follows by \eqref{eq_busy_p_12} and \eqref{eq_busy_p_15}.

\subsubsection{Proof of Lemma \ref{l_length}}\label{l_length_proof}

Recall that by definition (see \eqref{eq_review_3} and \eqref{eq_review_4})
\begin{equation}\label{eq_length_1}
L^r(\tau_n^r) = \max_{i\in\mI^H}\max_{j\in\mJ_i^{\leq,r}(\tau_n^r)} \frac{q_j^{*,r}(\tau_n^r) - Q_j^r(\tau_n^r)}{\la_j^r}.
\end{equation}
Without loss of generality, let $k\in\mJ^H$ be such that
\begin{equation}\label{eq_length_2}
L^r(\tau_n^r) = \frac{q_k^{*,r}(\tau_n^r) - Q_k^r(\tau_n^r)}{\la_k^r}.
\end{equation}
By \eqref{eq_good_set_4}, in the set $\{\tau_n^r\leq r^2T\}\cap\mA_{n-1}^r(\ep_2)$, we have
\begin{equation}\label{eq_length_3}
L^r(\tau_n^r) \leq \frac{C_8\ep_2 r}{\underline{\la}}.
\end{equation}

Similarly, let
\begin{equation}\label{eq_length_4}
\ot{L}^r(\tau_n^r):= \max_{i\in\mI^H}\max_{j\in\mJ_i^{\leq,r}(\tau_n^r)} \frac{q_j^{*,r}(\tau_n^r) - Q_j^r(\tau_n^r)}{\la_j}.
\end{equation}
The only difference between $L^r(\tau_n^r)$ and $\ot{L}^r(\tau_n^r)$ is that the latter one has the limiting arrival rate $\la_j$ instead of $\la_j^r$ in the denominator. Without loss of generality, let $m\in\mJ^H$ be such that
\begin{equation}\label{eq_length_5}
\ot{L}^r(\tau_n^r) = \frac{q_m^{*,r}(\tau_n^r) - Q_m^r(\tau_n^r)}{\la_m}.
\end{equation}
By \eqref{eq_review_2} and Lemma \ref{l_review}, we have
\begin{equation}\label{eq_length_6}
\ot{L}^r(\tau_n^r) = \sum_{j\in\mJ_i} \frac{Q_j^r(\tau_n^r) - q_j^{*,r}(\tau_n^r) + \la_j \ot{L}^r(\tau_n^r)}{\mu_j}.
\end{equation}

By \eqref{eq_length_1}, \eqref{eq_length_2}, \eqref{eq_length_4}, and \eqref{eq_length_5}, we have
\begin{align*}
&L^r(\tau_n^r) = \frac{q_k^{*,r}(\tau_n^r) - Q_k^r(\tau_n^r)}{\la_k^r} \geq  \frac{q_m^{*,r}(\tau_n^r) - Q_m^r(\tau_n^r)}{\la_m^r} = \frac{\la_m}{\la_m^r}\ot{L}^r(\tau_n^r),\\
&\ot{L}^r(\tau_n^r) = \frac{q_m^{*,r}(\tau_n^r) - Q_m^r(\tau_n^r)}{\la_m} \geq \frac{q_k^{*,r}(\tau_n^r) - Q_k^r(\tau_n^r)}{\la_k} = \frac{\la_k^r}{\la_k} L^r(\tau_n^r).
\end{align*}
Therefore,
\begin{equation*}
L^r(\tau_n^r) - \ot{L}^r(\tau_n^r) \leq \left(1- \frac{\la_k^r}{\la_k}\right) L^r(\tau_n^r)\qquad\text{and}\qquad \ot{L}^r(\tau_n^r) -L^r(\tau_n^r) \leq \left(\frac{\la_m^r}{\la_m}-1\right) L^r(\tau_n^r),
\end{equation*}
which implies
\begin{equation}\label{eq_length_7}
\left| L^r(\tau_n^r) - \ot{L}^r(\tau_n^r) \right| \leq \frac{\max_{j\in\mJ} \left|\la_j- \la_j^r\right|}{\underline{\la}} L^r(\tau_n^r) \leq  \frac{\max_{j\in\mJ} \left|\la_j- \la_j^r\right|C_8\ep_2r}{\underline{\la}^2},
\end{equation}
where the last inequality is by \eqref{eq_length_3}. Furthermore, for all $l\in\mJ^H$, we have
\begin{align}
\left| \la_l^rL^r(\tau_n^r) - \la_l\ot{L}^r(\tau_n^r) \right| &\leq \left| \la_l^rL^r(\tau_n^r) - \la_lL^r(\tau_n^r)+\la_lL^r(\tau_n^r)- \la_l\ot{L}^r(\tau_n^r) \right| \nonumber\\
&\leq  \left| \la_l^r - \la_l\right| L^r(\tau_n^r) + \la_l\left|L^r(\tau_n^r)- \ot{L}^r(\tau_n^r) \right|\nonumber\\
&\leq  \max_{j\in\mJ}\left|\la_j- \la_j^r\right| L^r(\tau_n^r) + \bar{\la} \frac{\max_{j\in\mJ} \left|\la_j- \la_j^r\right|}{\underline{\la}} L^r(\tau_n^r)\nonumber\\
&\leq \left( 1+\frac{\bar{\la} }{\underline{\la}}\right)\max_{j\in\mJ}\left|\la_j- \la_j^r\right|L^r(\tau_n^r),\nonumber\\
&\leq \left( 1+\frac{\bar{\la} }{\underline{\la}}\right)\max_{j\in\mJ}\left|\la_j- \la_j^r\right|\frac{C_8\ep_2 r}{\underline{\la}},\label{eq_length_8}
\end{align}
where the third inequality is by the first inequality in \eqref{eq_length_7}.

Let us fix an arbitrary $i\in\mI^H$. Then,
\begin{align}
&\left| L^r(\tau_n^r)-\sum_{j\in\mJ_i} \frac{Q_j^r(\tau_n^r) -q_j^{*,r}(\tau_n^r) + \la_j^r L^r(\tau_n^r)}{\mu_j}\right|\nonumber\\
&\hspace{0.5cm} =\left| L^r(\tau_n^r)-\sum_{j\in\mJ_i} \frac{Q_j^r(\tau_n^r) -q_j^{*,r}(\tau_n^r) + \la_j^r L^r(\tau_n^r)}{\mu_j} - \ot{L}^r(\tau_n^r) + \sum_{j\in\mJ_i} \frac{ Q_j^r(\tau_n^r)- q_j^{*,r}(\tau_n^r) +\la_j \ot{L}^r(\tau_n^r)}{\mu_j}\right|\nonumber\\
&\hspace{0.5cm} \leq \left| L^r(\tau_n^r) - \ot{L}^r(\tau_n^r)\right| + \left| \sum_{j\in\mJ_i} \frac{\la_j^r L^r(\tau_n^r) - \la_j \ot{L}^r(\tau_n^r)}{\mu_j} \right|\nonumber\\
&\hspace{0.5cm} \leq \frac{\max_{j\in\mJ} \left|\la_j- \la_j^r\right|C_8\ep_2r}{\underline{\la}^2} +  \frac{J}{\underline{\mu}}\left( 1+\frac{\bar{\la} }{\underline{\la}}\right)\max_{j\in\mJ}\left|\la_j- \la_j^r\right|\frac{C_8\ep_2 r}{\underline{\la}}\nonumber\\
&\hspace{0.5cm} \leq \left(\frac{1}{\underline{\la}^2} + \left( 1+\frac{\bar{\la} }{\underline{\la}}\right)\frac{J}{\underline{\la}\underline{\mu}} \right)\max_{j\in\mJ} \left|\la_j- \la_j^r\right|C_8\ep_2r,\label{eq_length_9}
\end{align}
where the equality is by \eqref{eq_length_6} and the second inequality is by \eqref{eq_length_7} and \eqref{eq_length_8}. 

By Assumption \ref{a_regime} Part 1, there exists $r_{13}\in\N_+$ independent of $n$ and $i$ such that if $r\geq r_{13}$,
\begin{equation*}
\max_{j\in\mJ} \left|\la_j- \la_j^r\right| \leq  \frac{J\ep_3}{\underline{\mu}\left(\frac{1}{\underline{\la}^2} + \left( 1+\frac{\bar{\la} }{\underline{\la}}\right)\frac{J}{\underline{\la}\underline{\mu}} \right)C_8\ep_2} = \frac{J}{4\underline{\mu}\left(\frac{1}{\underline{\la}^2} + \left( 1+\frac{\bar{\la} }{\underline{\la}}\right)\frac{J}{\underline{\la}\underline{\mu}} \right)\left(1+ 2J^2\frac{\bar{\mu}}{\underline{\mu}}\right)},
\end{equation*}
where the equality is by \eqref{eq_d_17_1}. Hence, if $r\geq r_{13}$, the term in \eqref{eq_length_9} is less than or equal to $J\ep_3r/\underline{\mu}$. Because $i\in\mI^H$ is arbitrarily chosen, we obtain the desired result.

\end{document}